\documentclass[a4paper,12pt,twoside]{article}
\usepackage{amsmath,amssymb,ifthen}
\usepackage[amsmath,thmmarks]{ntheorem}
\usepackage{paper}
\usepackage[all]{xy}
\usepackage{mathrsfs}
\usepackage{url}
\usepackage{stmaryrd}

\newcommand{\ad}{\mathrm{ad}}

\DeclareMathOperator{\Spa}{Spa}

\DeclareMathOperator{\Ind}{Ind}

\DeclareMathOperator{\PGL}{PGL}

\DeclareMathOperator{\diag}{diag}
\DeclareMathOperator{\SO}{SO}
\DeclareMathOperator{\Or}{O}

\DeclareMathOperator{\U}{U}

\DeclareMathOperator{\OGr}{OGr}
\DeclareMathOperator{\UGr}{UGr}
\DeclareMathOperator{\rad}{rad}
\DeclareMathOperator{\tr}{tr}
\DeclareMathOperator{\Ad}{Ad}
\DeclareMathOperator{\Res}{Res}
\DeclareMathOperator{\Fil}{Fil}
\DeclareMathOperator{\Stab}{Stab}

\SelectTips{cm}{11}
\title{On Siegel domains in $p$-adic period spaces}
\author{Yoichi Mieda}

\def\headertitle{On Siegel domains in $p$-adic period spaces}

\begin{document}
\maketitle

\begin{firstfootnote}
 Graduate School of Mathematical Sciences, The University of Tokyo, 3--8--1 Komaba, Meguro-ku, Tokyo, 153--8914, Japan

 E-mail address: \texttt{mieda@ms.u-tokyo.ac.jp}

\end{firstfootnote}

\begin{abstract}
 We consider an analogue of a Siegel domain in the $p$-adic period space of Rapoport--Zink. It is a quasi-compact open subset of the admissible locus
 whose translates under the group action cover the entire admissible locus. We prove that there exists a Siegel domain if and only if
 the admissible locus coincides with the weakly admissible locus, as conjectured by Chen--Fargues--Shen.
\end{abstract}

\section{Introduction}\label{sec:intro}
A classical Siegel domain is a ``small'' subset in a symmetric space whose translates under the action of the relevant arithmetic group
cover the entire space. In this paper, we consider an analogue of Siegel domains in $p$-adic period spaces.
Let $F$ be a $p$-adic field, $G$ a connected reductive group over $F$, $\mu\colon \mathbb{G}_m\to G_{\overline{F}}$ a minuscule cocharacter
defined over an algebraic closure $\overline{F}$ of $F$, and $b$ an element of $G(\widehat{F}^{\mathrm{ur}})$ whose Frobenius-conjugacy class
belongs to the set $B(G,\mu)$ defined in \cite[6.2]{MR1485921}.
Such a triple $(G,b,\mu)$ is called a local Shimura datum. 
We write $E$ for the field of definition of the conjugacy class of $\mu$.
In \cite{MR1393439}, Rapoport and Zink introduced a flag variety $\mathcal{F}(G,\mu)$ over $\widehat{E}^{\mathrm{ur}}$
and its rigid-analytic open subset $\mathcal{F}(G,b,\mu)^{\mathit{wa}}$, called the weakly admissible locus.
Roughly speaking, $\mathcal{F}(G,\mu)$ classifies Hodge filtrations of type $\mu$, and
such a filtration belongs to $\mathcal{F}(G,b,\mu)^{\mathit{wa}}$ if and only if, together with $b$, it gives rise to a weakly admissible $\varphi$-module
with $G$-structure.
For example, if $G=\GL_n$ and $\mu\colon \mathbb{G}_m\to \GL_n$; $z\mapsto (\underbrace{z,\ldots,z}_s,1,\ldots,1)$, then
$\mathcal{F}(G,\mu)$ classifies subspaces $\Fil\subset (\widehat{F}^{\mathrm{ur}})^n$ with $\dim \Fil=s$,
and $\Fil$ lies in $\mathcal{F}(G,b,\mu)^{\mathit{wa}}$
if and only if $((\widehat{F}^{\mathrm{ur}})^n,b,\Fil)$ gives a weakly admissible filtered $\varphi$-module.
Further, in \cite[1.37]{MR1393439}, Rapoport--Zink conjectured the existence of an open subset $\mathcal{F}(G,b,\mu)^a$
of $\mathcal{F}(G,b,\mu)^{\mathit{wa}}$ satisfying the following conditions:
\begin{itemize}
 \item for a finite extension $E'$ of $\widehat{E}^{\mathrm{ur}}$, we have $\mathcal{F}(G,b,\mu)^a(E')=\mathcal{F}(G,b,\mu)^\mathit{wa}(E')$.
 \item there exists a canonical $G$-valued local system over $\mathcal{F}(G,b,\mu)^a$ whose fiber at each $x\in \mathcal{F}(G,b,\mu)^\mathit{wa}(E')$
       is equal to the $G$-valued Galois representation attached to the weakly admissible filtered $\varphi$-module with $G$-structure corresponding to $x$.
\end{itemize}
Such an open subset $\mathcal{F}(G,b,\mu)^a$, which is nowadays called the admissible locus, was constructed by Faltings \cite{MR2658191} and
Hartl \cite{MR3091977} in some cases, and by Chen--Fargues--Shen \cite[Definition 3.1]{MR4325262} in general.
It is related to local Shimura varieties as follows.
Let $\mathcal{M}(G,b,\mu)_K$ denote the local Shimura variety with level $K\subset G(F)$
associated with the local Shimura datum $(G,b,\mu)$ (see \cite[Definition 24.1.3]{MR4446467}),
and $\pi_{\mathrm{dR}}\colon \mathcal{M}(G,b,\mu)_K\to \mathcal{F}(G,\mu)$
the de Rham period map (see \cite[p.~248]{MR4325262}; it is denoted by $\pi_{\mathrm{GM}}$ in \cite[\S 24.1]{MR4446467}).
Then $\pi_{\mathrm{dR}}$ is \'etale, and $\mathcal{F}(G,b,\mu)^a$ is equal to the image $\pi_{\mathrm{dR}}(\mathcal{M}(G,b,\mu)_K)$.
Hence the admissible locus is a genuine analogue of classical symmetric spaces, and the weakly admissible locus can be regarded as a linear-algebraic
approximation of the admissible locus.

We write $J_b(F)$ for the Frobenius-centralizer of $b$ in $G(\widehat{F}^{\mathrm{ur}})$. The group $J_b(F)$ naturally acts on $\mathcal{F}(G,\mu)$,
and this action preserves $\mathcal{F}(G,b,\mu)^{\mathit{wa}}$ and $\mathcal{F}(G,b,\mu)^a$.
Following \cite[\S 7]{MR4325262}, we view a quasi-compact open subset $U$ of $\mathcal{F}(G,b,\mu)^a$ satisfying 
$\mathcal{F}(G,b,\mu)^a=\bigcup_{g\in J_b(F)}gU$ as a $p$-adic analogue of a Siegel domain.
In the following, we assume that $b$ is basic.
In \cite{MR4325262}, Chen--Fargues--Shen proved that $\mathcal{F}(G,b,\mu)^a=\mathcal{F}(G,b,\mu)^{\mathit{wa}}$
if and only if $(G,\mu)$ is fully Hodge--Newton decomposable in the sense of \cite[Definition 3.1]{MR4060001}, and proposed the following conjecture:

\begin{conj}[{{\cite[Conjecture 7.2]{MR4325262}}}]\label{conj:CFS}
 The following are equivalent:
 \begin{enumerate}
  \item $\mathcal{F}(G,b,\mu)^a$ admits a Siegel domain; namely, there exists a quasi-compact open subset $U$ of $\mathcal{F}(G,b,\mu)^a$
	such that $\mathcal{F}(G,b,\mu)^a=\bigcup_{g\in J_b(F)}gU$.
  \item $\mathcal{F}(G,b,\mu)^a=\mathcal{F}(G,b,\mu)^{\mathit{wa}}$.
 \end{enumerate}
\end{conj}

The main theorem of this paper is the following:

\begin{thm}[{{Theorem \ref{thm:main}}}]
 Conjecture \ref{conj:CFS} is true.
\end{thm}

The strategy of our proof is as follows. First we explain the proof of the implication (i) $\Rightarrow$ (ii).
Suppose that $G=\GL_n$ and $\mu\colon \mathbb{G}_m\to \GL_n$; $z\mapsto (\underbrace{z,\ldots,z}_s,1,\ldots,1)$ with $s$ prime to $n$.
Then $J_b(F)$ is compact-mod-center, so $\mathcal{F}(G,b,\mu)^a=\bigcup_{g\in J_b(F)}gU$ is also quasi-compact. On the other hand, we have
$\mathcal{F}(G,b,\mu)^{\mathit{wa}}=\mathcal{F}(G,\mu)$. In particular, $\mathcal{F}(G,b,\mu)^{\mathit{wa}}$ is quasi-compact.
Further, $\mathcal{F}(G,b,\mu)^a$ and $\mathcal{F}(G,b,\mu)^{\mathit{wa}}$ have the same classical points.
In this situation, we can deduce the equality $\mathcal{F}(G,b,\mu)^a=\mathcal{F}(G,b,\mu)^{\mathit{wa}}$ by \cite[Corollary 4.3]{MR1207303}.
The key idea for extending this argument to the general case is to introduce the stable locus $\mathcal{F}(G,b,\mu)^{\mathrm{st}}$,
which is a rigid-analytic open subset of $\mathcal{F}(G,b,\mu)^{\mathit{wa}}$.
Recall that Totaro \cite{MR1388844} gave an interpretation of the weakly admissible locus in terms of Mumford's semistability \cite{MR1304906}.
The stable locus is obtained by replacing the semistability with stability in Totaro's description.
By using a classical result of Mumford on the stable locus, one can prove the inclusion $\mathcal{F}(G,b,\mu)^{\mathrm{st}}\subset \mathcal{F}(G,b,\mu)^a$
in a similar way as above. On the other hand, the complement $\mathcal{F}(G,b,\mu)^{\mathit{wa}}\smallsetminus \mathcal{F}(G,b,\mu)^{\mathrm{st}}$
(resp.\ $\mathcal{F}(G,b,\mu)^a\smallsetminus \mathcal{F}(G,b,\mu)^{\mathrm{st}}$) can be related to the weakly admissible (resp.\ admissible) loci
for proper Levi subgroups of $G$ by the Bialynicki-Birula contraction maps. Therefore, we can prove the inclusion
$\mathcal{F}(G,b,\mu)^{\mathit{wa}}\smallsetminus \mathcal{F}(G,b,\mu)^{\mathrm{st}}\subset \mathcal{F}(G,b,\mu)^a$ by induction on
the split semisimple rank of $G$.

The opposite implication (ii) $\Rightarrow$ (i) is proved by a case-by-case argument. By Chen--Fargues--Shen's theorem \cite[Theorem 6.1]{MR4325262},
the pair $(G,\mu)$ is fully Hodge--Newton decomposable. Such pairs are completely classified by \cite{MR4060001}.
Roughly speaking, we are reduced to considering 4 families: the Drinfeld case, the $\SO(2,n-2)$-case, the $\U(1,n-1)$-case, and the Hilbert--Blumenthal case.
In the first three cases we explicitly construct a quasi-compact open subset $U$ of $\mathcal{F}(G,b,\mu)^a=\mathcal{F}(G,b,\mu)^{\mathit{wa}}$, 
and check the condition $\mathcal{F}(G,b,\mu)^{\mathit{wa}}=\bigcup_{g\in J_b(F)}gU$ by elementary linear algebra.
In the Hilbert--Blumenthal case, we use an inductive argument inspired by the proof of (i) $\Rightarrow$ (ii).
Note that the existence of a Siegel domain in the Drinfeld case is well-known; for example, it can be constructed by using Deligne's formal model
of the Drinfeld upper half space. We also mention that Shen \cite{XuShen-RZ} discussed the existence of a Siegel domain in the $\U(1,n-1)$-case,
under the assumption that $F=\Q_p$, the quadratic extension giving the unitary group is unramified and $p>2$. 
However, the proof in \cite{XuShen-RZ} seems to have a gap; see Remark \ref{rem:Shen-gap}.

The author's original motivation for this research was to use Siegel domains to generalize the duality theory for
the $\ell$-adic \'etale cohomology of local Shimura varieties established in \cite{Fargues-Zel} and \cite{RZ-Zel}.
This remains a problem for future research.
It also turned out that the existence of a Siegel domain is very useful to establish the Steinness of local Shimura varieties associated with
fully Hodge--Newton decomposable local Shimura data, together with a technique developed in a recent preprint \cite{arXiv260922921}.
We will treat it in a forthcoming paper.

It seems natural to ask whether the weakly admissible locus always admits a Siegel domain.
After completing the proof in this paper, the author was informed by his student Kosuke Sugeno that
he had succeeded in proving the existence of a Siegel domain in $\mathcal{F}(G,b,\mu)^{\mathit{wa}}$ for every $(G,b,\mu)$.
His proof is uniform for all $(G,b,\mu)$, and does not rely on any classification.

The outline of this paper is as follows. In Section \ref{sec:period-space}, we briefly recall the definition and properties of
the weakly admissible locus and the admissible locus. In Section \ref{sec:reduction}, we give a proof of the implication (i) $\Rightarrow$ (ii)
of Conjecture \ref{conj:CFS}. Section \ref{sec:existence} is devoted to a proof of the converse (ii) $\Rightarrow$ (i).
We first give a classification of fully Hodge--Newton decomposable pairs by following \cite{MR4060001}.
After that, we construct a Siegel domain separately in each case.

\medbreak
\noindent{\bfseries Acknowledgment}\quad
A part of this work was carried out during the author's stay at Institut des Hautes \'Etudes Scientifiques.
The author is grateful to the institute for its hospitality.
This work was supported by JSPS KAKENHI Grant Number 25K00905.

\medbreak

\noindent{\bfseries Use of AI}\quad
The author used ChatGPT (GPT-6.2 Astra) as an auxiliary tool to identify typos, grammatical errors and possible gaps in proofs. 
He also used it to examine proofs in a published paper; see Remark \ref{rem:Shen-gap}.
The manuscript was written by the author.

\medbreak

\noindent{\bfseries Notation}\quad Let $p$ be a prime number and $F$ a finite extension of $\Q_p$.
We write $\mathcal{O}_F$ for the ring of integers of $F$, $\mathfrak{p}_F$ for the maximal ideal of $\mathcal{O}_F$,
and $k_F=\mathcal{O}_F/\mathfrak{p}_F$ for the residue field of $\mathcal{O}_F$.
Let $v_F$ denote the normalized valuation of $F$.
We fix a uniformizer $\varpi$ of $\mathcal{O}_F$. We also fix an algebraic closure $\overline{F}$ of $F$,
and write $\C_p$ for the completion of $\overline{F}$. Let $F^{\mathrm{ur}}$ denote the maximal unramified extension of $F$ contained in $\overline{F}$,
and $\breve{F}$ the completion of $F^{\mathrm{ur}}$. The residue field $\overline{k}_F$ of the ring of integers $\mathcal{O}_{\breve{F}}$ of $\breve{F}$
is an algebraic closure of $k_F$.

For an object $X$ over a ring $A$ and an $A$-algebra $B$, we write $X_B$ for the base change of $X$ to $B$.
For an $A$-module $M$, we write $\langle x_1,\ldots,x_r\rangle_A$ for the $A$-submodule of $M$
generated by $x_1,\ldots,x_r\in M$. 
For an algebraic variety $X$ over a non-archimedean field $k$, we write $X^{\ad}$ for the adic space over $\Spa(k,k^\circ)$ associated with $X$.
For an adic space $Y$ locally of finite type over $\Spa(k,k^\circ)$, we simply write $Y(k)$ for $Y(k,k^\circ)$.

For an algebraic group $G$ over $F$, let $X^*(G)$ (resp.\ $X_*(G)$) denote the set of characters $G_{\overline{F}}\to \mathbb{G}_m$
(resp.\ cocharacters $\mathbb{G}_m\to G_{\overline{F}}$) defined over $\overline{F}$.
We write $X^*_F(G)$ (resp.\ $X_{*,F}(G)$) for the subset of $X^*(G)$ (resp.\ $X_*(G)$)
consisting of characters (resp.\ cocharacters) defined over $F$.

\section{$p$-adic period spaces}\label{sec:period-space}
Let $G$ be a connected reductive group over $F$, $b$ an element of $G(\breve{F})$, and $\{\mu\}$ the conjugacy class of minuscule
cocharacters $\mathbb{G}_m\to G_{\overline{F}}$.
Let $E$ denote the reflex field of $\{\mu\}$, that is, the smallest subfield of $\overline{F}$ over which $\{\mu\}$ is defined.
In \cite[1.31]{MR1393439}, Rapoport--Zink considered an algebraic variety $\mathcal{F}(G,\{\mu\})$ over $E$ parametrizing
the quotient $\{\mu\}/{\sim}$ of $\{\mu\}$ by a suitable equivalence relation $\sim$.
If we fix an element $\mu\in \{\mu\}$, the base change $\mathcal{F}(G,\{\mu\})_{\overline{F}}$ is isomorphic
to the flag variety $G_{\overline{F}}/P(\mu)$,
where $P(\mu)=\{g\in G_{\overline{F}}\mid \text{$\lim_{t\to 0}\mu(t)g\mu(t)^{-1}$ exists}\}$ denotes the parabolic subgroup associated with the cocharacter $\mu$. Further, in \cite[Proposition 1.36]{MR1393439}, Rapoport--Zink constructed an open subset $\mathcal{F}(G,b,\{\mu\})^{\textit{wa}}$ of 
$\mathcal{F}(G,\{\mu\})_{\breve{E}}^{\ad}$, which is called the weakly admissible locus (to be precise, they use the framework of
Tate's rigid analytic spaces, but a similar construction can be done in the framework of adic spaces).
Roughly speaking, it parametrizes elements $\mu'\in \{\mu\}/{\sim}$ such that $(G,b,\mu')$ gives rise to a weakly admissible
filtered $\varphi$-module with $G$-structure.
It is also known that there exists a possibly smaller open subset $\mathcal{F}(G,b,\{\mu\})^a$, called the admissible locus,
equipped with a $G$-valued local system
whose fiber at each classical point $\mu'\in \mathcal{F}(G,b,\{\mu\})^a$ corresponds to $(G,b,\mu')$.
The existence of the admissible locus was conjectured by Rapoport--Zink (see \cite[1.37]{MR1393439}); 
the construction was given by Faltings \cite{MR2658191} and Hartl \cite{MR3091977} when $(G,b,\{\mu\})$ is of Hodge type,
and by Chen--Fargues--Shen \cite[Definition 3.1]{MR4325262} in general. Thanks to \cite{MR1779803}, for a finite extension $E'$ of $\breve{E}$,
we have $\mathcal{F}(G,b,\{\mu\})^{\textit{wa}}(E')=\mathcal{F}(G,b,\{\mu\})^a(E')$ (see \cite[Remarks A.5 (i)]{MR3861564-appendix}).

Let $J_b$ denote the Frobenius-centralizer of $b$ (see \cite[Proposition 1.12]{MR1393439}). Then, the group $J_b(F)$ acts on 
$\mathcal{F}(G,\{\mu\})_{\breve{E}}$, and the action preserves the open subsets $\mathcal{F}(G,b,\{\mu\})^{\textit{wa}}$, $\mathcal{F}(G,b,\{\mu\})^a$ of
$\mathcal{F}(G,\{\mu\})_{\breve{E}}^{\ad}$. Let $B(G)$ denote the set of Frobenius-conjugacy classes of $G(\breve{F})$.
The spaces $\mathcal{F}(G,\{\mu\})$, $\mathcal{F}(G,b,\{\mu\})^{\textit{wa}}$, $\mathcal{F}(G,b,\{\mu\})^a$
and the group $J_b$ essentially depend only on the image $[b]$ of $b$ in $B(G)$ (see \cite[Proposition 1.36]{MR1393439} for $\mathcal{F}(G,b,\{\mu\})^{\textit{wa}}$; the claim for $\mathcal{F}(G,b,\{\mu\})^a$ is clear from the definition).

In the following, we assume that $b$ is basic.
Based on a conjecture by Fargues and Rapoport, Chen--Fargues--Shen determined when $\mathcal{F}(G,b,\{\mu\})^{\textit{wa}}$
and $\mathcal{F}(G,b,\{\mu\})^a$ are equal:

\begin{thm}[{{\cite[Theorem 6.1]{MR4325262}; see also \cite[Theorem 1.1]{MR4060001}}}]\ \par
\noindent
 Assume that $b$ is basic and $[b]$ belongs to the subset $B(G,\{\mu\})$ of $B(G)$ (see \cite[6.2]{MR1485921} for the definition of
 $B(G,\{\mu\})$; recall that $B(G,\{\mu\})$ contains a unique basic element). Then, the following are equivalent:
 \begin{enumerate}
  \item $\mathcal{F}(G,b,\{\mu\})^{\textit{wa}}=\mathcal{F}(G,b,\{\mu\})^a$.
  \item the pair $(G,\{\mu\})$ is fully Hodge--Newton decomposable.
 \end{enumerate}
\end{thm}

See \cite[Definition 3.1]{MR4060001} for the full Hodge--Newton decomposability. Pairs $(G,\{\mu\})$ satisfying this condition are
completely classified in \cite[Theorem 3.5]{MR4060001}; see also Section \ref{subsec:HN-dec}.
Now we introduce a $p$-adic analogue of a Siegel domain.

\begin{defn}
 A Siegel domain in $\mathcal{F}(G,b,\{\mu\})^a$ is a quasi-compact open subset $U$ such that $\mathcal{F}(G,b,\{\mu\})^a=\bigcup_{g\in J_b(F)}gU$.
\end{defn}

The main theorem of this paper is the following:

\begin{thm}[{{\cite[Conjecture 7.2]{MR4325262}}}]\label{thm:main}
 Assume that $b$ is basic and $[b]$ belongs to $B(G,\{\mu\})$. Then, the following are equivalent:
 \begin{enumerate}
  \item $\mathcal{F}(G,b,\{\mu\})^a$ admits a Siegel domain.
  \item $\mathcal{F}(G,b,\{\mu\})^{\textit{wa}}=\mathcal{F}(G,b,\{\mu\})^a$.
 \end{enumerate}
\end{thm}

In the remaining part of this section, we collect some basic facts which are useful in the proof of Theorem \ref{thm:main}.

\begin{prop}\label{prop:passing-to-adjoint}
 Let $Z(G)$ denote the center of $G$, and put $G_{\ad}=G/Z(G)$. We write $b_{\ad}$ for the image of $b$ in $G_{\ad}(\breve{F})$, and
 $\{\mu_{\ad}\}$ for the conjugacy class of cocharacters $\mathbb{G}_m\to (G_{\ad})_{\overline{F}}$ induced by $\{\mu\}$.
 Let $E_{\mathrm{ad}}$ be the reflex field of $\{\mu_{\ad}\}$.
 Suppose that $[b]$ belongs to the subset $A(G,\{\mu\})$ of $B(G)$ (see \cite[Definition 1.3]{MR4612104} for the definition of the set $A(G,\{\mu\})$).
 \begin{enumerate}
  \item We have 
	\begin{align*}
	 \mathcal{F}(G,\{\mu\})&=\mathcal{F}(G_{\ad},\{\mu_{\ad}\})\otimes_{E_{\mathrm{ad}}}E,\\
	 \mathcal{F}(G,b,\{\mu\})^{\textit{wa}}&=\mathcal{F}(G_{\ad},b_{\ad},\{\mu_{\ad}\})^{\textit{wa}}\otimes_{\breve{E}_{\ad}}\breve{E},\\
	 \mathcal{F}(G,b,\{\mu\})^a&=\mathcal{F}(G_{\ad},b_{\ad},\{\mu_{\ad}\})^a\otimes_{\breve{E}_{\ad}}\breve{E}.
	\end{align*}
	
  \item The admissible locus $\mathcal{F}(G,b,\{\mu\})^a$ admits a Siegel domain if and only if so does $\mathcal{F}(G_{\ad},b_{\ad},\{\mu_{\ad}\})^a$.
 \end{enumerate}
\end{prop}

\begin{prf}
 (i) The first equality is trivial. For the second equality, see \cite[proof of Lemma A.11 (i)]{MR3861564-appendix} and \cite[Proposition 3.1]{MR3271247}.
 We prove the third equality. Let $(C,C^+)$ be a complete algebraically closed affinoid field over $(\C_p,\mathcal{O}_{\C_p})$
 (see \cite[Definition 1.1.5]{MR1734903} for the definition of an affinoid field).
 It suffices to prove the equality $\mathcal{F}(G,b,\{\mu\})^a(C,C^+)=\mathcal{F}(G_{\ad},b_{\ad},\{\mu_{\ad}\})^a(C,C^+)$.
 Since $\mathcal{F}(G,b,\{\mu\})^a$ and $\mathcal{F}(G_{\ad},b_{\ad},\{\mu_{\ad}\})^a$ are partially proper by \cite[Proposition 3.2 (1)]{MR4325262},
 we may assume $C^+=\mathcal{O}_C$.
 Clearly we have $\mathcal{F}(G,b,\{\mu\})^a(C,\mathcal{O}_C)\subset \mathcal{F}(G_{\ad},b_{\ad},\{\mu_{\ad}\})^a(C,\mathcal{O}_C)$.
 Therefore, it suffices to prove that every element $x\in \mathcal{F}(G_{\ad},b_{\ad},\{\mu_{\ad}\})^a(C,\mathcal{O}_C)$ belongs to
 $\mathcal{F}(G,b,\{\mu\})^a(C,\mathcal{O}_C)$.
 Let $X_{\mathrm{FF},C^\flat}$ denote the adic Fargues--Fontaine curve over $F$. We write $\mathcal{E}_b$ for the $G$-bundle corresponding to $b$,
 and $\mathcal{E}_{b,x}$ for the modification of $\mathcal{E}_b$ corresponding to $x$.
 We fix a quasi-split inner form $G^*$ of $G$ and an inner twist $\psi\colon G_{\overline{F}}\xrightarrow{\cong}G^*_{\overline{F}}$.
 Then we have the Newton map $\nu_{\mathcal{E}_{b,x}}\colon \mathbb{D}\to G^*$, which is defined over $F$.
 Since $x$ belongs to $\mathcal{F}(G_{\ad},b_{\ad},\{\mu_{\ad}\})^a$, $\nu_{\mathcal{E}_{b,x}}$ factors through $Z(G^*)$.
 We need to prove that $\nu_{\mathcal{E}_{b,x}}$ is trivial.
 To observe it, it suffices to show that $\chi\circ \nu_{\mathcal{E}_{b,x}}$ is trivial for every rational character $\chi\in X^*_F(G^*)$.
 We write $\chi'$ for the composite $\chi\circ\psi$, which gives an element of $X^*_F(G)$.
 By \cite[p.~231]{MR4325262}, we have $c_1(\mathcal{E}_{b,x})=\mu^\sharp-\kappa(b)$, where $\mu^\sharp$ denotes the image of $\mu$
 in $\pi_1(G)_{\Gal(\overline{F}/F)}$. Therefore, we have $c_1(\chi'_*\mathcal{E}_{b,x})=\langle \chi',\mu\rangle-v_F(\chi'(b))=\langle \chi,\psi\circ\mu\rangle-\langle \chi,\nu_b\rangle=\langle \chi,\mu^\diamond\rangle-\langle \chi,\nu_b\rangle=0$, where the last equality follows 
 from the assumption $[b]\in A(G,\{\mu\})$ (for the definition of $\mu^\diamond$, see \cite[\S 1.3]{MR4612104}). 
 This means that the $\mathbb{G}_m$-bundle
 $\chi'_*\mathcal{E}_{b,x}$ is trivial. Hence $\chi\circ \nu_{\mathcal{E}_{b,x}}=\nu_{\chi'_*\mathcal{E}_{b,x}}$ is also trivial, as desired.

 (ii) Note that $J_{b_{\mathrm{ad}}}\cong (J_b)_{\ad}=J_b/Z(G)$.
 A Siegel domain in $\mathcal{F}(G,b,\{\mu\})^a$ clearly gives a Siegel domain in $\mathcal{F}(G_{\ad},b_{\ad},\{\mu_{ad}\})^a$.
 To prove the converse, it suffices to show that the image of $J_b(F)\to (J_b)_{\ad}(F)$ has finite index in $(J_b)_{\ad}(F)$.
 This follows from the exact sequence $J_b(F)\to (J_b)_{\ad}(F)\to H^1(F,Z(G))$ and the finiteness of $H^1(F,Z(G))$ (see \cite[Theorem 6.14]{MR1278263}).
\end{prf}

\begin{prop}\label{prop:passing-to-inner-form}
 Assume that $G$ is adjoint, $b$ is basic and decent (see \cite[Definition 1.8]{MR1393439}; in particular, $b$ belongs to $G(F^{\mathrm{ur}})$ 
 by \cite[Corollary 1.9]{MR1393439}). Let $(J_b)_{F^{\mathrm{ur}}}\xrightarrow{\cong}G_{F^{\mathrm{ur}}}$ be
 an isomorphism under which the Frobenius action on $(J_b)_{F^{\mathrm{ur}}}$ corresponds to $g\mapsto b\sigma_G(g)b^{-1}$,
 where $\sigma_G$ denotes the Frobenius action on $G_{F^{\mathrm{ur}}}$. We regard $\{\mu\}$ as a conjugacy class of cocharacters 
 $\mathbb{G}_m\to(J_b)_{\overline{F}}$ by means of the isomorphism above.
 Then we have identifications $\mathcal{F}(G,b,\{\mu\})^{\textit{wa}}=\mathcal{F}(J_b,1,\{\mu\})^{\textit{wa}}$ and 
 $\mathcal{F}(G,b,\{\mu\})^a=\mathcal{F}(J_b,1,\{\mu\})^a$ compatible with $J_b(F)$-actions.
\end{prop}

\begin{prf}
 See \cite[Remarks 2.15 (3), 2.18 (3)]{MR4612104}.
\end{prf}

\begin{rem}
 Even if $b$ lies in $B(G,\{\mu\})$, $1$ does not necessarily belong to $B(J_b,\{\mu\})$. It belongs to $A(J_b,\{\mu\})$.
\end{rem}

\begin{prop}\label{prop:twin-tower-principle}
 Assume that $b$ is basic and $[b]$ belongs to $B(G,\{\mu\})$.
 Let $(J_b,b'',\{\mu^{-1}\})$ be the dual triple attached to $(G,b,\{\mu\})$ (see \cite[\S 5.1]{MR4325262}).
 Then, $\mathcal{F}(G,b,\{\mu\})^a$ admits a Siegel domain if and only if so does $\mathcal{F}(J_b,b'',\{\mu^{-1}\})^a$.
\end{prop}

\begin{prf}
 Let $\mathcal{M}(G,b,\{\mu\})_\infty$ be the local Shimura variety at infinite level associated with the triple $(G,b,\{\mu\})$.
 Let $\pi_{\mathrm{dR}}\colon \mathcal{M}(G,b,\{\mu\})_\infty\to \mathcal{F}(G,b,\{\mu\})^a$
 (resp.\ $\pi_{\mathrm{HT}}\colon \mathcal{M}(G,b,\{\mu\})_\infty\to \mathcal{F}(J_b,b'',\{\mu^{-1}\})^a$) be the de Rham
 (resp.\ Hodge--Tate) period map, which is known to be surjective (see \cite[p.~248]{MR4325262}).
 Suppose that $\mathcal{F}(G,b,\{\mu\})^a$ admits a Siegel domain. Then, by \cite[the comment after Conjecture 7.2]{MR4325262},
 we can find a quasi-compact open subset $V$ of $\mathcal{M}(G,b,\{\mu\})_\infty$ such that 
 $\mathcal{M}(G,b,\{\mu\})_\infty=\bigcup_{g\in J_b(F)\times G(F)}gV$. It is easy to see that $\pi_{\mathrm{HT}}(V)$ gives a Siegel domain in 
 $\mathcal{F}(J_b,b'',\{\mu^{-1}\})^a$. 
 We can also obtain the converse direction by replacing $(G,b,\{\mu\})$ with $(J_b,b'',\{\mu^{-1}\})$. 
\end{prf}

\begin{prop}\label{prop:base-change-to-C_p}
 If we define a Siegel domain in $\mathcal{F}(G,b,\{\mu\})^a_{\C_p}$ in the same way as that in $\mathcal{F}(G,b,\{\mu\})^a$,
 then the following are equivalent:
 \begin{enumerate}
  \item $\mathcal{F}(G,b,\{\mu\})^a$ admits a Siegel domain.
  \item $\mathcal{F}(G,b,\{\mu\})^a_{\C_p}$ admits a Siegel domain.
 \end{enumerate}
\end{prop}

\begin{prf}
 Clearly (i) implies (ii). We prove the converse. Let $U$ be a Siegel domain in $\mathcal{F}(G,b,\{\mu\})^a_{\C_p}$.
 By \cite[Remark 2.4.3 ii)]{MR1734903}, we can find a finite Galois extension $E'$ of $\breve{E}$ contained in $\C_p$
 and a quasi-compact open subset $U'$ of $\mathcal{F}(G,b,\{\mu\})^a_{E'}$
 whose inverse image under $\mathcal{F}(G,b,\{\mu\})^a_{\C_p}\to \mathcal{F}(G,b,\{\mu\})^a_{E'}$ is equal to
 $U$. The open subset $\bigcup_{\sigma\in \Gal(E'/\breve{E})}\sigma U'$ of $\mathcal{F}(G,b,\{\mu\})^a_{E'}$ is quasi-compact and
 descends to a quasi-compact open subset $U''$ of $\mathcal{F}(G,b,\{\mu\})^a$. Clearly $U''$ is a Siegel domain in $\mathcal{F}(G,b,\{\mu\})^a$.
\end{prf}

By Proposition \ref{prop:base-change-to-C_p}, to prove Theorem \ref{thm:main}, we can replace $\mathcal{F}(G,b,\{\mu\})^{\textit{wa}}$ and
$\mathcal{F}(G,b,\{\mu\})^a$ with $\mathcal{F}(G,b,\{\mu\})_{\C_p}^{\textit{wa}}$ and $\mathcal{F}(G,b,\{\mu\})_{\C_p}^a$, respectively.
In the following, we simply write $\mathcal{F}(G,\{\mu\})$, $\mathcal{F}(G,b,\{\mu\})^{\textit{wa}}$ and $\mathcal{F}(G,b,\{\mu\})^a$ for
the base changes $\mathcal{F}(G,\{\mu\})_{\C_p}$, $\mathcal{F}(G,b,\{\mu\})_{\C_p}^{\textit{wa}}$ and $\mathcal{F}(G,b,\{\mu\})_{\C_p}^a$,
respectively. The original spaces $\mathcal{F}(G,\{\mu\})$, $\mathcal{F}(G,b,\{\mu\})^{\textit{wa}}$ and $\mathcal{F}(G,b,\{\mu\})^a$ over $\breve{E}$
are referred to as $\mathcal{F}_0(G,\{\mu\})$, $\mathcal{F}_0(G,b,\{\mu\})^{\textit{wa}}$ and $\mathcal{F}_0(G,b,\{\mu\})^a$.

\begin{defn}
 We define a subset $\mathcal{F}(G,\{\mu\})^{\ad}(\cl)$ of $\mathcal{F}(G,\{\mu\})^{\ad}(\C_p)$
 by 
 \[
  \mathcal{F}(G,\{\mu\})^{\ad}(\cl)=\bigcup_{E'}\mathcal{F}_0(G,\{\mu\})^{\ad}(E'),
 \]
 where $E'$ runs through finite extensions of $\breve{E}$ contained in $\C_p$. We call an element of $\mathcal{F}(G,\{\mu\})^{\ad}(\cl)$
 a classical point of $\mathcal{F}(G,\{\mu\})^{\ad}$. For a locally closed adic subspace $Z$ of $\mathcal{F}(G,\{\mu\})^{\ad}$, we put
 $Z(\cl)=Z(\C_p)\cap \mathcal{F}(G,\{\mu\})^{\ad}(\cl)$.
\end{defn}

Note that we have $\mathcal{F}(G,b,\{\mu\})^a(\cl)=\mathcal{F}(G,b,\{\mu\})^{\mathit{wa}}(\cl)$.
The following simple lemma is useful in the proof of (i) $\Rightarrow$ (ii) of Theorem \ref{thm:main}.

\begin{lem}\label{lem:classical-points}
 \begin{enumerate}
  \item Let $U$ and $V$ be quasi-compact open subsets of $\mathcal{F}(G,\{\mu\})^{\ad}$. If $U(\cl)=V(\cl)$, then $U=V$.
  \item Let $E'$ be a finite extension of $\breve{E}$ contained in $\C_p$, and $Z'$ a locally closed adic subspace of $\mathcal{F}_0(G,\{\mu\})^{\ad}_{E'}$.
	If $Z'_{\C_p}(\cl)=\varnothing$, then $Z'_{\C_p}=\varnothing$.
 \end{enumerate}
\end{lem}

\begin{prf}
 (i) In the same way as in the proof of Proposition \ref{prop:base-change-to-C_p}, we can find a finite extension $E'$ of $\breve{E}$ contained in $\C_p$
 and quasi-compact open subsets $U'$, $V'$ of $\mathcal{F}_0(G,\{\mu\})^{\ad}_{E'}$ such that $U'_{\C_p}=U$ and $V'_{\C_p}=V$.
 By the assumption and \cite[Corollary 4.3]{MR1207303}, we have $U'=V'$. Therefore we conclude that $U=V$.

 (ii) Since $Z'$ is locally of finite type over $\Spa(E',\mathcal{O}_{E'})$, the claim is clear.
\end{prf}

\section{Reduction theory}\label{sec:reduction}
In this section, we prove the implication (i) $\Rightarrow$ (ii) of Theorem \ref{thm:main}.
By Proposition \ref{prop:passing-to-adjoint}, we may assume that $G$ is adjoint. By \cite[5.2]{MR809866} (see also \cite[1.11]{MR1393439}), 
we can replace $b$ by its Frobenius-conjugate so that $b$ is decent. Therefore, by Proposition \ref{prop:passing-to-inner-form},
it suffices to show the following theorem:

\begin{thm}\label{thm:A-version}
 Let $G$ be a connected reductive group over $F$ and $\{\mu\}$ a conjugacy class of minuscule cocharacters $\mathbb{G}_m\to G_{\overline{F}}$.
 If $\mathcal{F}(G,1,\{\mu\})^a$ admits a Siegel domain, then $\mathcal{F}(G,1,\{\mu\})^{\textit{wa}}=\mathcal{F}(G,1,\{\mu\})^a$ holds.
\end{thm}

We prove this theorem by induction on the split semisimple rank of $G$, so assume that Theorem \ref{thm:A-version} holds
for any connected reductive group whose split semisimple rank is smaller than that of $G$.
If $1\notin A(G,\{\mu\})$, then $\mathcal{F}(G,1,\{\mu\})^{\textit{wa}}=\varnothing$ by \cite[Proposition 3.1]{MR3271247},
so the equality $\mathcal{F}(G,1,\{\mu\})^{\textit{wa}}=\mathcal{F}(G,1,\{\mu\})^a$ obviously holds.
In the following, we consider the case $1\in A(G,\{\mu\})$.
Then, by Proposition \ref{prop:passing-to-adjoint}, we may assume that $G$ is adjoint.

The key to our proof is to consider the stable locus inside $\mathcal{F}(G,1,\{\mu\})^{\textit{wa}}$. To define it, let us fix a maximal split torus $S$
of $G$ and a maximal torus $T$ of $G_{\overline{F}}$ containing $S_{\overline{F}}$. We fix a cocharacter $\mu\in \{\mu\}$ whose image lies in $T$;
so we may regard $\mu$ as an element of $X_*(T)$. Recall that this choice gives an identification $\mathcal{F}(G,\{\mu\})=G_{\C_p}/P(\mu)_{\C_p}$.
In the same way as in \cite[\S 2]{MR2198221}, we can define an element $\mu^*$ of $X^*(P(\mu))\otimes_\Z\Q$
(it depends on the choice of an invariant inner product of $G$).
We take an integer $d>0$ such that $d\mu^*$ lies in $X^*(P(\mu))$, and write $\mathcal{L}$
for the $G_{\C_p}$-equivariant ample line bundle on $\mathcal{F}(G,\{\mu\})$ corresponding to $-d\mu^*$.
For an $F$-rational cocharacter $\lambda\in X_{*,F}(G)$ and $x\in \mathcal{F}(G,\{\mu\})$,
we can define an integer $\mu^{\mathcal{L}}(x,\lambda)$ (see \cite[Definition 2.2]{MR1304906}).
For an element $\lambda\in X_{*,F}(S)\otimes_\Z\Q$, we take an integer $n_\lambda>0$ such that
$n_\lambda\lambda\in X_{*,F}(S)$ and define $\mu^{\mathcal{L}}(x,\lambda)=n_\lambda^{-1}\mu^{\mathcal{L}}(x,n_\lambda\lambda)\in \Q$;
it is independent of the choice of $n_\lambda$.

\begin{defn}
 We fix a minimal parabolic subgroup $P_0$ containing $S$, and write $\Delta\subset X_F^*(S)$ for the set of relative simple roots with respect to $P_0$.
 Let $\{\omega_\alpha\}_{\alpha\in \Delta}\subset X_{*,F}(S)\otimes_\Z\Q$ be the dual basis to $\Delta$. We take an integer $n_\alpha>0$ such that
 $n_\alpha\omega_\alpha\in X_{*,F}(S)$, and put $P_\alpha=P(n_\alpha\omega_\alpha)$. For $\alpha\in \Delta$, we put
 \begin{itemize}
  \item $Z_{\alpha,<0}=\{x\in\mathcal{F}(G,\{\mu\})\mid \mu^{\mathcal{L}}(x,\omega_\alpha)<0\}$, and
  \item $Z_{\alpha,\le 0}=\{x\in\mathcal{F}(G,\{\mu\})\mid \mu^{\mathcal{L}}(x,\omega_\alpha)\le 0\}$.
 \end{itemize}
 These are closed subvarieties of $\mathcal{F}(G,\{\mu\})$. Further, we set
 \begin{itemize}
  \item $\mathcal{F}(G,1,\{\mu\})^{\mathrm{ss}}=\mathcal{F}(G,\{\mu\})^{\ad}\smallsetminus \bigcup_{\alpha\in\Delta,g\in G(F)}gZ_{\alpha,<0}^{\ad}$, and
  \item $\mathcal{F}(G,1,\{\mu\})^{\mathrm{st}}=\mathcal{F}(G,\{\mu\})^{\ad}\smallsetminus \bigcup_{\alpha\in\Delta,g\in G(F)}gZ_{\alpha,\le 0}^{\ad}$.
 \end{itemize}
 These are $G(F)$-invariant subsets of $\mathcal{F}(G,\{\mu\})^{\ad}$.
\end{defn}

\begin{prop}\label{prop:ss=wa}
 \begin{enumerate}
  \item For $\alpha\in \Delta$, $Z_{\alpha,\le 0}$ and $Z_{\alpha,<0}$ are $P_{\alpha,\C_p}$-invariant.
  \item The subsets $\mathcal{F}(G,1,\{\mu\})^{\mathrm{ss}}$ and $\mathcal{F}(G,1,\{\mu\})^{\mathrm{st}}$ are open in $\mathcal{F}(G,\{\mu\})^{\ad}$;
	in particular, they are naturally regarded as open adic subspaces of $\mathcal{F}(G,\{\mu\})^{\ad}$.
  \item We have $\mathcal{F}(G,1,\{\mu\})^{\mathrm{ss}}=\mathcal{F}(G,1,\{\mu\})^{\mathit{wa}}$.
 \end{enumerate}
\end{prop}

\begin{prf}
 (i) By \cite[Definition 2.2 (i), Proposition 2.7]{MR1304906}, for $x\in\mathcal{F}(G,\{\mu\})$ and $g\in P_{\alpha,\C_p}$
 we have $\mu^{\mathcal{L}}(gx,\omega_\alpha)=\mu^{\mathcal{L}}(x,\Ad(g)^{-1}\circ\omega_\alpha)=\mu^{\mathcal{L}}(x,\omega_\alpha)$.
 The claim immediately follows from it.

 (ii) Consider a morphism $f\colon Z_{\alpha,<0}\times G_{\C_p}\to \mathcal{F}(G,\{\mu\})\times G_{\C_p}$; $(x,g)\mapsto (gx,g)$,
 which is a closed immersion. By (i), we can define an action of $P_{\alpha,\C_p}$ on $Z_{\alpha,<0}\times G_{\C_p}$ by $h(x,g)=(hx,gh^{-1})$.
 By taking the quotient by $P_{\alpha,\C_p}$, we have a closed immersion 
 $\overline{f}\colon (Z_{\alpha,<0}\times G_{\C_p})/P_{\alpha,\C_p}\to \mathcal{F}(G,\{\mu\})\times G_{\C_p}/P_{\alpha,\C_p}$.
 By applying \cite[1.32]{MR1393439} to the image of $\overline{f}$, we conclude that $\mathcal{F}(G,1,\{\mu\})^{\mathrm{ss}}$ is open 
 in $\mathcal{F}(G,\{\mu\})^{\ad}$ (see also \cite[proof of Lemma 3.2]{MR2198221} for an argument in the context of adic spaces).
 The openness of $\mathcal{F}(G,1,\{\mu\})^{\mathrm{st}}$ can be proven exactly in the same way.

 (iii) By the construction in \cite[Proposition 1.36]{MR1393439}, it turns out that there exists a closed subvariety $Z$ of $\mathcal{F}(G,\{\mu\})$
 defined over a finite extension of $\breve{E}$
 such that $\mathcal{F}(G,1,\{\mu\})^{\textit{wa}}=\mathcal{F}(G,\{\mu\})^{\ad}\smallsetminus\bigcup_{g\in G(F)}gZ^{\ad}$
 (see also \cite[p.~230]{MR4325262}). By \cite[Corollary 2.4]{MR2198221},
 we have 
 \[
  \bigcup_{g\in G(F)}gZ^{\ad}(\C_p)=\bigcup_{g\in G(F)}gZ(\C_p)=\bigcup_{\alpha\in \Delta,g\in G(F)}gZ_{\alpha,<0}(\C_p)=\bigcup_{\alpha\in \Delta,g\in G(F)}gZ_{\alpha,<0}^{\ad}(\C_p)
 \]
 (although $G$ is assumed to be quasi-split in \cite{MR2198221}, the proof of \cite[Corollary 2.4]{MR2198221} is valid for an arbitrary $G$).
 Therefore, $Z^{\ad}\cap\mathcal{F}(G,1,\{\mu\})^{\mathrm{ss}}$, which is an open adic subspace of $Z^{\ad}$ by (ii), has no $\C_p$-valued point.
 Hence we have $Z^{\ad}\cap\mathcal{F}(G,1,\{\mu\})^{\mathrm{ss}}=\varnothing$. As $\mathcal{F}(G,1,\{\mu\})^{\mathrm{ss}}$ is $G(F)$-invariant,
 we have $(\bigcup_{g\in G(F)}gZ^{\ad})\cap\mathcal{F}(G,1,\{\mu\})^{\mathrm{ss}}=\varnothing$, namely,
 $\mathcal{F}(G,1,\{\mu\})^{\mathrm{ss}}\subset \mathcal{F}(G,1,\{\mu\})^{\mathit{wa}}$.
 For the reverse inclusion, note that $Z_{\alpha,<0}^{\ad}\cap\mathcal{F}(G,1,\{\mu\})^{\mathit{wa}}$,
 which is an open adic subspace of $Z_{\alpha,<0}^{\ad}$, has no $\C_p$-valued point for every $\alpha\in \Delta$.
 By the same argument as above, we conclude that
 $(\bigcup_{\alpha\in\Delta, g\in G(F)}gZ_{\alpha,<0}^{\ad})\cap\mathcal{F}(G,1,\{\mu\})^{\mathit{wa}}=\varnothing$,
 that is, $\mathcal{F}(G,1,\{\mu\})^{\mathit{wa}}\subset \mathcal{F}(G,1,\{\mu\})^{\mathrm{ss}}$.
\end{prf}

\begin{rem}\label{rem:ss=wa-without-adjointness}
 The proof of Proposition \ref{prop:ss=wa} works for $G$ which is not necessarily adjoint, if we replace $S$ with a maximal split torus of
 $G_{\mathrm{der}}$. In Section \ref{sec:existence}, we use Proposition \ref{prop:ss=wa} (iii) in this form.
\end{rem}

In the remaining part of this section, we always use the invariant inner product of $G$ determined by the adjoint representation of $G$ to define $\mu^*$.
To show $\mathcal{F}(G,1,\{\mu\})^a=\mathcal{F}(G,1,\{\mu\})^{\mathit{wa}}$ under the assumption in Theorem \ref{thm:A-version}, it suffices to prove
$\mathcal{F}(G,1,\{\mu\})^{\mathrm{st}}\subset \mathcal{F}(G,1,\{\mu\})^a$ and
$\mathcal{F}(G,1,\{\mu\})^{\mathit{wa}}\smallsetminus \mathcal{F}(G,1,\{\mu\})^{\mathrm{st}}\subset \mathcal{F}(G,1,\{\mu\})^a$.
First we treat the former inclusion.

\begin{defn}
 Let $\mathcal{F}(G,\{\mu\})^{\mathrm{alg}\text{-}\mathrm{ss}}$ (resp.\ $\mathcal{F}(G,\{\mu\})^{\mathrm{alg}\text{-}\mathrm{st}}$) be
 the semistable (resp.\ stable) locus of $\mathcal{F}(G,\{\mu\})$ with respect to the action of $S_{\C_p}$.
 These are open subsets of $\mathcal{F}(G,\{\mu\})$. By \cite[Theorem 2.1]{MR1304906}, we have
 \begin{itemize}
  \item $\mathcal{F}(G,\{\mu\})^{\mathrm{alg}\text{-}\mathrm{ss}}=\{x\in\mathcal{F}(G,\{\mu\})\mid \mu^{\mathcal{L}}(x,\lambda)\ge 0\ (\lambda\in X_*(S_{\C_p})\smallsetminus \{0\})\}$, and
  \item $\mathcal{F}(G,\{\mu\})^{\mathrm{alg}\text{-}\mathrm{st}}=\{x\in\mathcal{F}(G,\{\mu\})\mid \mu^{\mathcal{L}}(x,\lambda)> 0\ (\lambda\in X_*(S_{\C_p})\smallsetminus \{0\})\}$.
 \end{itemize}
 Here we use the definition of a stable point in \cite[\S 8.1]{MR2004511}; it is referred to as a properly stable point in \cite[Definition 1.8]{MR1304906}.
\end{defn}

\begin{prop}\label{prop:alg-st}
 We have $\mathcal{F}(G,1,\{\mu\})^{\mathrm{ss}}\subset \mathcal{F}(G,\{\mu\})^{\mathrm{alg}\text{-}\mathrm{ss},\ad}$ and 
 $\mathcal{F}(G,1,\{\mu\})^{\mathrm{st}}\subset \mathcal{F}(G,\{\mu\})^{\mathrm{alg}\text{-}\mathrm{st},\ad}$.
\end{prop}

\begin{prf}
 Let $W(G,S)=N_G(S)(F)/Z_G(S)(F)$ denote the relative Weyl group.
 We fix a lift $\dot{w}\in N_G(S)(F)$ of each element $w\in W(G,S)$.
 Let us prove that $\{x\in \mathcal{F}(G,\{\mu\})\mid \mu^{\mathcal{L}}(x,\lambda)<0\}\subset \bigcup_{w\in W(G,S),\alpha\in\Delta}\dot{w}^{-1}Z_{\alpha,<0}$ for every $\lambda\in X_*(S_{\C_p})\smallsetminus\{0\}$. 
 Since $S$ is split, $\lambda\in X_*(S_{\C_p})\smallsetminus \{0\}$ is an $F$-rational cocharacter $\mathbb{G}_m\to S$.
 We take $w\in W(G,S)$ such that $w\lambda$ is dominant with respect to $P_0$;
 then there exists $(c_\alpha)_{\alpha\in\Delta}\in (\prod_{\alpha\in\Delta}\Q_{\ge 0})\smallsetminus \{0\}$ such that
 $w\lambda=\sum_{\alpha\in\Delta}c_\alpha\omega_\alpha$ (here we use the assumption that $G$ is adjoint, hence semisimple).
 By \cite[Definition 2.2 (i)]{MR1304906} and \cite[Proposition 2.3]{MR2198221}, we have
 \[
  \mu^{\mathcal{L}}(x,\lambda)=\mu^{\mathcal{L}}(\dot{w}x,w\lambda)=\mu^{\mathcal{L}}\Bigl(\dot{w}x,\sum_{\alpha\in\Delta}c_\alpha\omega_\alpha\Bigr)=\sum_{\alpha\in\Delta}c_\alpha\mu^{\mathcal{L}}(\dot{w}x,\omega_\alpha).
 \]
 Therefore, if $x\in\mathcal{F}(G,\{\mu\})$ satisfies $\mu^{\mathcal{L}}(x,\lambda)<0$, then there exists $\alpha\in\Delta$ such that
 $\mu^{\mathcal{L}}(\dot{w}x,\omega_\alpha)<0$. In other words, $x$ belongs to $\dot{w}^{-1}Z_{\alpha,<0}$.
 Hence we have $\{x\in \mathcal{F}(G,\{\mu\})\mid \mu^{\mathcal{L}}(x,\lambda)<0\}\subset \bigcup_{w\in W(G,S),\alpha\in\Delta}\dot{w}^{-1}Z_{\alpha,<0}$. 
 In particular, we have
 \begin{align*}
  \mathcal{F}(G,1,\{\mu\})^{\mathrm{ss}}&\subset \mathcal{F}(G,\{\mu\})^{\ad}\smallsetminus\bigcup_{\alpha\in\Delta,w\in W(G,S)}\dot{w}^{-1}Z_{\alpha,<0}^{\ad}\\
  &=\Bigl(\mathcal{F}(G,\{\mu\})\smallsetminus\bigcup_{\alpha\in\Delta,w\in W(G,S)}\dot{w}^{-1}Z_{\alpha,<0}\Bigr)^{\ad}\subset \mathcal{F}(G,\{\mu\})^{\mathrm{alg}\text{-}\mathrm{ss},\ad},
 \end{align*}
 as desired.

 Next we prove $\mathcal{F}(G,1,\{\mu\})^{\mathrm{st}}\subset \mathcal{F}(G,\{\mu\})^{\mathrm{alg}\text{-}\mathrm{st},\ad}$.
 For $x\in \mathcal{F}(G,\{\mu\})$ and $\lambda\in X_*(S_{\C_p})\smallsetminus \{0\}$, suppose $\mu^{\mathcal{L}}(x,\lambda)\le 0$.
 Let $w\in W(G,S)$ and $(c_\alpha)_{\alpha\in\Delta}\in (\prod_{\alpha\in\Delta}\Q_{\ge 0})\smallsetminus \{0\}$ be as above, and
 put $\Delta'=\{\alpha\in \Delta\mid c_\alpha>0\}$. Since $\lambda\neq 0$, $\Delta'$ is non-empty.
 By \cite[Definition 2.2 (i)]{MR1304906} and \cite[Proposition 2.3]{MR2198221}, we have
 \[
  0\ge \mu^{\mathcal{L}}(x,\lambda)=\mu^{\mathcal{L}}(\dot{w}x,w\lambda)=\mu^{\mathcal{L}}\Bigl(\dot{w}x,\sum_{\alpha\in\Delta'}c_\alpha\omega_\alpha\Bigr)=\sum_{\alpha\in\Delta'}c_\alpha\mu^{\mathcal{L}}(\dot{w}x,\omega_\alpha).
 \]
 Therefore, for some $\alpha\in\Delta'$ we have $\mu^{\mathcal{L}}(\dot{w}x,\omega_\alpha)\le 0$, namely, $x\in \dot{w}^{-1}Z_{\alpha,\le 0}$.
 Hence we have $\{x\in \mathcal{F}(G,\{\mu\})\mid \mu^{\mathcal{L}}(x,\lambda)\le 0\}\subset \bigcup_{w\in W(G,S),\alpha\in\Delta}\dot{w}^{-1}Z_{\alpha,\le 0}$, from which we conclude $\mathcal{F}(G,1,\{\mu\})^{\mathrm{st}}\subset \mathcal{F}(G,\{\mu\})^{\mathrm{alg}\text{-}\mathrm{st},\ad}$ in the same way
 as above.
\end{prf}

\begin{prop}\label{prop:algebraic-proper-action}
 The morphism 
 \[
  f\colon S_{\C_p}\times \mathcal{F}(G,\{\mu\})^{\mathrm{alg}\text{-}\mathrm{st}}\to \mathcal{F}(G,\{\mu\})^{\mathrm{alg}\text{-}\mathrm{ss}}\times \mathcal{F}(G,\{\mu\})^{\mathrm{alg}\text{-}\mathrm{st}};\ (g,x)\mapsto (gx,x)
 \]
 of algebraic varieties over $\C_p$ is proper.
\end{prop}

\begin{prf}
 By \cite[Corollary 2.5]{MR1304906}, the morphism 
 \[
  f'\colon S_{\C_p}\times \mathcal{F}(G,\{\mu\})^{\mathrm{alg}\text{-}\mathrm{st}}\to \mathcal{F}(G,\{\mu\})^{\mathrm{alg}\text{-}\mathrm{st}}\times \mathcal{F}(G,\{\mu\})^{\mathrm{alg}\text{-}\mathrm{st}};\ (g,x)\mapsto (gx,x)
 \]
 is proper. We apply \cite[Theorem 1.10]{MR1304906} to the action of $S_{\C_p}$ on $\mathcal{F}(G,\{\mu\})$; there exists a uniform categorical quotient
 $\phi\colon \mathcal{F}(G,\{\mu\})^{\mathrm{alg}\text{-}\mathrm{ss}}\to Y$, and $Y$ is quasi-projective, hence separated over $\C_p$.
 Further, there exists an open subscheme $Y^{\mathrm{st}}$ of $Y$ such that 
 $\phi^{-1}(Y^{\mathrm{st}})=\mathcal{F}(G,\{\mu\})^{\mathrm{alg}\text{-}\mathrm{st}}$.
 Since $Y$ is separated over $\C_p$, the diagonal map $\delta\colon Y\hookrightarrow Y\times Y$ is a closed immersion, so its base change
 $\delta'\colon Y^{\mathrm{st}}\hookrightarrow Y\times Y^{\mathrm{st}}$ is also a closed immersion.
 Let $D$ be the closed subscheme of $\mathcal{F}(G,\{\mu\})^{\mathrm{alg}\text{-}\mathrm{ss}}\times \mathcal{F}(G,\{\mu\})^{\mathrm{alg}\text{-}\mathrm{st}}$
 obtained as the pullback of $\delta'(Y^{\mathrm{st}})$ under $\phi\times\phi\colon \mathcal{F}(G,\{\mu\})^{\mathrm{alg}\text{-}\mathrm{ss}}\times \mathcal{F}(G,\{\mu\})^{\mathrm{alg}\text{-}\mathrm{st}}\to Y\times Y^{\mathrm{st}}$.
 Since $\delta'(Y^{\mathrm{st}})\subset Y^{\mathrm{st}}\times Y^{\mathrm{st}}$, $D$ is contained in 
 $\mathcal{F}(G,\{\mu\})^{\mathrm{alg}\text{-}\mathrm{st}}\times \mathcal{F}(G,\{\mu\})^{\mathrm{alg}\text{-}\mathrm{st}}$, and
 $f'$ factors as $S_{\C_p}\times \mathcal{F}(G,\{\mu\})^{\mathrm{alg}\text{-}\mathrm{st}}\xrightarrow{f''}D\hookrightarrow \mathcal{F}(G,\{\mu\})^{\mathrm{alg}\text{-}\mathrm{st}}\times \mathcal{F}(G,\{\mu\})^{\mathrm{alg}\text{-}\mathrm{st}}$ by definition.
 Since $f'$ is proper, so is $f''$. Therefore $f$ is also proper, for it is the composite of
 $S_{\C_p}\times \mathcal{F}(G,\{\mu\})^{\mathrm{alg}\text{-}\mathrm{st}}\xrightarrow{f''}D\hookrightarrow \mathcal{F}(G,\{\mu\})^{\mathrm{alg}\text{-}\mathrm{ss}}\times \mathcal{F}(G,\{\mu\})^{\mathrm{alg}\text{-}\mathrm{st}}$.
\end{prf}

\begin{prop}\label{prop:rigid-proper-action}
 Let $V\subset \mathcal{F}(G,1,\{\mu\})^{\mathrm{st}}$ and $U\subset \mathcal{F}(G,1,\{\mu\})^{\mathrm{ss}}=\mathcal{F}(G,1,\{\mu\})^{\mathit{wa}}$ be
 quasi-compact open subsets.
 \begin{enumerate}
  \item The set $\{\lambda\in X_{*,F}(S)\mid \lambda(\varpi)V\cap U\neq\varnothing\}$ is finite.
  \item The subset $\{g\in G(F)\mid gV\cap U\neq \varnothing\}$ of $G(F)$ is compact.
 \end{enumerate}
\end{prop}

\begin{prf}
 By Proposition \ref{prop:algebraic-proper-action}, the morphism 
 \[
  f^{\ad}\colon S^{\ad}_{\C_p}\times \mathcal{F}(G,\{\mu\})^{\mathrm{alg}\text{-}\mathrm{st},\ad}\to \mathcal{F}(G,\{\mu\})^{\mathrm{alg}\text{-}\mathrm{ss},\ad}\times \mathcal{F}(G,\{\mu\})^{\mathrm{alg}\text{-}\mathrm{st},\ad}
 \]
 of adic spaces over $\C_p$ is proper, hence quasi-compact. By Proposition \ref{prop:alg-st}, we have $V\subset \mathcal{F}(G,\{\mu\})^{\mathrm{alg}\text{-}\mathrm{st},\ad}$ and $U\subset \mathcal{F}(G,\{\mu\})^{\mathrm{alg}\text{-}\mathrm{ss},\ad}$.
 Therefore $(f^{\ad})^{-1}(U\times V)$ is a quasi-compact open subset of $S^{\ad}_{\C_p}\times \mathcal{F}(G,\{\mu\})^{\mathrm{alg}\text{-}\mathrm{st},\ad}$.
 Let $W$ be the image of $(f^{\ad})^{-1}(U\times V)$ under the first projection $\pr_1\colon S^{\ad}_{\C_p}\times \mathcal{F}(G,\{\mu\})^{\mathrm{alg}\text{-}\mathrm{st},\ad}\to S^{\ad}_{\C_p}$, which is also quasi-compact. We fix an isomorphism $S\cong \mathbb{G}_m^r$, where $r$ denotes the rank of $S$, 
 and define $(\mathbb{G}_m^\ad)^r_n=\{(x_1,\ldots,x_r)\in (\mathbb{G}_m^\ad)^r\mid \lvert\varpi^n\rvert\le \lvert x_i\rvert\le \lvert\varpi^{-n}\rvert\ (1\le i\le r)\}$ for each integer $n\ge 0$. Since $\{(\mathbb{G}_m^\ad)^r_n\}_{n\ge 0}$ is an open covering of $(\mathbb{G}_m^\ad)^r$, there exists an integer $n\ge 0$ such that $W\subset (\mathbb{G}_m^\ad)^r_n$. If $\lambda\in X_{*,F}(S)$ satisfies $\lambda(\varpi)V\cap U\neq \varnothing$, then $\lambda(\varpi)$ belongs to $W$, hence $\lambda(\varpi)$ lies in $(\mathbb{G}_m^\ad)^r_n$. Clearly there exist only finitely many such $\lambda$.

 (ii) By \cite[Th\'eor\`eme V.5.1 (4)]{MR2567785}, there exist a compact open subgroup $K$ of $G(F)$ and a finite subset $\Xi$ of $G(F)$
 such that $G(F)=\bigcup_{\lambda\in X_{*,F}(S),\xi\in\Xi}K\lambda(\varpi)\xi K$.
 We put $V'=\bigcup_{k\in K}kV$ and $U'=\bigcup_{k\in K}kU$. Since the stabilizers of $V$ and $U$ in $K$ are open
 by \cite[Lemma 2.2]{MR3861564}, $V'$ and $U'$ are quasi-compact.
 For each $\xi\in \Xi$, we define $\Lambda_\xi=\{\lambda\in X_{*,F}(S)\mid \lambda(\varpi)\xi V'\cap U'\neq \varnothing\}$, which is a finite set by (i).
 Then we have $\{g\in G(F)\mid gV\cap U\neq \varnothing\}\subset \bigcup_{\xi\in \Xi}\bigcup_{\lambda\in \Lambda_\xi}K\lambda(\varpi)\xi K$.
 Since $K\lambda(\varpi)\xi K$ is compact for each $\xi\in \Xi$ and $\lambda\in \Lambda_\xi$, the set
 $\{g\in G(F)\mid gV\cap U\neq \varnothing\}$ has finite volume with respect to a Haar measure of $G(F)$. On the other hand,
 by \cite[Lemma 2.2]{MR3861564}, the set $\{g\in G(F)\mid gV\cap U\neq \varnothing\}$ is bi-invariant under some compact open subgroup of $G(F)$.
 Therefore we conclude that $\{g\in G(F)\mid gV\cap U\neq \varnothing\}$ is compact.
\end{prf}

\begin{prop}\label{prop:st-subset-adm}
 Assume that $\mathcal{F}(G,1,\{\mu\})^a$ admits a Siegel domain. Then we have $\mathcal{F}(G,1,\{\mu\})^{\mathrm{st}}\subset \mathcal{F}(G,1,\{\mu\})^a$.
\end{prop}

\begin{prf}
 Let $U$ be a Siegel domain in $\mathcal{F}(G,1,\{\mu\})^a$. By \cite[Lemma 2.2]{MR3861564}, there exists a compact open subgroup $K'$ of $G(F)$
 such that $U$ is $K'$-invariant. Let $x\in \mathcal{F}(G,1,\{\mu\})^{\mathrm{st}}$ be an arbitrary point, and take a quasi-compact open neighborhood $V$
 of $x$ contained in $\mathcal{F}(G,1,\{\mu\})^{\mathrm{st}}$ (see Proposition \ref{prop:ss=wa} (ii)).
 Since the set $\{g\in G(F)\mid gV\cap U\neq\varnothing\}$ is left $K'$-invariant, by Proposition \ref{prop:rigid-proper-action} (ii)
 we can find $g_1,\ldots,g_r\in G(F)$ such that $\{g\in G(F)\mid gV\cap U\neq\varnothing\}=\bigcup_{i=1}^rK'g_i$.

 We put $V'=\bigcup_{i=1}^rV\cap g_i^{-1}U$. Since $\mathcal{F}(G,\{\mu\})^{\ad}$ is separated, $V\cap g_i^{-1}U$ is quasi-compact
 for every $i$. Therefore, $V'$ is a quasi-compact open subset of $V$. We prove that every classical point $y\in V(\cl)$ belongs to $V'$.
 Since $V(\cl)\subset \mathcal{F}(G,1,\{\mu\})^{\mathit{wa}}(\cl)=\mathcal{F}(G,1,\{\mu\})^a(\cl)$, $y$ lies in $\mathcal{F}(G,1,\{\mu\})^a$.
 Therefore, there exists $g\in G(F)$ such that $gy\in U$, as $U$ is a Siegel domain in $\mathcal{F}(G,1,\{\mu\})^a$.
 In particular, $gV\cap U$ is non-empty. Hence there exist $k\in K'$ and $1\le i\le r$ such that $g=kg_i$.
 Then we have $y\in V\cap g^{-1}U=V\cap g_i^{-1}k^{-1}U=V\cap g_i^{-1}U\subset V'$.
 Therefore, by Lemma \ref{lem:classical-points} (i) we conclude that $V=V'$.
 In particular, we have $x\in V=V'\subset \mathcal{F}(G,1,\{\mu\})^a$, as desired.
\end{prf}

Next we prove the inclusion $\mathcal{F}(G,1,\{\mu\})^{\mathit{wa}}\smallsetminus \mathcal{F}(G,1,\{\mu\})^{\mathrm{st}}\subset \mathcal{F}(G,1,\{\mu\})^a$
by using the induction hypothesis. For $\alpha\in\Delta$, let $M_\alpha$ denote the centralizer of $n_\alpha\omega_\alpha$,
which is a Levi subgroup of $P_\alpha$. We recall the construction in \cite[p.~233]{MR4325262} and \cite[p.~2895]{MR4612104}.
Let $W=W(G_{\overline{F}},T)$ be the absolute Weyl group. For $w\in W$, we put $\mathcal{F}(G,\{\mu\})_{\alpha,w}=P_{\alpha,\C_p}wP(\mu)_{\C_p}/P(\mu)_{\C_p}$, which is a locally closed subvariety of $\mathcal{F}(G,\{\mu\})$ defined over a finite extension of $\breve{E}$. 
Since $P_{\alpha,\C_p}wP(\mu)_{\C_p}/P(\mu)_{\C_p}\cong P_{\alpha,\C_p}/P_{\alpha,\C_p}\cap P(w\mu)_{\C_p}$,
the quotient morphism $P_{\alpha,\C_p}\to M_{\alpha,\C_p}$ induces a $P_{\alpha,\C_p}$-equivariant morphism 
$\pr_{\alpha,w}\colon \mathcal{F}(G,\{\mu\})_{\alpha,w}\to \mathcal{F}(M_\alpha,\{w\mu\})$, which is clearly defined over a finite extension of $\breve{E}$.
This is an affine fibration,
so in particular a smooth surjection. There is a decomposition $\mathcal{F}(G,\{\mu\})=\coprod_{w\in W_{P_{\alpha,\overline{F}}}\backslash W/W_{P(\mu)}}\mathcal{F}(G,\{\mu\})_{\alpha,w}$ into $P_{\alpha,\C_p}$-orbits.
Since $Z_{\alpha,=0}=Z_{\alpha,\le 0}\smallsetminus Z_{\alpha,<0}$ is $P_{\alpha,\C_p}$-invariant by Proposition \ref{prop:ss=wa} (i),
there exists a left $W_{P_{\alpha,\overline{F}}}$-invariant and right $W_{P(\mu)}$-invariant subset $W_\alpha$ of $W$ such that
$Z_{\alpha,=0}=\coprod_{w\in W_{P_{\alpha,\overline{F}}}\backslash W_\alpha/W_{P(\mu)}}\mathcal{F}(G,\{\mu\})_{\alpha,w}$ holds.

\begin{prop}\label{prop:alpha-w-closed}
 For every $\alpha\in\Delta$ and $w\in W_\alpha$, $\mathcal{F}(G,\{\mu\})_{\alpha,w}$ is a closed subvariety of $Z_{\alpha,=0}$.
 Moreover, $\mathcal{F}(G,1,\{\mu\})_{\alpha,w}^{\mathit{wa}}=\mathcal{F}(G,\{\mu\})_{\alpha,w}^{\ad}\cap \mathcal{F}(G,1,\{\mu\})^{\mathit{wa}}$ is
 a closed adic subspace of $\mathcal{F}(G,1,\{\mu\})^{\mathit{wa}}$,
 and $\mathcal{F}(G,1,\{\mu\})_{\alpha,w}^a=\mathcal{F}(G,\{\mu\})_{\alpha,w}^{\ad}\cap \mathcal{F}(G,1,\{\mu\})^a$ is
 a closed adic subspace of $\mathcal{F}(G,1,\{\mu\})^a$.
\end{prop}

\begin{prf}
 We take an integer $N>0$ such that $\mathcal{L}^{\otimes N}$ is very ample, and consider the embedding $\mathcal{F}(G,\{\mu\})\hookrightarrow\P(V)$
 corresponding to $\mathcal{L}^{\otimes N}$. Here $V$ is a finite-dimensional representation of $G_{\C_p}$, so $\mathbb{G}_m$ acts on $V$
 through $n_\alpha\omega_\alpha\colon \mathbb{G}_m\to G$. Let $V=\bigoplus_{m\in\Z}V_m$ be the weight decomposition with respect to this action,
 and put $V_{\ge 0}=\bigoplus_{m\ge 0}V_m$ and $V_{>0}=\bigoplus_{m>0}V_m$. By the definition of $\mu^{\mathcal{L}}(x,\omega_\alpha)$,
 we have $Z_{\alpha,=0}\subset \P(V_{\ge 0})\smallsetminus \P(V_{>0})$ (see \cite[proof of Proposition 2.2]{MR1304906}).
 The linear map $V_{\ge 0}=V_{>0}\oplus V_0\xrightarrow{\pr_2}V_0\hookrightarrow V$ induces a morphism 
 $\P(V_{\ge 0})\smallsetminus \P(V_{>0})\to \P(V)$ of algebraic varieties. We write $r\colon Z_{\alpha,=0}\to \P(V)$
 for the restriction of this morphism to $Z_{\alpha,=0}$. For $w\in W_\alpha$, it is easy to observe that the morphism
 $r\vert_{\mathcal{F}(G,\{\mu\})_{\alpha,w}}\colon \mathcal{F}(G,\{\mu\})_{\alpha,w}\to \P(V)$ coincides with the composite
 \begin{align*}
  \mathcal{F}(G,\{\mu\})_{\alpha,w}&\xrightarrow{\pr_{\alpha,w}}\mathcal{F}(M_\alpha,\{w\mu\})\cong M_{\alpha,\C_p}wP(\mu)_{\C_p}/P(\mu)_{\C_p}\hookrightarrow \mathcal{F}(G,\{\mu\})_{\alpha,w}\\
  &\hookrightarrow \P(V).
 \end{align*}
 We show $r^{-1}(M_{\alpha,\C_p}wP(\mu)_{\C_p}/P(\mu)_{\C_p})=\mathcal{F}(G,\{\mu\})_{\alpha,w}$.
 The inclusion $\supset$ is clear. We prove the reverse inclusion.
 Let $x$ be a point of $Z_{\alpha,=0}$, and assume that $r(x)\in M_{\alpha,\C_p}wP(\mu)_{\C_p}/P(\mu)_{\C_p}$.
 We take $w'\in W_\alpha$ such that $x\in \mathcal{F}(G,\{\mu\})_{\alpha,w'}$. 
 Then we have $r(x)\in (M_{\alpha,\C_p}wP(\mu)_{\C_p}/P(\mu)_{\C_p})\cap (M_{\alpha,\C_p}w'P(\mu)_{\C_p}/P(\mu)_{\C_p})\subset \mathcal{F}(G,\{\mu\})_{\alpha,w}\cap \mathcal{F}(G,\{\mu\})_{\alpha,w'}$. This implies $\mathcal{F}(G,\{\mu\})_{\alpha,w}=\mathcal{F}(G,\{\mu\})_{\alpha,w'}$, 
 hence $x\in \mathcal{F}(G,\{\mu\})_{\alpha,w}$. Therefore we obtain
 $r^{-1}(M_{\alpha,\C_p}wP(\mu)_{\C_p}/P(\mu)_{\C_p})\subset\mathcal{F}(G,\{\mu\})_{\alpha,w}$.

 Since $M_{\alpha,\C_p}wP(\mu)_{\C_p}/P(\mu)_{\C_p}\cong \mathcal{F}(M_\alpha,\{w\mu\})$ is projective, it is a closed subvariety of $\P(V)$.
 Therefore $\mathcal{F}(G,\{\mu\})_{\alpha,w}=r^{-1}(M_{\alpha,\C_p}wP(\mu)_{\C_p}/P(\mu)_{\C_p})$ is a closed subvariety of $Z_{\alpha,=0}$.

 By Proposition \ref{prop:ss=wa} (iii), we have $Z_{\alpha,=0}^{\ad}\cap \mathcal{F}(G,1,\{\mu\})^{\mathit{wa}}=Z_{\alpha,\le 0}^{\ad}\cap \mathcal{F}(G,1,\{\mu\})^{\mathit{wa}}$. Therefore, $Z_{\alpha,=0}^{\ad}\cap \mathcal{F}(G,1,\{\mu\})^{\mathit{wa}}$ is a closed adic subspace of $\mathcal{F}(G,1,\{\mu\})^{\mathit{wa}}$. As $\mathcal{F}(G,\{\mu\})_{\alpha,w}^{\ad}$ is a closed adic subspace of $Z_{\alpha,=0}^{\ad}$, we conclude that
 $\mathcal{F}(G,1,\{\mu\})_{\alpha,w}^{\mathit{wa}}=\mathcal{F}(G,\{\mu\})_{\alpha,w}^{\ad}\cap \mathcal{F}(G,1,\{\mu\})^{\mathit{wa}}$ is a closed adic subspace of $\mathcal{F}(G,1,\{\mu\})^{\mathit{wa}}$. Taking the intersection with $\mathcal{F}(G,1,\{\mu\})^a$, we also conclude that
 $\mathcal{F}(G,1,\{\mu\})_{\alpha,w}^a$ is a closed adic subspace of $\mathcal{F}(G,1,\{\mu\})^a$.
\end{prf}

\begin{prop}\label{prop:admissible-locus-reduction}
 Let $\alpha\in \Delta$ and $w\in W_\alpha$.
 \begin{enumerate}
  \item We have $\mathcal{F}(G,1,\{\mu\})^a_{\alpha,w}=(\pr_{\alpha,w}^{\ad})^{-1}(\mathcal{F}(M_\alpha,1,\{w\mu\})^a)$.
	Moreover, the image of $\mathcal{F}(G,1,\{\mu\})^a_{\alpha,w}$ under $\pr_{\alpha,w}^{\ad}\colon \mathcal{F}(G,\{\mu\})^{\ad}_{\alpha,w}\to \mathcal{F}(M_\alpha,\{w\mu\})^{\ad}$ is equal to $\mathcal{F}(M_\alpha,1,\{w\mu\})^a$. 
  \item We have $\mathcal{F}(G,1,\{\mu\})^{\mathit{wa}}_{\alpha,w}=(\pr_{\alpha,w}^{\ad})^{-1}(\mathcal{F}(M_\alpha,1,\{w\mu\})^{\mathit{wa}})$.
	Moreover, the image of $\mathcal{F}(G,1,\{\mu\})^{\mathit{wa}}_{\alpha,w}$ under $\pr_{\alpha,w}^{\ad}\colon \mathcal{F}(G,\{\mu\})^{\ad}_{\alpha,w}\to \mathcal{F}(M_\alpha,\{w\mu\})^{\ad}$ is equal to $\mathcal{F}(M_\alpha,1,\{w\mu\})^{\mathit{wa}}$. 
 \end{enumerate}
\end{prop}

To prove this proposition, we need some preparations.

\begin{lem}\label{lem:slope-0-condition}
 Let $\alpha\in \Delta$ and $w\in W$.
 \begin{enumerate}
  \item Let $\mathfrak{g}$ denote the Lie algebra of $G$, and consider the adjoint representation $\Ad\colon G\to \GL(\mathfrak{g})$.
	Let $\mathfrak{g}=\bigoplus_{a\in\Z}V_a$ be the weight decomposition with respect to $\Ad\circ (n_\alpha\omega_\alpha)$. Since $M_\alpha$
	centralizes $n_\alpha\omega_\alpha$, each $V_a$ can be regarded as a representation of $M_\alpha$. We write $\rho_a\colon M_\alpha\to \GL(V_a)$
	for it, and put $\chi_\alpha=\prod_{a\in\Z}(\det\rho_a)^a\in X_F^*(M_\alpha)$. Then, for every $x\in \mathcal{F}(G,\{\mu\})_{\alpha,w}$,
	we have $\mu^{\mathcal{L}}(x,\omega_\alpha)=-dn_\alpha^{-1}\langle\chi_\alpha,w\mu\rangle$.
  \item If $w\in W_\alpha$, we have $\langle \chi,w\mu\rangle=0$ for every $\chi\in X_F^*(M_\alpha)\otimes_\Z\Q$.
	In particular, for every $a\in \Z$ we have $\langle \det\rho_a,w\mu\rangle=0$.
 \end{enumerate}
\end{lem}

\begin{prf}
 (i) The action of $\mathbb{G}_m$ on $\mathfrak{g}_{\overline{F}}$ given by $w\mu\colon \mathbb{G}_m\to T\hookrightarrow G_{\overline{F}}\xrightarrow{\Ad}\GL(\mathfrak{g}_{\overline{F}})$ preserves the weight decomposition $\mathfrak{g}_{\overline{F}}=\bigoplus_{a\in\Z}V_{a,\overline{F}}$.
 Let $V_{a,\overline{F}}=\bigoplus_{b\in \Z}V_{a,b,\overline{F}}$ be the weight decomposition with respect to this action.
 By applying \cite[Lemma 2.2 (i)]{MR2198221} to the faithful representation $\Ad$ of $G$, we have
 $\mu^{\mathcal{L}}(x,n_\alpha\omega_\alpha)=-d\sum_{a,b\in\Z}ab\dim V_{a,b,\overline{F}}$.
 Here we use the assumption that $G$ is adjoint. Note also that our $\mathcal{L}$ is attached to the character $-d\mu^*$.
 On the other hand, we have $\langle \det\rho_a,w\mu\rangle=\sum_{b\in \Z}b\dim V_{a,b,\overline{F}}$. Hence we obtain
 $\langle \chi_\alpha,w\mu\rangle=\sum_{a,b\in \Z}ab\dim V_{a,b,\overline{F}}=-d^{-1}\mu^{\mathcal{L}}(x,n_\alpha\omega_\alpha)=-d^{-1}n_\alpha\mu^{\mathcal{L}}(x,\omega_\alpha)$. This concludes the proof.

 (ii) By (i) we have $\langle \chi_\alpha,w\mu\rangle=0$. Since $\dim X_F^*(M_\alpha)\otimes_\Z\Q=1$, it suffices to show $\chi_\alpha\neq 0$.
 Assume $\chi_\alpha=0$. Then we have $0=\langle \chi_\alpha,n_\alpha\omega_\alpha\rangle=\sum_{a\in\Z}a^2\dim V_a$. This implies $\mathfrak{g}=V_0$,
 which contradicts the non-triviality of $\omega_\alpha$. Therefore we conclude $\chi_\alpha\neq 0$.
\end{prf}

\begin{lem}\label{lem:FF-curve-bundle}
 Let $C$ be a complete algebraically closed extension of $\C_p$.
 As in the proof of Proposition \ref{prop:passing-to-adjoint} (i), $X_{\mathrm{FF},C^\flat}$ denotes the adic Fargues--Fontaine curve over $F$.
 \begin{enumerate}
  \item Let $H$ be a connected reductive group over $F$ and $\rho\colon H\to \GL_n$ a faithful representation of $H$.
	An $H$-bundle $\mathcal{E}$ over $X_{\mathrm{FF},C^\flat}$ satisfies $\nu_{\mathcal{E}}=0$ if and only if the vector bundle
	$\rho_*\mathcal{E}$ of rank $n$ is semistable of slope $0$.
  \item Let $\mathcal{V}$ be a vector bundle over $X_{\mathrm{FF},C^\flat}$, and $\mathcal{V}=\mathcal{V}_0\supset \cdots\supset \mathcal{V}_r=0$
	a filtration by vector subbundles. Assume that $\deg(\mathcal{V}_a/\mathcal{V}_{a+1})=0$ for every $0\le a\le r-1$.
	Then, $\mathcal{V}$ is semistable of slope $0$ if and only if $\mathcal{V}_a/\mathcal{V}_{a+1}$ is semistable of slope $0$
	for every $0\le a\le r-1$.
 \end{enumerate}
\end{lem}

\begin{prf}
 (i) By \cite[Th\'eor\`eme 8.2.10]{MR3917141}, the conditions $\nu_{\rho_*\mathcal{E}}=0$, $\rho_*\mathcal{E}$ is trivial, and
 $\rho_*\mathcal{E}$ is semistable of slope $0$ are all equivalent. Therefore, it suffices to prove that $\nu_{\mathcal{E}}=0$
 if and only if $\nu_{\rho_*\mathcal{E}}=0$.
 By \cite[8.2.1.2]{MR3917141} and \cite[\S 1.4]{MR4325262}, we may replace $F$ with a finite extension of it.
 Therefore, we may assume that $H$ is quasi-split.
 In this case, by \cite[Proposition 6.6]{MR4179595} we have $\nu_{\mathcal{E}}=w_0(-\nu_{[b]})$,
 where $w_0$ denotes the longest element of the relative Weyl group of $H$ and $[b]\in B(H)$ the element corresponding to the isomorphism class
 of $\mathcal{E}$ (see \cite[Th\'eor\`eme 5.1]{MR4179595}).
 In particular, $\nu_{\mathcal{E}}=0$ if and only if $\nu_{[b]}=0$. Similarly,
 we can see that $\nu_{\rho_*\mathcal{E}}=0$ if and only if $\nu_{[\rho(b)]}=0$. 
 Since $\nu_{[\rho(b)]}=\rho\circ\nu_{[b]}$ and $\rho$ is faithful, $\nu_{[b]}=0$ if and only if $\nu_{[\rho(b)]}=0$.
 This concludes the proof.

 (ii) Assume that $\mathcal{V}$ is semistable of slope $0$. By induction we can show that $\deg \mathcal{V}_a=0$ for every $0\le a\le r$.
 Therefore, $\mathcal{V}_a$ is also semistable of slope $0$ by definition.
 Hence $\mathcal{V}_a/\mathcal{V}_{a+1}$ is semistable of slope $0$ for every $0\le a\le r-1$
 by \cite[Th\'eor\`eme 5.5.4]{MR3917141}. The converse follows from \cite[Th\'eor\`eme 5.5.4]{MR3917141}.
\end{prf}

\begin{prf}[of Proposition \ref{prop:admissible-locus-reduction}]
 (i) Let $(C,C^+)$ be a complete algebraically closed affinoid field over $(\C_p,\mathcal{O}_{\C_p})$ and
 $x$ an element in $\mathcal{F}(G,\{\mu\})^{\ad}_{\alpha,w}(C,C^+)$. We prove that
 $x\in \mathcal{F}(G,1,\{\mu\})^a_{\alpha,w}(C,C^+)$ if and only if $\pr^{\ad}_{\alpha,w}(x)\in \mathcal{F}(M_\alpha,1,\{w\mu\})^a(C,C^+)$.
 Since $\mathcal{F}(G,\{\mu\})^{\ad}_{\alpha,w}$, $\mathcal{F}(G,1,\{\mu\})_{\alpha,w}^a$ and $\mathcal{F}(M_\alpha,1,\{w\mu\})^a$
 are partially proper (see Proposition \ref{prop:alpha-w-closed} and \cite[Proposition 3.2 (1)]{MR4325262}), we may assume $C^+=\mathcal{O}_C$.

 Let $\mathcal{E}_1$ denote the trivial $G$-bundle over $X_{\mathrm{FF},C^\flat}$, and $\mathcal{E}_{1,x}$ the modification of $\mathcal{E}_1$
 corresponding to $x$. Let $(\mathcal{E}_{1,x})_{P_\alpha}$ be the $P_\alpha$-reduction of $\mathcal{E}_{1,x}$ corresponding to the trivial
 $P_\alpha$-reduction of $\mathcal{E}_1$ under the bijection in \cite[Lemma 2.4]{MR4325262}. By \cite[Lemma 2.6]{MR4325262}, we have
 $(\mathcal{E}_{1,x})_{P_\alpha}\times^{P_\alpha}M_\alpha\cong \mathcal{E}_{1,\pr_{\alpha,w}(x)}$.

 Let $P_{\alpha,\GL}$ be the parabolic subgroup of $\GL(\mathfrak{g})$ associated with $\Ad\circ n_\alpha\omega_\alpha$.
 It is the stabilizer of the filtration $\{\bigoplus_{m\ge a}V_m\}_{a\in\Z}$ of $\mathfrak{g}$. Let $\mathcal{V}_x=\Ad_*\mathcal{E}_{1,x}$ be
 the vector bundle on $X_{\mathrm{FF},C^\flat}$ induced from $\Ad\colon G\to \GL(\mathfrak{g})$.
 The $P_\alpha$-reduction of $\mathcal{E}_{1,x}$ gives a $P_{\alpha,\GL}$-reduction of $\mathcal{V}_x$, which corresponds to
 a filtration $\cdots \supset \mathcal{V}_{x,-1}\supset \mathcal{V}_{x,0}\supset \mathcal{V}_{x,1}\supset \cdots$ by vector subbundles of $\mathcal{V}_x$.
 We have $(\Ad\rvert_{M_\alpha})_*\mathcal{E}_{1,\pr_{\alpha,w}(x)}\cong \bigoplus_{a\in\Z}\mathcal{V}_{x,a}/\mathcal{V}_{x,a+1}$ and
 $(\rho_a)_*\mathcal{E}_{1,\pr_{\alpha,w}(x)}\cong\mathcal{V}_{x,a}/\mathcal{V}_{x,a+1}$ for every $a\in \Z$.
 By \cite[p.~231]{MR4325262}, we have $c_1^{M_\alpha}(\mathcal{E}_{1,\pr_{\alpha,w}(x)})=(w\mu)^\sharp$, where $(w\mu)^\sharp$
 denotes the image of $w\mu$ in $\pi_1(M_\alpha)_{\Gal(\overline{F}/F)}$. Therefore, Lemma \ref{lem:slope-0-condition} (ii) tells us that
 $\deg(\mathcal{V}_{x,a}/\mathcal{V}_{x,a+1})=\deg((\rho_a)_*\mathcal{E}_{1,\pr_{\alpha,w}(x)})=\langle \det\rho_a,w\mu\rangle=0$.
 Now we can obtain the desired equivalence as follows:
 \begin{align*}
  &x\in \mathcal{F}(G,1,\{\mu\})^a_{\alpha,w}(C,\mathcal{O}_C)\stackrel{(1)}{\iff}
  \nu_{\mathcal{E}_{1,x}}=0\\
  &\qquad\stackrel{(2)}{\iff}\text{$\mathcal{V}_x=\Ad_*\mathcal{E}_{1,x}$ is semistable of slope $0$}\\
  &\qquad\stackrel{(3)}{\iff}\text{$\mathcal{V}_{x,a}/\mathcal{V}_{x,a+1}$ is semistable of slope $0$ for every $a\in \Z$}\\
  &\qquad\stackrel{(4)}{\iff}\text{$(\Ad\rvert_{M_\alpha})_*\mathcal{E}_{1,\pr_{\alpha,w}(x)}\cong \bigoplus_{a\in\Z}\mathcal{V}_{x,a}/\mathcal{V}_{x,a+1}$ is semistable of slope $0$}\\
  &\qquad\stackrel{(5)}{\iff}\nu_{\mathcal{E}_{1,\pr_{\alpha,w}(x)}}=0\stackrel{(6)}{\iff}\pr_{\alpha,w}(x)\in \mathcal{F}(M_\alpha,1,\{w\mu\})^a(C,\mathcal{O}_C).
 \end{align*}
 Here (1), (6) are nothing but the definition of the admissible locus (see \cite[Definition 3.1]{MR4325262}),
 (2), (5) follow from Lemma \ref{lem:FF-curve-bundle} (i), and (3), (4) are consequences of Lemma \ref{lem:FF-curve-bundle} (ii).
 Hence we have $\mathcal{F}(G,1,\{\mu\})^a_{\alpha,w}=(\pr_{\alpha,w}^{\ad})^{-1}(\mathcal{F}(M_\alpha,1,\{w\mu\})^a)$.
 By this equality and the surjectivity of $\pr_{\alpha,w}^{\ad}$, we also conclude
 $\pr_{\alpha,w}^{\ad}(\mathcal{F}(G,1,\{\mu\})^a_{\alpha,w})=\mathcal{F}(M_\alpha,1,\{w\mu\})^a$. 

 (ii) Let $x$ be a point of $\mathcal{F}(G,\{\mu\})^{\ad}_{\alpha,w}$. We prove that $x\in \mathcal{F}(G,1,\{\mu\})^{\mathit{wa}}_{\alpha,w}$
 if and only if $\pr_{\alpha,w}(x)\in \mathcal{F}(M_\alpha,1,\{w\mu\})^{\mathit{wa}}$. To prove the only if direction, we assume that 
 $x\in \mathcal{F}(G,1,\{\mu\})^{\mathit{wa}}_{\alpha,w}$ and $\pr_{\alpha,w}^{\ad}(x)\notin \mathcal{F}(M_\alpha,1,\{w\mu\})^{\mathit{wa}}$.
 As in the proof of Proposition \ref{prop:ss=wa} (iii), there exists a closed subvariety $Z_{M_\alpha}$ of $\mathcal{F}(M_\alpha,\{w\mu\})$
 defined over a finite extension of $\breve{E}$ such that 
 $\mathcal{F}(M_\alpha,1,\{w\mu\})^{\mathit{wa}}=\mathcal{F}(M_\alpha,\{w\mu\})^{\ad}\smallsetminus \bigcup_{m\in M_\alpha(F)}mZ_{M_\alpha}^{\ad}$.
 We can find $m\in M_\alpha(F)$ such that $\pr_{\alpha,w}^{\ad}(x)\in mZ_{M_\alpha}^{\ad}$. Replacing $x$ with $m^{-1}x$, we may assume that
 $\pr_{\alpha,w}^{\ad}(x)\in Z_{M_\alpha}^{\ad}$. 
 Then $x$ lies in $\mathcal{F}(G,1,\{\mu\})_{\alpha,w}^{\mathit{wa}}\cap (\pr_{\alpha,w}^{\ad})^{-1}Z_{M_\alpha}^{\ad}$,
 which is an open adic subspace of $(\pr_{\alpha,w}^{\ad})^{-1}Z_{M_\alpha}^{\ad}=(\pr_{\alpha,w}^{-1}Z_{M_\alpha})^{\ad}$.
 By (i), the set of classical points of this adic space satisfies
 \begin{align*}
  &(\mathcal{F}(G,1,\{\mu\})_{\alpha,w}^{\mathit{wa}}\cap (\pr_{\alpha,w}^{\ad})^{-1}Z_{M_\alpha}^{\ad})(\cl)=\mathcal{F}(G,1,\{\mu\})_{\alpha,w}^a(\cl)\cap ((\pr_{\alpha,w}^{\ad})^{-1}Z_{M_\alpha}^{\ad})(\cl)\\
  &\qquad=((\pr_{\alpha,w}^{\ad})^{-1}\mathcal{F}(M_\alpha,1,\{w\mu\})^a)(\cl)\cap ((\pr_{\alpha,w}^{\ad})^{-1}Z_{M_\alpha}^{\ad})(\cl)\\
  &\qquad=((\pr_{\alpha,w}^{\ad})^{-1}\mathcal{F}(M_\alpha,1,\{w\mu\})^a\cap (\pr_{\alpha,w}^{\ad})^{-1}Z_{M_\alpha}^{\ad})(\cl)\\
  &\qquad=((\pr_{\alpha,w}^{\ad})^{-1}(\mathcal{F}(M_\alpha,1,\{w\mu\})^a\cap Z_{M_\alpha}^{\ad}))(\cl)\\
  &\qquad\subset ((\pr_{\alpha,w}^{\ad})^{-1}(\mathcal{F}(M_\alpha,1,\{w\mu\})^{\mathit{wa}}\cap Z_{M_\alpha}^{\ad}))(\cl)=\varnothing;
 \end{align*}
 here we use the fact that $\pr_{\alpha,w}$ is defined over a finite extension of $\breve{E}$.
 Therefore, Lemma \ref{lem:classical-points} (ii) tells us that
 $\mathcal{F}(G,1,\{\mu\})_{\alpha,w}^{\mathit{wa}}\cap (\pr_{\alpha,w}^{\ad})^{-1}Z_{M_\alpha}^{\ad}=\varnothing$, 
 which contradicts $x\in \mathcal{F}(G,1,\{\mu\})_{\alpha,w}^{\mathit{wa}}\cap (\pr_{\alpha,w}^{\ad})^{-1}Z_{M_\alpha}^{\ad}$.

 Next, to prove the if direction, we assume that 
 $x\notin \mathcal{F}(G,1,\{\mu\})^{\mathit{wa}}_{\alpha,w}$ and $\pr_{\alpha,w}^{\ad}(x)\in \mathcal{F}(M_\alpha,1,\{w\mu\})^{\mathit{wa}}$.
 As in the proof of Proposition \ref{prop:ss=wa} (iii), there exists a closed subvariety $Z$ of $\mathcal{F}(G,\{\mu\})$
 defined over a finite extension of $\breve{E}$
 such that $\mathcal{F}(G,1,\{\mu\})^{\mathit{wa}}=\mathcal{F}(G,\{\mu\})^{\ad}\smallsetminus \bigcup_{g\in G(F)}gZ^{\ad}$.
 We can find $g\in G(F)$ such that $x\in gZ^{\ad}$. 
 Then $x$ lies in $gZ^{\ad}\cap (\pr_{\alpha,w}^{\ad})^{-1}\mathcal{F}(M_\alpha,1,\{w\mu\})^{\mathit{wa}}$,
 which is an open adic subspace of $(gZ\cap \mathcal{F}(G,\{\mu\})_{\alpha,w})^{\ad}$.
 By (i), the set of classical points of this adic space satisfies
 \begin{align*}
  &(gZ^{\ad}\cap (\pr_{\alpha,w}^{\ad})^{-1}\mathcal{F}(M_\alpha,1,\{w\mu\})^{\mathit{wa}})(\cl)\\
  &\qquad=gZ^{\ad}(\cl)\cap ((\pr_{\alpha,w}^{\ad})^{-1}\mathcal{F}(M_\alpha,1,\{w\mu\})^a)(\cl)\\
  &\qquad=gZ^{\ad}(\cl)\cap \mathcal{F}(G,1,\{\mu\})^a_{\alpha,w}(\cl)\subset gZ^{\ad}(\cl)\cap \mathcal{F}(G,1,\{\mu\})^{\mathit{wa}}(\cl)=\varnothing;
 \end{align*}
 here we again use the fact that $\pr_{\alpha,w}$ is defined over a finite extension of $\breve{E}$.
 Therefore, Lemma \ref{lem:classical-points} (ii) tells us that
 $gZ^{\ad}\cap (\pr_{\alpha,w}^{\ad})^{-1}\mathcal{F}(M_\alpha,1,\{w\mu\})^{\mathit{wa}}=\varnothing$, 
 which contradicts $x\in gZ^{\ad}\cap (\pr_{\alpha,w}^{\ad})^{-1}\mathcal{F}(M_\alpha,1,\{w\mu\})^{\mathit{wa}}$.
 This completes the proof of $\mathcal{F}(G,1,\{\mu\})^{\mathit{wa}}_{\alpha,w}=(\pr_{\alpha,w}^{\ad})^{-1}(\mathcal{F}(M_\alpha,1,\{w\mu\})^{\mathit{wa}})$.
 By this equality and the surjectivity of $\pr_{\alpha,w}^{\ad}$, we conclude that 
 $\pr_{\alpha,w}^{\ad}(\mathcal{F}(G,1,\{\mu\})^{\mathit{wa}}_{\alpha,w})=\mathcal{F}(M_\alpha,1,\{w\mu\})^{\mathit{wa}}$. 
 \end{prf}

 \begin{prop}\label{prop:nonst-subset-adm}
  Assume that Theorem \ref{thm:A-version} holds for every connected reductive group whose split semisimple rank is smaller than that of $G$,
  and that $\mathcal{F}(G,1,\{\mu\})^a$ admits a Siegel domain. 
  Then we have $\mathcal{F}(G,1,\{\mu\})^{\mathit{wa}}\smallsetminus \mathcal{F}(G,1,\{\mu\})^{\mathrm{st}}\subset \mathcal{F}(G,1,\{\mu\})^a$.
 \end{prop}

 \begin{prf}
  Fix a Siegel domain $U$ in $\mathcal{F}(G,1,\{\mu\})^a$. 
  Let $x$ be a point in $\mathcal{F}(G,1,\{\mu\})^{\mathit{wa}}\smallsetminus \mathcal{F}(G,1,\{\mu\})^{\mathrm{st}}$.
  By the definition of $\mathcal{F}(G,1,\{\mu\})^{\mathrm{st}}$, there exist $\alpha\in\Delta$, $w\in W_\alpha$ and $g\in G(F)$ such that
  $x\in g\mathcal{F}(G,1,\{\mu\})_{\alpha,w}^{\mathit{wa}}$. 
  Replacing $x$ with $g^{-1}x$, we may assume that $x\in \mathcal{F}(G,1,\{\mu\})_{\alpha,w}^{\mathit{wa}}$.
  We take a compact open subgroup $K$ of $G(F)$ such that $G(F)=P_\alpha(F)K$ (see \cite[Th\'eor\`eme V.5.1 (1)]{MR2567785}).
  Replacing $U$ with $\bigcup_{k\in K}kU$, we may assume that $U$ is $K$-invariant. 
  Then we have $\mathcal{F}(G,1,\{\mu\})^a=\bigcup_{g\in P_\alpha(F)}gU$.

  By Proposition \ref{prop:alpha-w-closed}, $U\cap \mathcal{F}(G,1,\{\mu\})^a_{\alpha,w}$ is a quasi-compact open subset of
  $\mathcal{F}(G,1,\{\mu\})^a_{\alpha,w}$. Let $U_{M_\alpha}$ denote the image of $U\cap \mathcal{F}(G,1,\{\mu\})^a_{\alpha,w}$ under
  the morphism $\pr^{\ad}_{\alpha,w}\colon \mathcal{F}(G,\{\mu\})^{\ad}_{\alpha,w}\to \mathcal{F}(M_\alpha,\{w\mu\})^{\ad}$.
  By Proposition \ref{prop:admissible-locus-reduction} (i), it is a quasi-compact open subset of $\mathcal{F}(M_\alpha,1,\{w\mu\})^a$.
  We prove that $U_{M_\alpha}$ is a Siegel domain in $\mathcal{F}(M_\alpha,1,\{w\mu\})^a$.
  Let $y$ be a point of $\mathcal{F}(M_\alpha,1,\{w\mu\})^a$. By Proposition \ref{prop:admissible-locus-reduction} (i),
  we can find $y'\in \mathcal{F}(G,1,\{\mu\})^a_{\alpha,w}$ such that $\pr_{\alpha,w}^{\ad}(y')=y$.
  As $\mathcal{F}(G,1,\{\mu\})^a=\bigcup_{g\in P_\alpha(F)}gU$, there exists $g\in P_\alpha(F)$ such that $y'\in gU$.
  Since $\mathcal{F}(G,\{\mu\})^{\ad}_{\alpha,w}$ is $P_{\alpha,\C_p}$-invariant,
  we have $g^{-1}y'\in U\cap \mathcal{F}(G,1,\{\mu\})^a_{\alpha,w}$.
  Let $\overline{g}$ denote the image of $g$ in $M_\alpha(F)$. Then we have
  $\overline{g}^{\,-1}y=\overline{g}^{\,-1}\pr_{\alpha,w}^{\ad}(y')=\pr_{\alpha,w}^{\ad}(g^{-1}y')\in U_{M_\alpha}$;
  hence $y\in \overline{g}U_{M_\alpha}\subset \bigcup_{m\in M_\alpha(F)}mU_{M_\alpha}$, as desired.

  Therefore, by the induction hypothesis on $M_\alpha$, we obtain $\mathcal{F}(M_\alpha,1,\{w\mu\})^a=\mathcal{F}(M_\alpha,1,\{w\mu\})^{\mathit{wa}}$.
  By Proposition \ref{prop:admissible-locus-reduction} (ii),
  we have $\pr_{\alpha,w}(x)\in \mathcal{F}(M_\alpha,1,\{w\mu\})^{\mathit{wa}}=\mathcal{F}(M_\alpha,1,\{w\mu\})^a$.
  Hence Proposition \ref{prop:admissible-locus-reduction} (i) tells us that $x\in \mathcal{F}(G,1,\{\mu\})^a_{\alpha,w}$, 
  which completes the proof.
 \end{prf}

 Theorem \ref{thm:A-version} immediately follows from Propositions \ref{prop:st-subset-adm} and \ref{prop:nonst-subset-adm}.

\section{Existence of Siegel domains}\label{sec:existence}
\subsection{Classification of fully Hodge--Newton decomposable pairs}\label{subsec:HN-dec}
In this section, we prove the implication (ii) $\Rightarrow$ (i) of Theorem \ref{thm:main} by a case-by-case argument.
By Proposition \ref{prop:passing-to-adjoint}, we may assume that $G$ is adjoint. By decomposing $G$ into $F$-simple factors,
we may also assume that $G$ is $F$-simple. If $\{\mu\}$ contains the trivial cocharacter,
then $\mathcal{F}(G,\{\mu\})$ is a single point, so clearly $\mathcal{F}(G,b,\{\mu\})^a$ admits a Siegel domain.
Therefore, we may assume that $\{\mu\}$ does not contain the trivial cocharacter. 
First, we give a classification of such pairs $(G,\{\mu\})$ which are fully Hodge--Newton decomposable.

\begin{prop}\label{prop:G-mu-FHND}
 Let $G$ be an adjoint $F$-simple connected reductive group over $F$. We write $G=\Res_{F'/F}G'$, where $F'$ is a finite extension of $F$ and
 $G'$ is an adjoint absolutely simple connected reductive group over $F'$. We put $\Sigma=\Hom_{\textnormal{$F$-alg}}(F',\overline{F})$
 and $G'_\tau=G'\otimes_{F',\tau}\tau(F')$ for each $\tau\in\Sigma$, so that we have $G_{\overline{F}}\cong \prod_{\tau\in\Sigma}(G'_\tau)_{\overline{F}}$.
 Let $\mu\in X_*(G)$ be a minuscule cocharacter, and write $(\mu_\tau)_\tau$ for the element of $\prod_{\tau\in \Sigma}X_*(G'_\tau)$ corresponding to $\mu$
 under the isomorphism above. If $\mu\neq 0$, then $(G,\mu)$ is fully Hodge--Newton decomposable if and only if one of the following holds:
 \begin{enumerate}
  \item There exists a unique $\tau\in\Sigma$ such that $\mu_\tau\neq 0$, and 
	$(G'_\tau,\mu_\tau)$ is one of the pairs in the list of \cite[\S 3.6]{MR4060001} other than $(A_{n-1},\omega_1^\vee+\omega_{n-1}^\vee)$,
	$(C_n,\omega_1^\vee)$ and $(C\text{-}BC_2,\omega_2^\vee)$.
  \item We have $G'=\PGL_n$, and there exist $\tau_1,\tau_2\in\Sigma$ with $\tau_1\neq \tau_2$ such that $\mu_{\tau_1}=\omega_1^\vee$,
	$\mu_{\tau_2}=\omega_{n-1}^\vee$ and $\mu_{\tau'}=0$ for $\tau'\in\Sigma\smallsetminus \{\tau_1,\tau_2\}$.
 \end{enumerate}
\end{prop}

\begin{prf}
 First note that $\mu$ is minuscule if and only if $\mu_\tau$ is minuscule for every $\tau\in\Sigma$.

 We fix a maximal torus $T'$ of $G'$ such that $T'\otimes_{F'}\breve{F}'$ is contained in a Borel subgroup $\breve{B}'$ of $G'\otimes_{F'}\breve{F}'$,
 and put $T'_\tau=T'\otimes_{F',\tau}\tau(F')$ for $\tau\in\Sigma$. Changing $\mu_\tau$ by conjugacy, we may assume that $\mu_\tau\in X_*(T'_\tau)$.
 We take $\mu_{\tau,+}$ in the conjugacy class of $\mu_\tau$ as in \cite[\S 2.2, p.~1733]{MR4476128} and write $\underline{\mu}{}_\tau$
 for its image in $X_*(T'_\tau)_{I_{\tau(F')}}$; here $I_{\tau(F')}$ denotes the inertia subgroup of
 the absolute Galois group $\Gamma_{\tau(F')}=\Gal(\overline{F}/\tau(F'))$.

 We would like to compare $\underline{\mu}{}_\tau$ for different $\tau$. 
 We write $F'_0$ for the maximal unramified extension of $F$ contained in $F'$, put $d=[F'_0:F]$ and
 fix an $F$-embedding $\iota\colon F'_0\hookrightarrow \overline{F}$. Let $\varphi\in \Gal(F'_0/F)$ denote the (arithmetic) Frobenius lift.
 For each $0\le i\le d-1$, we put $\Sigma_i=\{\tau\in\Sigma\mid \tau\vert_{F'_0}=\iota\circ\varphi^i\}$ and fix an element $\tau_i\in\Sigma_i$.
 For every $\tau\in\Sigma_i$, there exists $g_\tau\in I_F$ such that $\tau=g_\tau\circ\tau_i$. Indeed, if we regard $F'_0$ as a subfield of $\overline{F}$
 by $\iota$, we can find extensions $\widetilde{\tau}_i,\widetilde{\tau}\colon F^{\mathrm{ur}}F'\to \overline{F}$ of $\tau_i$, $\tau$
 such that $\widetilde{\tau}_i\vert_{F^{\mathrm{ur}}}=\widetilde{\tau}\vert_{F^{\mathrm{ur}}}=\varphi^i$.
 An element $g_\tau\in \Gamma_F$ such that $\widetilde{\tau}=g_\tau\circ\widetilde{\tau}_i$ belongs to $I_F$ and satisfies 
 $\tau=g_\tau\circ\tau_i$. By using this element, we define an isomorphism $X_*(T'_\tau)\xrightarrow{\cong}X_*(T'_{\tau_i})$ by
 $\lambda\mapsto g_\tau^{-1}\lambda:=(\id_{T'}\otimes g_\tau^{-1})\circ\lambda\circ(\id_{\mathbb{G}_m}\otimes g_\tau)$.
 We can check that it induces an isomorphism $X_*(T'_\tau)_{I_{\tau(F')}}\xrightarrow{\cong}X_*(T'_{\tau_i})_{I_{\tau_i(F')}}$,
 and the induced isomorphism is independent of the choice of $g_\tau$.
 Therefore, we obtain an element $(\underline{\mu}{}_\tau)_{\tau\in\Sigma_i}\in \prod_{\tau\in\Sigma_i}X_*(T'_{\tau_i})_{I_{\tau_i(F')}}$
 for each $0\le i\le d-1$.

 On the other hand, by applying the same argument to $T=\Res_{F'/F'_0}T'\subset \Res_{F'/F'_0}G'$ and $\breve{B}=\Res_{\breve{F}'/\breve{F}}\breve{B}'$,
 we obtain $\mu_i\in X_*(T_{\iota\circ\varphi^i})$ and
 $\underline{\mu}{}_i\in X_*(T_{\iota\circ\varphi^i})_{I_{\iota(F'_0)}}$ for each $0\le i\le d-1$.
 By $X_*(T_{\iota\circ\varphi^i})\cong \Ind_{\Gamma_{\tau_i(F')}}^{\Gamma_{\iota(F'_0)}}X_*(T'_{\tau_i})$ and 
 $\Gamma_{\iota(F'_0)}=\Gamma_{\tau_i(F')}I_{\iota(F'_0)}$, we have 
 $X_*(T_{\iota\circ\varphi^i})_{I_{\iota(F'_0)}}\cong X_*(T'_{\tau_i})_{I_{\tau_i(F')}}$
 and $\underline{\mu}{}_i=\sum_{\tau\in\Sigma_i}\underline{\mu}{}_\tau$, where the right hand side is regarded as an element of
 $X_*(T_{\iota\circ\varphi^i})_{I_{\iota(F'_0)}}$ under this isomorphism (see \cite[proof of Lemma 2.2]{MR4476128}). 

 Let $\widetilde{\Delta}$ denote the absolute local Dynkin diagram of $G'$ and $\delta$ the diagram automorphism (see \cite[Definition 2.6 (i)]{MR4476128}).
 Since the Weil restriction with respect to a totally ramified extension does not change the local Dynkin diagram (see \cite[1.12]{MR546588}),
 the Tits datum attached to $(G,\mu)$ is given by $(\Res_{F'_0/F}(\widetilde{\Delta},\delta),(\underline{\mu}{}_0,\ldots,\underline{\mu}{}_{d-1}))$.
 As $\mu\neq 0$, \cite[Theorem 3.4]{MR4476128} tells us that $(G,\mu)$ is fully Hodge--Newton decomposable if and only if one of the following holds:
 \begin{enumerate}
  \item[(1)] There exists a unique $i$ such that $\underline{\mu}{}_i\neq 0$, and the Coxeter datum attached to 
	     $(\widetilde{\Delta},\delta,\underline{\mu}{}_i)$ is one of the entries in the table of \cite[Theorem 3.3]{MR4476128}\footnote{Note that $(\widetilde{D}_n,\Ad(\tau_1),\omega_1^\vee)$ is missing in this table. See \cite[Remark 4.11]{arXiv240215310}.}.
  \item[(2)] We have $(\widetilde{\Delta},\delta)=(\widetilde{A}_{n-1},\id)$, and there exist $i$, $j$ with $i\neq j$ such that
	     $\underline{\mu}{}_i=\omega_1^\vee$, $\underline{\mu}{}_j=\omega_{n-1}^\vee$ and $\underline{\mu}{}_{i'}=0$ for $i'\neq i,j$.
 \end{enumerate}
 First we consider the case (1) with $\underline{\mu}{}_i$ minuscule. By applying \cite[Lemma 2.2]{MR4476128} to $\Res_{F'/F'_0}G'$ and
 $\mu_i\in X_*(T_{\iota\circ\varphi^i})$, we conclude that there exists a unique $\tau\in\Sigma_i$ such that $\mu_\tau\neq 0$,
 and this $\mu_\tau$ is minuscule. Moreover, since $\underline{\mu}{}_\tau=\underline{\mu}{}_i$, the Coxeter datum attached to $(G'_\tau,\mu_\tau)$ is
 one of the entries in the table of \cite[Theorem 3.3]{MR4476128}. Such a pair $(G'_\tau,\mu_\tau)$ is classified in \cite[\S 3.6]{MR4060001}.
 Since $\mu_\tau$ is minuscule, $(A_{n-1},\omega_1^\vee+\omega_{n-1}^\vee)$ and $(C_n,\omega_1^\vee)$ do not occur.
 The pair $(C\text{-}BC_2,\omega_2^\vee)$ also does not occur. Indeed, the group of type $C\text{-}BC_2$ has $4$ minuscule cocharacters
 $\omega_1^\vee,\ldots,\omega_4^\vee$, among which $\omega_1^\vee$ and $\omega_4^\vee$ are minute
 (see \cite[Definition 3.2]{MR4060001} for the definition of minute cocharacters, and \cite[Theorem 3.3]{MR4060001} for the relation with
 full Hodge--Newton decomposability);
 by \cite[p.~525, Table 2]{MR4086580}, the images of $\omega_1^\vee$ and $\omega_4^\vee$ in $X_*(T'_\tau)_{I_{\tau(F')}}$ are $\omega_1^\vee$.
 Let $0\le i'\le d-1$ be an integer with $i'\neq i$. Since $\underline{\mu}{}_{i'}=0$, by applying \cite[Lemma 2.3]{MR4476128} to
 $\Res_{F'/F'_0}G'$ and $\mu_{i'}\in X_*(T_{\iota\circ\varphi^{i'}})$, we have $\mu_{\tau'}=0$ for every $\tau'\in\Sigma_{i'}$.
 Therefore the condition (i) in the proposition is satisfied.

 Next we consider the case (1) with $\underline{\mu}{}_i$ non-minuscule. Then $(\widetilde{\Delta},\delta,\underline{\mu}{}_i)$ is
 either $(\widetilde{A}_{n-1},\id,\omega_1^\vee+\omega_{n-1}^\vee)$ or $(\widetilde{C}_n,\id,\omega_1^\vee)$.
 If $(\widetilde{\Delta},\delta,\underline{\mu}{}_i)=(\widetilde{A}_{n-1},\id,\omega_1^\vee+\omega_{n-1}^\vee)$, 
 the type of $G'$ is $A_{n-1}$ or $C\text{-}BC_1$. If the type of $G'$ is $A_{n-1}$, the inertia group $I_{F'}$ acts on $X_*(T')$ trivially,
 hence we have $\omega_1^\vee+\omega_{n-1}^\vee=\sum_{\tau\in\Sigma_i}\mu_{\tau,+}$. Therefore, there exist $\tau_1,\tau_2\in\Sigma_i$
 with $\tau_1\neq\tau_2$ such that $\mu_{\tau_1,+}=\omega_1^\vee$, $\mu_{\tau_2,+}=\omega_{n-1}^\vee$ and $\mu_{\tau'}=0$ for $\tau'\in \Sigma_i\smallsetminus \{\tau_1,\tau_2\}$. Moreover, we have $\mu_{\tau'}=0$ for every $\tau'\in \Sigma_{i'}$ with $i'\neq i$ by \cite[Lemma 2.3]{MR4476128}.
 Therefore the condition (ii) in the proposition is satisfied.
 Assume that the type of $G'$ is $C\text{-}BC_1$. In particular we have $n=2$. By \cite[p.~525, Table 2]{MR4086580}, if $\mu_\tau\neq 0$, then
 $\underline{\mu}{}_\tau$ is equal to $2\omega_1^\vee$. Hence the equality $\omega_1^\vee+\omega_{n-1}^\vee=\sum_{\tau\in\Sigma_i}\underline{\mu}{}_\tau$
 implies that there exists a unique $\tau\in \Sigma_i$ such that $\mu_\tau\neq 0$.
 Moreover, we have $\mu_{\tau'}=0$ for every $\tau'\in \Sigma_{i'}$ with $i'\neq i$ by \cite[Lemma 2.3]{MR4476128}.
 Therefore the condition (i) in the proposition is satisfied.

 If $(\widetilde{\Delta},\delta,\underline{\mu}{}_i)=(\widetilde{C}_n,\id,\omega_1^\vee)$,
 the type of $G'$ is one of $C_n$, $C\text{-}BC_n$ and $C\text{-}B_n$. If the type of $G'$ is $C_n$, the inertia group $I_{F'}$ acts on $X_*(T')$ trivially,
 hence we have $\omega_1^\vee=\sum_{\tau\in\Sigma_i}\mu_{\tau,+}$. Since $\mu_{\tau,+}$ is minuscule, it is equal to $\omega_n^\vee$ or $0$.
 Therefore we have $\omega_1^\vee\in \Z\omega_n^\vee$, which is a contradiction.
 Assume that the type of $G'$ is either $C\text{-}BC_n$ or $C\text{-}B_n$.
 Note that \cite[p.~525, Table 2]{MR4086580} gives the list of minuscule cocharacters of such $G'$ and their images in $X_*(T')_{I_{F'}}$.
 By this table and the linear independence of $\omega_1^\vee,\ldots,\omega_n^\vee\in X_*(T')_{I_{F'}}$, we conclude from the equality
 $\omega_1^\vee=\sum_{\tau\in \Sigma_i}\underline{\mu}{}_\tau$ that there exists a unique $\tau\in\Sigma_i$ such that $\mu_\tau\neq 0$.
 Moreover, we have $\mu_{\tau'}=0$ for every $\tau'\in \Sigma_{i'}$ with $i'\neq i$ by \cite[Lemma 2.3]{MR4476128}.
 Therefore the condition (i) in the proposition is satisfied.

 Finally, we consider the case (2). By the same argument as in the case (1) with $\underline{\mu}{}_i$ minuscule, we can prove that
 there exist $\tau_1\in\Sigma_i$ and $\tau_2\in\Sigma_j$ such that $\mu_{\tau_1}=\omega_1^\vee$, $\mu_{\tau_2}=\omega_{n-1}^\vee$ and
 $\mu_{\tau'}=0$ for $\tau'\in \Sigma\smallsetminus \{\tau_1,\tau_2\}$. Therefore the condition (ii) in the proposition is satisfied.

 Now we conclude that if $(G,\mu)$ is fully Hodge--Newton decomposable, then one of the conditions (i) and (ii) is satisfied.
 The converse is also clear from the description above.
\end{prf}

We consider the case (i) of Proposition \ref{prop:G-mu-FHND} further. By \cite[p.~35]{MR1393439},
we may replace $(F,G,\mu)$ with $(\tau(F'),G'_\tau,\mu_\tau)$. The list in \cite[\S 3.6]{MR4060001}, excluding $(A_{n-1},\omega_1^\vee+\omega_{n-1}^\vee)$,
$(C_n,\omega_1^\vee)$ and $(C\text{-}BC_2,\omega_2^\vee)$, is as follows:
\begin{itemize}
 \item $(A_n,\omega_1^\vee)$: the Lubin-Tate case.
 \item $({}^d\!A_{d-1},\omega_1^\vee)$: the dual of $(A_{d-1},\omega_1^\vee)$.
 \item $({}^2\!A'_n,\omega_1^\vee)$: the unramified quasi-split $U(1,n)$ case.
 \item $({}^2\!A''_{2m-1},\omega_1^\vee)$: the dual of $({}^2\!A'_{2m-1},\omega_1^\vee)$.
 \item $(A_3,\omega_2^\vee)$: the same as $(D_3,\omega_1^\vee)$.
 \item $({}^2\!A_3,\omega_2^\vee)$: the same as $({}^2D'_3,\omega_1^\vee)$.
 \item $({}^2\!A'_3,\omega_2^\vee)$: the same as $({}^2D_3,\omega_1^\vee)$.
 \item $(B_n,\omega_1^\vee)$: the split $\SO(2,2n-1)$ case.
 \item $(B\text{-}C_n,\omega_1^\vee)$: the ramified quasi-split $U(1,2n-1)$ case.
 \item $({}^2\!B_n,\omega_1^\vee)$: the dual of $(B_n,\omega_1^\vee)$.
 \item $({}^2\!B\text{-}C_n,\omega_1^\vee)$: the dual of $(B\text{-}C_n,\omega_1^\vee)$.
 \item $(C\text{-}B_n,\omega_1^\vee)$: the ramified quasi-split $\SO(2,2n)$ case.
 \item $(C\text{-}BC_n,\omega_1^\vee)$: the ramified quasi-split $U(1,2n)$ case.
 \item $(C_2,\omega_2^\vee)$: the same as $(B_2,\omega_1^\vee)$.
 \item $(C\text{-}B_2,\omega_2^\vee)$: the same as $(B\text{-}C_2,\omega_1^\vee)$.
 \item $({}^2C_2,\omega_2^\vee)$: the same as $({}^2\!B_2,\omega_1^\vee)$.
 \item $({}^2C\text{-}B_2,\omega_2^\vee)$: the same as $({}^2\!B\text{-}C_2,\omega_1^\vee)$.
 \item $(D_n,\omega_1^\vee)$: the split $\SO(2,2n-2)$ case.
 \item $({}^2\!D_n,\omega_1^\vee)$: the unramified quasi-split $\SO(2,2n-2)$ case.
 \item $({}^2\!D'_n,\omega_1^\vee)$: the dual of $(D_n,\omega_1^\vee)$.
\end{itemize}
Therefore, by Proposition \ref{prop:twin-tower-principle}, we have only to consider the following 8 cases:
\begin{itemize}
 \item $(A_n,\omega_1^\vee)$: the Lubin-Tate case.
 \item $({}^2\!A'_n,\omega_1^\vee)$: the unramified quasi-split $U(1,n)$ case.
 \item $(B_n,\omega_1^\vee)$: the split $\SO(2,2n-1)$ case.
 \item $(B\text{-}C_n,\omega_1^\vee)$: the ramified quasi-split $U(1,2n-1)$ case.
 \item $(C\text{-}B_n,\omega_1^\vee)$: the ramified quasi-split $\SO(2,2n)$ case.
 \item $(C\text{-}BC_n,\omega_1^\vee)$: the ramified quasi-split $U(1,2n)$ case.
 \item $(D_n,\omega_1^\vee)$: the split $\SO(2,2n-2)$ case.
 \item $({}^2\!D_n,\omega_1^\vee)$: the unramified quasi-split $\SO(2,2n-2)$ case.
\end{itemize}
Assume that $(G,\mu)$ is one of the above, and let $b$ be a basic element of $G(\breve{F})$ whose image $[b]$ in $B(G)$ lies in $B(G,\mu)$.
Let $(J_b,b'',\{\mu^{-1}\})$ be the dual triple attached to $(G,b,\{\mu\})$. By Proposition \ref{prop:twin-tower-principle} again,
it suffices to find a Siegel domain in $\mathcal{F}(J_b,b'',\{\mu^{-1}\})^{\mathit{wa}}$.
Further, by Proposition \ref{prop:passing-to-inner-form} we have
$\mathcal{F}(J_b,b'',\{\mu^{-1}\})^{\mathit{wa}}\cong \mathcal{F}(G,1,\{\mu^{-1}\})^{\mathit{wa}}$.
Therefore, to show (ii) $\Rightarrow$ (i) of Theorem \ref{thm:main}, it suffices to prove the following proposition:

\begin{prop}\label{prop:case-by-case}
 Let $(G,\mu)$ be one of the following:
 \begin{enumerate}
  \item (the Drinfeld case) $G=\PGL_n$ with $n\ge 2$, and $\mu=\omega_1^\vee$.
  \item (the orthogonal case) $G$ is the adjoint quotient of a quasi-split special orthogonal group in $n$ variables with $n\ge 5$
	(if $n$ is even, there are the split, the unramified, and the ramified cases), and
	$\mu=\omega_1^\vee$.
  \item (the unitary case) $G$ is the adjoint quotient of a quasi-split unitary group in $n$ variables with $n\ge 3$
	(it depends on a quadratic extension $E/F$), and $\mu=\omega_1^\vee$.
  \item (the Hilbert--Blumenthal case) $G=\Res_{E/F}\PGL_n$ with $n\ge 2$ and a finite extension $E/F$, and $\mu=(\mu_\tau)_{\tau\in\Sigma}$ with
	$\Sigma=\Hom_{F\text{-}\mathrm{alg}}(E,\overline{F})$ such that there exist $\tau_1,\tau_2\in \Sigma$
	satisfying $\tau_1\neq\tau_2$, $\mu_{\tau_1}=\omega_1^\vee$, $\mu_{\tau_2}=\omega_{n-1}^\vee$
	and $\mu_{\tau'}=0$ for $\tau'\in \Sigma\smallsetminus \{\tau_1,\tau_2\}$.
 \end{enumerate}
 Then $\mathcal{F}(G,1,\{\mu^{-1}\})^{\mathit{wa}}$ admits a Siegel domain.
\end{prop}

\subsection{The Drinfeld case}\label{subsec:Drinfeld}
In this subsection, we treat the Drinfeld case. Let $G$ be $\PGL_n$ with $n\ge 2$.
The cocharacter $\mu=\omega_1^\vee$ is given by $z\mapsto [\diag(z,1,\ldots,1)]$.
Since the automorphism $G\to G$; $g\mapsto {}^t\!g^{-1}$ maps $\mu^{-1}$ to $\mu$, it suffices to consider $\mathcal{F}(G,1,\{\mu\})^{\mathit{wa}}$
in place of $\mathcal{F}(G,1,\{\mu^{-1}\})^{\mathit{wa}}$. For simplicity we write $\mathcal{F}$ and $\mathcal{F}^{\mathit{wa}}$ for 
$\mathcal{F}(G,\{\mu\})$ and $\mathcal{F}(G,1,\{\mu\})^{\mathit{wa}}$, respectively.
Since $P(\mu)$ is the stabilizer of $\langle e_1\rangle_F\subset F^n$ (here $e_1,\ldots,e_n$ denote the standard basis of $F^n$),
we have $\mathcal{F}=\P^{n-1}$. By using Proposition \ref{prop:ss=wa} (iii), we can observe that
$\mathcal{F}^{\mathit{wa}}=(\P^{n-1})^{\ad}\smallsetminus \bigcup_{\text{$H\subset \P^{n-1}$: $F$-rational hyperplane}}H^{\ad}_{\C_p}$;
namely, $\mathcal{F}^{\mathit{wa}}$ is the $(n-1)$-dimensional Drinfeld upper half space.
The existence of a Siegel domain in the Drinfeld upper half space is well-known (for example we can prove it by using Deligne's formal model).
However we give an elementary proof, which is helpful for understanding the arguments in the orthogonal and the unitary cases.

Let $(C,C^+)$ be a complete affinoid field over $(\C_p,\mathcal{O}_{\C_p})$, and let $\lvert-\rvert$ denote the continuous valuation on $C$
corresponding to the closed point of $\Spa(C,C^+)$.
Let $(-,-)$ denote the symmetric pairing $F^n\times F^n\to F$; $((a_i),(b_i))\mapsto \sum_{i=1}^n a_ib_i$. We use the same symbol
to denote the base change of the pairing to $C$. For $x\in C^n$ and $a\in F^n$, we put $\llparenthesis x,a\rrparenthesis=\lvert (x,a)\rvert$.

In the following we regard an element of $(\P^{n-1})^{\ad}(C,C^+)=\P^{n-1}(C)$ as a $1$-dimensional subspace of $C^n$.
For $L\in \P^{n-1}(C)$ (resp.\ $W\in\P^{n-1}(F)$), we put $\Lambda_L=(C^+)^n\cap L$ (resp.\ $\Lambda_W=\mathcal{O}_F^n\cap W$),
which is a free $C^+$-module (resp.\ $\mathcal{O}_F$-module) of rank $1$.
For $L\in \P^{n-1}(C)$ and $W\in\P^{n-1}(F)$, we define
$\llparenthesis L,W\rrparenthesis=\max\{\llparenthesis x,a\rrparenthesis \mid x\in \Lambda_L,a\in \Lambda_W\}$.
Since $\Lambda_L$ and $\Lambda_W$ are finitely generated, the maximum always exists.

\begin{prop}\label{prop:Drinfeld-tubular-nbhd}
 \begin{enumerate}
  \item For $W\in\P^{n-1}(F)$ and an integer $m\ge 1$, there exists a unique quasi-compact open subset $U_{m,W}$ of $\mathcal{F}^{\ad}$
	satisfying the following: for every complete affinoid field $(C,C^+)$ over $(\C_p,\mathcal{O}_{\C_p})$, we have 
	\[
	 U_{m,W}(C,C^+)=\{L\in (\P^{n-1})^{\ad}(C,C^+)\mid \llparenthesis L,W\rrparenthesis\ge \lvert\varpi^m\rvert\}.
	\]
  \item The subset $U_m=\bigcap_{W\in \P^{n-1}(F)}U_{m,W}$ of $\mathcal{F}^{\ad}$ is a quasi-compact open subset of $\mathcal{F}^{\mathit{wa}}$.
  \item We have $\mathcal{F}^{\mathit{wa}}=\bigcup_{m\ge 1}U_m$.
  \item For $g\in\GL_n(\mathcal{O}_F)$, we have $gU_{m,W}=U_{m,{}^t\!g^{-1}W}$ and $gU_m=U_m$.
 \end{enumerate}
\end{prop}

\begin{prf}
 (i) Let $\mathcal{E}_n=(\A^n)^{\ad}\times\mathcal{F}^{\ad}$ be the trivial vector bundle of rank $n$ over $\mathcal{F}^{\ad}=(\P^{n-1})^{\ad}$,
 and $\mathcal{E}^\circ_n=\mathbb{D}^n\times\mathcal{F}^{\ad}$ an open adic subspace of it. The universal line bundle $\mathcal{L}$ over
 $\mathcal{F}^{\ad}$ can be regarded as a subbundle of $\mathcal{E}_n$.
 We also consider the trivial vector bundle $\mathcal{E}_W=W_{\C_p}^{\ad}\times \mathcal{F}^{\ad}\subset \mathcal{E}_n$.
 We put $\mathcal{L}^\circ=\mathcal{L}\cap \mathcal{E}_n^\circ$ and $\mathcal{E}_W^\circ=\mathcal{E}_W\cap \mathcal{E}_n^\circ$.

 The pairing $(-,-)$ induces a morphism $\mathcal{E}_n^\circ\times_{\mathcal{F}^{\ad}}\mathcal{E}_n^\circ\to \mathcal{E}_1^\circ$.
 Let $\mathcal{U}_{m,W}$ be the inverse image of $\{z\in \mathbb{D}^1\mid \lvert z\rvert\ge \lvert\varpi^m\rvert\}\times \mathcal{F}^{\ad}\subset \mathcal{E}_1^\circ$ under the restriction $\mathcal{L}^\circ\times_{\mathcal{F}^{\ad}}\mathcal{E}_W^\circ\to \mathcal{E}_1^\circ$ of this morphism.
 Then, $\mathcal{U}_{m,W}$ is a quasi-compact open subset of $\mathcal{L}^\circ\times_{\mathcal{F}^{\ad}}\mathcal{E}_W^\circ$.
 Let $U_{m,W}$ be the image of $\mathcal{U}_{m,W}$ under the projection $\pr\colon \mathcal{L}^\circ\times_{\mathcal{F}^{\ad}}\mathcal{E}_W^\circ\to \mathcal{F}^{\ad}$. It is a quasi-compact open subset of $\mathcal{F}^{\ad}$, since $\pr$ is smooth. By construction, $U_{m,W}$ obviously satisfies the condition in (i). The uniqueness of $U_{m,W}$ is also clear.

 (ii) Let $W$ and $W'$ be elements of $\P^{n-1}(F)$. We prove that $U_{m,W}=U_{m,W'}$ holds if the images of $\Lambda_W$ and $\Lambda_{W'}$ in 
 $(\mathcal{O}_F/\varpi^{m+1}\mathcal{O}_F)^n$ are the same. Let $(C,C^+)$ be a complete affinoid field over $(\C_p,\mathcal{O}_{\C_p})$.
 For $L\in U_{m,W}(C,C^+)$, we take $x\in\Lambda_L$ and $a\in\Lambda_W$ satisfying $\llparenthesis x,a\rrparenthesis\ge \lvert\varpi^m\rvert$. 
 By assumption, there exist $a'\in\Lambda_{W'}$ and $b\in \mathcal{O}_F^n$ such that $a=a'+\varpi^{m+1}b$.
 Then we have $\lvert(x,\varpi^{m+1}b)\rvert\le \lvert\varpi^{m+1}\rvert<\lvert\varpi^m\rvert\le \lvert(x,a)\rvert$,
 hence $\llparenthesis x,a'\rrparenthesis=\lvert (x,a)-(x,\varpi^{m+1}b)\rvert=\lvert(x,a)\rvert\ge \lvert\varpi^m\rvert$.
 Therefore we obtain $L\in U_{m,W'}(C,C^+)$. Namely, we have $U_{m,W}\subset U_{m,W'}$. The reverse inclusion $U_{m,W}\supset U_{m,W'}$ can be proven
 in the same way. Since $(\mathcal{O}_F/\varpi^{m+1}\mathcal{O}_F)^n$ is a finite set, the intersection $\bigcap_{W\in\P^{n-1}(F)}U_{m,W}$ is
 essentially finite. Therefore $U_m$ is a quasi-compact open subset of $\mathcal{F}^{\ad}$.
 For $W\in \P^{n-1}(F)$, we write $W^\perp\subset F^n$ for the orthogonal complement of it with respect to the pairing $(-,-)$.
 As $\P(W^\perp)^{\ad}\cap U_{m,W}=\varnothing$, we have
 $U_m\subset \bigcap_{W\in\P^{n-1}(F)}(\mathcal{F}^{\ad}\smallsetminus \P(W^\perp)^{\ad})=\mathcal{F}^{\mathit{wa}}$.
 This completes the proof.

 (iii) Let $(C,C^+)$ be a complete affinoid field over $(\C_p,\mathcal{O}_{\C_p})$ and $L\in \mathcal{F}^{\mathit{wa}}(C,C^+)$.
 For $W\in \P^{n-1}(F)$, we have $\llparenthesis L,W\rrparenthesis\neq 0$; hence there exists an integer $m_W\ge 1$ such that
 $\llparenthesis L,W\rrparenthesis\ge \lvert\varpi^{m_W}\rvert$.
 The subset $V_W=\{W'\in\P^{n-1}(F)\mid \Lambda_{W'}\equiv \Lambda_W\pmod{\varpi^{m_W+1}}\}$ of $\P^{n-1}(F)$ is an open neighborhood of $W$
 with respect to the $\varpi$-adic topology on $\P^{n-1}(F)$.
 Since $\P^{n-1}(F)$ is compact, there exist $W_1,\ldots,W_r\in \P^{n-1}(F)$ such that $\P^{n-1}(F)=\bigcup_{i=1}^rV_{W_i}$.
 We put $m=\max\{m_{W_1},\ldots,m_{W_r}\}$ and prove $L\in U_m(C,C^+)$. For $W\in \P^{n-1}(F)$, we can find $1\le i\le r$ satisfying $W\in V_{W_i}$.
 By the proof of (ii), we have $U_{m_{W_i},W}=U_{m_{W_i},W_i}$. 
 Hence we have $L\in U_{m_{W_i},W_i}(C,C^+)=U_{m_{W_i},W}(C,C^+)\subset U_{m,W}(C,C^+)$, as desired.

 (iv) The claim immediately follows from the equality $\llparenthesis g^{-1}L,W\rrparenthesis=\llparenthesis L,{}^t\!g^{-1}W\rrparenthesis$
 for $L\in\mathcal{F}^{\ad}(C,C^+)$ and $W\in \P^{n-1}(F)$.
\end{prf}

To show Proposition \ref{prop:case-by-case} (i), it suffices to prove the following proposition.

\begin{prop}\label{prop:Drinfeld-case}
 We have $\mathcal{F}^{\mathit{wa}}=\bigcup_{g\in\GL_n(F)}gU_2$.
\end{prop}

In the proof we use the following lemma.

\begin{lem}\label{lem:O_C-basis}
 Let $(C,C^+)$ be a complete affinoid field over $(\C_p,\mathcal{O}_{\C_p})$. Let $x=(x_1,\ldots,x_n)$ be an element of $(C^+)^n$
 such that $x_i\in (C^+)^\times$ for some $1\le i\le n$. Then, there exist $g\in \GL_n(\mathcal{O}_F)$ and integers $s$, $s'$
 with $1\le s\le s'\le n$ such that $gx=y=(y_1,\ldots,y_n)$ satisfies the following:
 \begin{itemize}
  \item $y_{s+1},\ldots,y_{s'}$ lie in $\varpi C^+$, and the images of $y_1,\ldots,y_s,\varpi^{-1}y_{s+1},\ldots,\varpi^{-1}y_{s'}$ in
	$C^+/\varpi C^+$ are linearly independent over $k_F$, and
  \item $y_{s'+1},\ldots,y_n$ lie in $\varpi^2C^+$.
 \end{itemize}
\end{lem}

\begin{prf}
 We write $W$ for the image of $\langle x_1,\ldots,x_n\rangle_{\mathcal{O}_F}$ in $C^+/\varpi C^+$ and put $s=\dim_{k_F}W$.
 By the assumption we have $1\le s\le n$. By permuting $x_1,\ldots,x_n$, we may assume that the images of $x_1,\ldots,x_s$
 in $C^+/\varpi C^+$ form a basis of $W$. Further, by Gaussian elimination, we may assume that
 $x_{s+1},\ldots,x_n\in\varpi C^+$. We write $W'$ for the image of $\langle\varpi^{-1}x_{s+1},\ldots,\varpi^{-1}x_n\rangle_{\mathcal{O}_F}$
 in $C^+/\varpi C^+$ and put $s'=\dim_{k_F}(W+W')$. By permuting $x_{s+1},\ldots,x_n$, we may assume that the images of
 $x_1,\ldots,x_s,\varpi^{-1}x_{s+1},\ldots,\varpi^{-1}x_{s'}$ in $C^+/\varpi C^+$ form a basis of $W+W'$.
 By Gaussian elimination again, we may assume that $x_{s'+1},\ldots,x_n$ lie in $\varpi^2C^+$. This concludes the proof.
\end{prf}

\begin{prf}[of Proposition \ref{prop:Drinfeld-case}]
 Let $(C,C^+)$ be a complete affinoid field over $(\C_p,\mathcal{O}_{\C_p})$ and $L\in \mathcal{F}^{\mathit{wa}}(C,C^+)$.
 We show that there exists $g\in \GL_n(F)$ such that $gL\in U_2(C,C^+)$. By Proposition \ref{prop:Drinfeld-tubular-nbhd} (iii),
 there exists an integer $m\ge 1$ such that $L\in U_m(C,C^+)$. 
 We have only to find $g\in \GL_n(F)$ satisfying $gL\in U_{m-1}(C,C^+)$ under the assumption $m\ge 3$.
 We set $I_L=\{W\in \P^{n-1}(F)\mid L\notin U_{m-1,W}(C,C^+)\}$. If $I_L=\varnothing$, then $L\in U_{m-1}(C,C^+)$, so we may assume $I_L\neq \varnothing$.
 For $W\in I_L$, we can find $g\in \GL_n(\mathcal{O}_F)$ such that $gW=\langle e_n\rangle_F$.
 By Proposition \ref{prop:Drinfeld-tubular-nbhd} (iv), we have ${}^t\!g^{-1}L\notin {}^t\!g^{-1}U_{m-1,W}(C,C^+)=U_{m-1,gW}(C,C^+)$.
 Therefore, replacing $L$ with ${}^t\!g^{-1}L$, we may assume that $W=\langle e_n\rangle_F$. Let $x=(x_1,\ldots,x_n)$ be a $C^+$-basis of $\Lambda_L$.
 By $L\notin U_{m-1,\langle e_n\rangle_F}(C,C^+)$ we have $\lvert x_n\rvert<\lvert\varpi^{m-1}\rvert\le \lvert\varpi^2\rvert$.
 By applying Lemma \ref{lem:O_C-basis} to $(x_1,\ldots,x_{n-1})$, after replacing $L$ with $\diag(h,1)L$ for a suitable $h\in\GL_{n-1}(\mathcal{O}_F)$,
 we may assume that $(x_1,\ldots,x_{n-1})$ satisfies the following: for some integers $s$, $s'$ with $1\le s\le s'\le n-1$, 
 \begin{itemize}
  \item $x_{s+1},\ldots,x_{s'}$ lie in $\varpi C^+$, and the images of $x_1,\ldots,x_s,\varpi^{-1}x_{s+1},\ldots,\varpi^{-1}x_{s'}$ in
	$C^+/\varpi C^+$ are linearly independent over $k_F$, and
  \item $x_{s'+1},\ldots,x_{n-1}$ lie in $\varpi^2 C^+$.
 \end{itemize}
 We put $t=\diag(\underbrace{1,\ldots,1}_{s},\underbrace{\varpi,\ldots,\varpi}_{n-s})\in \GL_n(F)$. As 
 \[
  t^{-1}x=(x_1,\ldots,x_s,\varpi^{-1}x_{s+1},\ldots,\varpi^{-1}x_n),
 \]
 we have $\Lambda_{t^{-1}L}=\langle t^{-1}x\rangle_{C^+}$.
 We prove that for every $W\in\P^{n-1}(F)$ we have $t^{-1}L\in U_{m-1,W}(C,C^+)$, that is,
 $\llparenthesis t^{-1}L,W\rrparenthesis\ge \lvert\varpi^{m-1}\rvert$.
 Let $a=(a_1,\ldots,a_n)$ be an $\mathcal{O}_F$-basis of $\Lambda_W$.

 First suppose that $a_{i_0}\in\mathcal{O}_F^\times$ for some $1\le i_0\le s'$. 
 We have $(t^{-1}x,a)\in \sum_{i=1}^sa_ix_i+\sum_{i=s+1}^{s'}a_i(\varpi^{-1}x_i)+\varpi C^+$
 (note that $\lvert\varpi^{-1}x_i\rvert\le \lvert\varpi\rvert$ for every $s'+1\le i\le n$).
 Since the images of $x_1,\ldots,x_s,\varpi^{-1}x_{s+1},\ldots,\varpi^{-1}x_{s'}$ in $C^+/\varpi C^+$
 are linearly independent over $k_F$, we have $(t^{-1}x,a)\in C^+\smallsetminus\varpi C^+$, in other words,
 $\llparenthesis t^{-1}x,a\rrparenthesis>\lvert\varpi\rvert$.
 Therefore we have $\llparenthesis t^{-1}L,W\rrparenthesis>\lvert\varpi\rvert\ge \lvert\varpi^{m-1}\rvert$, as desired.

 Next we consider the case where $a_i\in\mathfrak{p}_F$ for every $1\le i\le s'$.
 In this case we have $\varpi t^{-1}a=(\varpi a_1,\ldots,\varpi a_s,a_{s+1},\ldots,a_n)\in \Lambda_{t^{-1}W}$, and
 one of $a_{s'+1},\ldots,a_n$ lies in $\mathcal{O}_F^\times$. Therefore we have $\Lambda_{t^{-1}W}=\langle \varpi t^{-1}a\rangle_{\mathcal{O}_F}$.
 Hence
 \[
  \llparenthesis t^{-1}L,W\rrparenthesis=\llparenthesis t^{-1}x,a\rrparenthesis=\llparenthesis x,t^{-1}a\rrparenthesis
 =\lvert\varpi^{-1}\rvert\cdot \llparenthesis x,\varpi t^{-1}a\rrparenthesis=\lvert\varpi^{-1}\rvert\cdot \llparenthesis L,t^{-1}W\rrparenthesis\ge \lvert\varpi^{m-1}\rvert,
 \]
 as desired. In the last inequality we use $L\in U_m(C,C^+)$.
\end{prf}

\subsection{The orthogonal case}\label{subsec:orthogonal}
Next we consider the orthogonal case.

\begin{defn}
 Let $r\ge 2$ be an integer. We define a quadratic form $q_r\colon \mathcal{O}_F^{2r}\to \mathcal{O}_F$ by
 $q_r(x_1,\ldots,x_{2r})=\sum_{i=1}^rx_ix_{2r+1-i}$. Let $(\Lambda_{\mathrm{an}},q_{\mathrm{an}})$ be one of the following quadratic spaces
 over $\mathcal{O}_F$:
 \begin{itemize}
  \item $\Lambda_{\mathrm{an}}=0$ and $q_{\mathrm{an}}=0$; in this case we assume $r\ge 3$.
  \item $\Lambda_{\mathrm{an}}=\mathcal{O}_F$ and $q_{\mathrm{an}}(x)=x^2$.
  \item $\Lambda_{\mathrm{an}}=\mathcal{O}_E$ and $q_{\mathrm{an}}=N_{E/F}$ for a quadratic extension $E/F$.
 \end{itemize}
 We put $\Lambda=\mathcal{O}_F^{2r}\oplus \Lambda_{\mathrm{an}}$ and $q=q_r\oplus q_{\mathrm{an}}$.

 We set $V_{\mathrm{an}}=\Lambda_{\mathrm{an}}\otimes_{\mathcal{O}_F}F$ and $V=\Lambda\otimes_{\mathcal{O}_F}F$. We also write $q_{\mathrm{an}}$
 (resp.\ $q$) for the quadratic form on $V_{\mathrm{an}}$ (resp.\ $V$) induced from $q_{\mathrm{an}}$ (resp.\ $q$).
 We put $\Or(V,q)=\{g\in \GL(V)\mid q(gx)=q(x)\ (x\in V)\}$, $\SO(V,q)=\{g\in \Or(V,q)\mid \det g=1\}$, 
 $\Or(\Lambda,q)=\Or(V,q)\cap \GL(\Lambda)$ and $\SO(\Lambda,q)=\SO(V,q)\cap \GL(\Lambda)$.

 Let $b_q\colon \Lambda\times \Lambda\to \mathcal{O}_F$ denote the symmetric form defined by $b_q(x,y)=q(x+y)-q(x)-q(y)$.
 We also write $b_q$ for the base change $V\times V\to F$. We define $b_{q_{\mathrm{an}}}\colon \Lambda_{\mathrm{an}}\times \Lambda_{\mathrm{an}}\to \mathcal{O}_F$, $V_{\mathrm{an}}\times V_{\mathrm{an}}\to F$ in the same way.
\end{defn}

\begin{defn}
 An $F$-subspace $W$ of $V$ is said to be totally isotropic if $q\vert_W=0$.
 It is well-known that the dimension of any maximal totally isotropic subspace
 is equal to $r$. We write $\OGr(V,q)$ or simply $\OGr(V)$ for the set of maximal totally isotropic subspaces of $V$.
\end{defn}

The group $\SO(V,q)$ is quasi-split over $F$ and its type is as follows:
\begin{itemize}
 \item $D_r$ if $\Lambda_{\mathrm{an}}=0$.
 \item $B_r$ if $\Lambda_{\mathrm{an}}=\mathcal{O}_F$.
 \item ${}^2\!D_{r+1}$ if $\Lambda_{\mathrm{an}}=\mathcal{O}_E$ and $E/F$ is unramified.
 \item $C\text{-}B_r$ if $\Lambda_{\mathrm{an}}=\mathcal{O}_E$ and $E/F$ is totally ramified.
\end{itemize}
Therefore, to prove Proposition \ref{prop:case-by-case} (ii) we have only to consider the case $G=\SO(V,q)$ and
$\mu\colon \mathbb{G}_m\to G$; $z\mapsto \diag(z,\underbrace{1,\ldots,1}_{2r-2},z^{-1})\oplus \id_{V_{\mathrm{an}}}$
(see also Proposition \ref{prop:passing-to-adjoint}). Since $\mu^{-1}$ is conjugate to $\mu$, 
it suffices to consider $\mathcal{F}(G,1,\{\mu\})^{\mathit{wa}}$
in place of $\mathcal{F}(G,1,\{\mu^{-1}\})^{\mathit{wa}}$. For simplicity we write $\mathcal{F}$ and $\mathcal{F}^{\mathit{wa}}$ for 
$\mathcal{F}(G,\{\mu\})$ and $\mathcal{F}(G,1,\{\mu\})^{\mathit{wa}}$, respectively.
Since $P(\mu)$ is the stabilizer of $\langle e_1\rangle_F\subset V$ (here $e_1,\ldots,e_{2r}$ denote the standard basis of $F^{2r}$),
we have 
\[
 \mathcal{F}=\{L\in \P(V)_{\C_p}\mid \text{$L\subset V_{\C_p}$ is totally isotropic}\}.
\]
For $W\in \OGr(V)$, we put $\mathcal{F}_W=\mathcal{F}\cap \P(W)_{\C_p}$.
By using Proposition \ref{prop:ss=wa} (iii) (see also Remark \ref{rem:ss=wa-without-adjointness}), we can observe that
$\mathcal{F}^{\mathit{wa}}=\mathcal{F}^{\ad}\smallsetminus\bigcup_{W\in\OGr(V)}\mathcal{F}_W^{\ad}$.

Let $(C,C^+)$ be a complete affinoid field over $(\C_p,\mathcal{O}_{\C_p})$, and let $\lvert-\rvert$ denote the continuous valuation on $C$
corresponding to the closed point of $\Spa(C,C^+)$. For $x\in V_C$ and $a\in V$, we put $\llparenthesis x,a\rrparenthesis=\lvert b_q(x,a)\rvert$.
For $L\in \mathcal{F}(C)$ (resp.\ an $F$-subspace $W$ of $V$), we put $\Lambda_L=\Lambda_{C^+}\cap L$ (resp.\ $\Lambda_W=\Lambda\cap W$),
which is a free $C^+$-module (resp.\ $\mathcal{O}_F$-module) of rank $1$ (resp.\ $\dim_FW$).
For $L\in \mathcal{F}(C)$ and $W\in\OGr(V)$, we define
$\llparenthesis L,W^\perp\rrparenthesis=\max\{\llparenthesis x,a\rrparenthesis \mid x\in \Lambda_L,a\in \Lambda_{W^\perp}\}$,
where $W^\perp$ denotes the orthogonal complement of $W$ with respect to $b_q$.
Since $\Lambda_L$ and $\Lambda_{W^\perp}$ are finitely generated, the maximum always exists.
The following proposition can be proven in the same way as Proposition \ref{prop:Drinfeld-tubular-nbhd}.

\begin{prop}\label{prop:orthogonal-tubular-nbhd}
 \begin{enumerate}
  \item For $W\in\OGr(V)$ and an integer $m\ge 1$, there exists a unique quasi-compact open subset $U_{m,W}$ of $\mathcal{F}^{\ad}$
	satisfying the following: for every complete affinoid field $(C,C^+)$ over $(\C_p,\mathcal{O}_{\C_p})$, we have 
	\[
	 U_{m,W}(C,C^+)=\{L\in \mathcal{F}^{\ad}(C,C^+)\mid \llparenthesis L,W^\perp\rrparenthesis\ge \lvert\varpi^m\rvert\}.
	\]
  \item The subset $U_m=\bigcap_{W\in \OGr(V)}U_{m,W}$ of $\mathcal{F}^{\ad}$ is a quasi-compact open subset of $\mathcal{F}^{\mathit{wa}}$.
  \item We have $\mathcal{F}^{\mathit{wa}}=\bigcup_{m\ge 1}U_m$.
  \item For $g\in\Or(\Lambda,q)$, we have $gU_{m,W}=U_{m,gW}$ and $gU_m=U_m$.
 \end{enumerate}
\end{prop}

To show Proposition \ref{prop:case-by-case} (ii), it suffices to prove the following proposition.

\begin{prop}\label{prop:orthogonal-case}
 Let $\alpha$ be the non-negative integer defined as follows:
 \begin{itemize}
  \item $\alpha=0$ if $\Lambda_{\mathrm{an}}=0$.
  \item $\alpha=v_F(2)$ if $\Lambda_{\mathrm{an}}=\mathcal{O}_F$.
  \item $\alpha=0$ if $\Lambda_{\mathrm{an}}=\mathcal{O}_E$ and $E/F$ is unramified.
  \item $\alpha=\lceil v_F(\mathfrak{d}_{E/F})/2\rceil$ if $\Lambda_{\mathrm{an}}=\mathcal{O}_E$ and $E/F$ is ramified,
	where $\mathfrak{d}_{E/F}$ denotes the discriminant of $E/F$.
 \end{itemize}
 Then we have $\mathcal{F}^{\mathit{wa}}=\bigcup_{g\in\SO(V,q)}gU_{3+\alpha}$.
\end{prop}

We need some preparation for the proof of this proposition.

\begin{lem}\label{lem:max-isotropic-transitive}
 For every $W_1,W_2\in \OGr(V)$, there exists $g\in\Or(\Lambda,q)$ such that $g\Lambda_{W_1}=\Lambda_{W_2}$.
\end{lem}

\begin{prf}
 We allow the cases $r=0,1$ (and $r=2$ when $\Lambda_{\mathrm{an}}=0$), and prove the claim by induction on $r$.
 If $r=0$, the claim is clear since $W_1=W_2=0$. Assume $r\ge 1$ and that the claim holds for $r-1$.
 We may assume that $W_2=\langle e_1,\ldots,e_r\rangle_F$. We take an $\mathcal{O}_F$-basis $f_1,\ldots,f_r$ of $\Lambda_{W_1}$.
 By \cite[\S 9.1, Remark 2]{MR0195863}, we can find $g\in \Or(\Lambda,q)$ such that $g(f_1)=e_1$. Replacing $W_1$ with $gW_1$,
 we may assume $f_1=e_1$. We put $\Lambda'=\langle e_2,\ldots,e_{2r-1}\rangle_{\mathcal{O}_F}\oplus \Lambda_{\mathrm{an}}$ and 
 $V'=\langle e_2,\ldots,e_{2r-1}\rangle_F\oplus V_{\mathrm{an}}$. Then we have orthogonal decompositions
 $\Lambda=\langle e_1,e_{2r}\rangle_{\mathcal{O}_F}\stackrel{\perp}{\oplus}\Lambda'$ and $V=\langle e_1,e_{2r}\rangle_F\stackrel{\perp}{\oplus}V'$.
 We put $W_1'=W_1\cap V'$ and prove that $W_1=\langle e_1\rangle_F\oplus W_1'$ and $\Lambda_{W_1}=\langle e_1\rangle_{\mathcal{O}_F}\oplus \Lambda_{W_1'}$.
 First we have $\langle e_1\rangle_F\cap W'_1\subset \langle e_1\rangle_F\cap V'=0$. For every $x=(x_1,\ldots,x_{2r},x_{\mathrm{an}})\in W_1$,
 we have $x_{2r}=b_q(x,e_1)=0$, which implies that $x-x_1e_1\in W_1'$. Therefore we obtain $W_1=\langle e_1\rangle_F\oplus W'_1$.
 Further, if $x$ lies in $\Lambda_{W_1}$, we have $x_1,\ldots,x_{2r-1}\in\mathcal{O}_F$ and $x_{\mathrm{an}}\in \Lambda_{\mathrm{an}}$,
 which implies that $x-x_1e_1\in \Lambda_{W_1'}$. Hence we conclude that $\Lambda_{W_1}=\langle e_1\rangle_{\mathcal{O}_F}\oplus\Lambda_{W'_1}$.
 Note that $W_1'\in \OGr(V')$, as it is totally isotropic with dimension $r-1$.
 By applying the induction hypothesis to $W_1'$ and $W_2'=\langle e_2,\ldots,e_r\rangle_F$, 
 we can find an element $g'\in \Or(\Lambda',q_{r-1}\oplus q_{\mathrm{an}})$ such that $g'\Lambda_{W_1'}=\Lambda_{W_2'}$.
 The element $g=\id_{\langle e_1,e_{2r}\rangle_{\mathcal{O}_F}}\oplus g'\in \Or(\Lambda,q)$ satisfies $g\Lambda_{W_1}=\Lambda_{W_2}$, as desired.
\end{prf}

\begin{lem}\label{lem:matrix-SO}
 Let $W$ be an element of $\OGr(V,q)$. We take an $\mathcal{O}_F$-basis $y_1,\ldots,y_r$ of $\Lambda_W$ and write 
 $y_j=(y_{1j},\ldots,y_{2r,j},y_{\mathrm{an},j})$ for each $1\le j\le r$. Let $s$ be an integer satisfying $1\le s\le r$, and assume that 
 $y_{ij}\in\mathfrak{p}_F$ for every $2r+1-s\le i\le 2r$ and $1\le j\le r$.
 Then, the rank of the matrix $(y_{ij}\bmod\mathfrak{p}_F)_{1\le i\le s,1\le j\le r}$ is equal to $s$.
 In particular, we can replace $y_1,\ldots,y_r$ so that $(y_{ij})_{1\le i,j\le s}$ is the identity matrix and
 $y_{ij}=0$ for every $1\le i\le s$ and $s+1\le j\le r$.
\end{lem}

\begin{prf}
 We allow the cases $r=0,1$ (and $r=2$ when $\Lambda_{\mathrm{an}}=0$).
 We prove the claim by induction on $s\ge 1$. Assume that the claim holds for $s-1$.

 First assume that $y_{11},\ldots,y_{1r}\in\mathfrak{p}_F$. We put $\overline{\Lambda}_{\mathrm{an}}=\Lambda_{\mathrm{an}}\otimes_{\mathcal{O}_F}k_F$
 and write $\overline{\Lambda}_W$ for the image of $\Lambda_W\subset \Lambda=\mathcal{O}_F^{2r}\oplus\Lambda_{\mathrm{an}}$ in $k_F^{2r}\oplus \overline{\Lambda}_{\mathrm{an}}$. By the assumption, $\overline{\Lambda}_W$ is contained in $0\oplus k_F^{2r-2}\oplus 0\oplus \overline{\Lambda}_{\mathrm{an}}$.
 Let $\overline{q}_{\mathrm{an}}$ (resp.\ $\overline{q}$) denote the quadratic form on $\overline{\Lambda}_{\mathrm{an}}$
 (resp.\ $k_F^{2r-2}\oplus \overline{\Lambda}_{\mathrm{an}}=0\oplus k_F^{2r-2}\oplus 0\oplus \overline{\Lambda}_{\mathrm{an}}$)
 induced from $q_{\mathrm{an}}$ (resp.\ $q$). Then, $\overline{\Lambda}_W$ gives a totally isotropic subspace of 
 $(k_F^{2r-2}\oplus \overline{\Lambda}_{\mathrm{an}},\overline{q})$. If $(\Lambda_{\mathrm{an}},q_{\mathrm{an}})$ is not $(\mathcal{O}_E,N_{E/F})$ 
 for a ramified quadratic extension $E$ of $F$, then $\overline{q}_{\mathrm{an}}$ is regular in the sense of \cite[p.~42]{MR2427530},
 and all maximal totally isotropic subspaces of $(k_F^{2r-2}\oplus \overline{\Lambda}_{\mathrm{an}},\overline{q})$ are $(r-1)$-dimensional
 (see \cite[Corollary 8.12]{MR2427530}).
 Therefore we have $r=\dim \overline{\Lambda}_W\le r-1$, which is a contradiction. If $(\Lambda_{\mathrm{an}},q_{\mathrm{an}})=(\mathcal{O}_E,N_{E/F})$
 for a ramified quadratic extension $E$ of $F$, the quadratic radical 
 $\rad(\overline{q}_{\mathrm{an}})=\{x\in \overline{\Lambda}_{\mathrm{an}}\mid b_{\overline{q}_{\mathrm{an}}}(x,y)=0\ (y\in \overline{\Lambda}_{\mathrm{an}}), \overline{q}_{\mathrm{an}}(x)=0\}$ of $\overline{q}_{\mathrm{an}}$ is equal to $\mathfrak{p}_E/\mathfrak{p}_F\mathcal{O}_E$,
 and all maximal totally isotropic subspaces of $(k_F^{2r-2}\oplus \overline{\Lambda}_{\mathrm{an}})/(0\oplus \rad(\overline{q}_{\mathrm{an}}))$
 are $(r-1)$-dimensional. Since $\dim \overline{\Lambda}_W=r$, $\overline{\Lambda}_W$ must contain $\rad(\overline{q}_{\mathrm{an}})$;
 otherwise, the image of $\overline{\Lambda}_W$ in $(k_F^{2r-2}\oplus \overline{\Lambda}_{\mathrm{an}})/(0\oplus \rad(\overline{q}_{\mathrm{an}}))$
 gives an $r$-dimensional totally isotropic subspace. Therefore, there exists $x=(x_1,\ldots,x_{2r},x_{\mathrm{an}})\in \Lambda_W\cap (\mathfrak{p}_F^{2r}\oplus \mathfrak{p}_E)$ such that $v_E(x_{\mathrm{an}})=1$. Since $W$ is totally isotropic, we have
 $0=q(x)=\sum_{i=1}^rx_ix_{2r+1-i}+N_{E/F}(x_{\mathrm{an}})$.
 Therefore $N_{E/F}(x_{\mathrm{an}})=-\sum_{i=1}^rx_ix_{2r+1-i}\in \mathfrak{p}_F^2$, which contradicts $v_E(x_{\mathrm{an}})=1$.

 Hence we conclude that one of $y_{11},\ldots,y_{1r}$ belongs to $\mathcal{O}_F^\times$. By permuting $y_1,\ldots,y_r$, we may assume that
 $y_{11}\in\mathcal{O}_F^\times$. By replacing $y_1,\ldots,y_r$ further, we may assume that $y_{11}=1$ and $y_{12}=\cdots=y_{1r}=0$.
 Now we can apply the induction hypothesis to $\langle (y_{2j},\ldots,y_{2r-1,j},y_{\mathrm{an},j})\mid 2\le j\le r\rangle_F\in\OGr(F^{2r-2}\oplus V_{\mathrm{an}},q_{r-1}\oplus q_{\mathrm{an}})$ to conclude the claim.
\end{prf}

\begin{lem}\label{lem:integral-pairing}
 Let $\alpha$ be the integer defined in Proposition \ref{prop:orthogonal-case}.
 For a complete affinoid field $(C,C^+)$ over $(\C_p,\mathcal{O}_{\C_p})$, we have
 $\{x\in V_{\mathrm{an},C}\mid \lvert b_{q_{\mathrm{an}}}(x,y)\rvert\le \lvert\varpi^\alpha\rvert\ (y\in\Lambda_{\mathrm{an}})\}\subset \Lambda_{\mathrm{an},C^+}$.
\end{lem}

\begin{prf}
 If $\Lambda_{\mathrm{an}}=0$, the claim is clear. If $\Lambda_{\mathrm{an}}=\mathcal{O}_F$, we have
 \[
  \{x\in V_{\mathrm{an},C}\mid \lvert b_{q_{\mathrm{an}}}(x,y)\rvert\le \lvert\varpi^\alpha\rvert\ (y\in\Lambda_{\mathrm{an}})\}
 =\{x\in C\mid \lvert 2xy\rvert\le \lvert 2\rvert\ (y\in\mathcal{O}_F)\}=C^+,
 \]
 as desired. We consider the case $\Lambda_{\mathrm{an}}=\mathcal{O}_E$. Let $c$ denote the non-trivial element of $\Gal(E/F)$.
 Then we have $b_{q_{\mathrm{an}}}(x,y)=(x+y)c(x+y)-xc(x)-yc(y)=xc(y)+c(x)y=\tr_{E/F}(xc(y))$ for $x,y\in E$. Note that we have
 $\{x\in E\mid \tr_{E/F}(xy)\in\mathcal{O}_F\ (y\in \mathcal{O}_E)\}=\mathfrak{D}_{E/F}^{-1}\subset\mathfrak{p}_F^{-\alpha}\mathcal{O}_E$,
 where $\mathfrak{D}_{E/F}$ denotes the different for $E/F$. Indeed, if $E/F$ is unramified,
 we have $\mathfrak{D}_{E/F}^{-1}=\mathcal{O}_E=\mathfrak{p}_F^{-\alpha}\mathcal{O}_E$.
 If $E/F$ is ramified, we have $v_E(\mathfrak{D}_{E/F})=v_F(N_{E/F}(\mathfrak{D}_{E/F}))=v_F(\mathfrak{d}_{E/F})$, hence
 $\mathfrak{D}_{E/F}^{-1}=\mathfrak{p}_E^{-v_F(\mathfrak{d}_{E/F})}\subset \mathfrak{p}_E^{-2\alpha}=\mathfrak{p}_F^{-\alpha}\mathcal{O}_E$.
 We fix an $\mathcal{O}_F$-basis $f_1,f_2$ of $\mathcal{O}_E$, and let $f_1^\vee,f_2^\vee\in E$ be the dual basis of it with respect to
 $b_{q_{\mathrm{an}}}$. By the above inclusion, we have $f_1^\vee,f_2^\vee\in \mathfrak{p}_F^{-\alpha}\mathcal{O}_E=\varpi^{-\alpha}\Lambda_{\mathrm{an}}$.
 For $x=x_1f_1^\vee+x_2f_2^\vee\in V_{\mathrm{an},C}$ $(x_1,x_2\in C)$ and $y=a_1f_1+a_2f_2\in \Lambda_{\mathrm{an}}$ $(a_1,a_2\in \mathcal{O}_F)$,
 we have $\lvert b_q(x,y)\rvert=\lvert a_1x_1+a_2x_2\rvert$. It is smaller than or equal to $\lvert\varpi^\alpha\rvert$ for every $a_1,a_2\in\mathcal{O}_F$
 if and only if $\lvert x_1\rvert, \lvert x_2\rvert\le \lvert\varpi^\alpha\rvert$, in other words, $x_1,x_2\in \varpi^\alpha C^+$.
 Together with $f_1^\vee,f_2^\vee\in \varpi^{-\alpha}\Lambda_{\mathrm{an}}$, we conclude that $x=x_1f_1^\vee+x_2f_2^\vee\in \Lambda_{\mathrm{an},C^+}$,
 as desired.
\end{prf}

\begin{prf}[of Proposition \ref{prop:orthogonal-case}]
 Let $(C,C^+)$ be a complete affinoid field over $(\C_p,\mathcal{O}_{\C_p})$ and $L\in \mathcal{F}^{\mathit{wa}}(C,C^+)$.
 We show that there exists $g\in \SO(V,q)$ such that $gL\in U_{3+\alpha}(C,C^+)$. By Proposition \ref{prop:orthogonal-tubular-nbhd} (iii),
 there exists an integer $m\ge 1$ such that $L\in U_m(C,C^+)$. 
 We have only to find $g\in \SO(V,q)$ satisfying $gL\in U_{m-1}(C,C^+)$ under the assumption $m\ge 4+\alpha$.
 We set $I_L=\{W\in \OGr(V)\mid L\notin U_{m-1,W}(C,C^+)\}$. If $I_L=\varnothing$, then $L\in U_{m-1}(C,C^+)$, so we may assume $I_L\neq \varnothing$.
 For $W\in I_L$, we can find $g\in \Or(\Lambda,q)$ such that $gW=\langle e_1,\ldots,e_r\rangle_F$ by Lemma \ref{lem:max-isotropic-transitive}.
 By Proposition \ref{prop:orthogonal-tubular-nbhd} (iv), we have $gL\notin gU_{m-1,W}(C,C^+)=U_{m-1,gW}(C,C^+)$.
 
 First we consider the case where $g\in \SO(\Lambda,q)$. Replacing $L$ with $gL$, we may assume that $W=\langle e_1,\ldots,e_r\rangle_F$.
 Then we have $\Lambda_{W^\perp}=\langle e_1,\ldots,e_r\rangle_{\mathcal{O}_F}\oplus \Lambda_{\mathrm{an}}$.
 Let $x=(x_1,\ldots,x_{2r},x_{\mathrm{an}})$ be a $C^+$-basis of $\Lambda_L$.
 By $L\notin U_{m-1,\langle e_1,\ldots,e_r\rangle_F}(C,C^+)$ we have $\lvert x_i\rvert<\lvert\varpi^{m-1}\rvert\le \lvert\varpi^{3+\alpha}\rvert\le \lvert\varpi^3\rvert$ for every $r+1\le i\le 2r$. Further, Lemma \ref{lem:integral-pairing} tells us that $x_{\mathrm{an}}\in \varpi^3\Lambda_{\mathrm{an},C^+}$.
 Therefore, one of $x_1,\ldots,x_r$ belongs to $(C^+)^\times$.
 Let $J_r$ be the $r\times r$ matrix $(\delta_{i,r+1-j})_{i,j}$, where $\delta_{i,r+1-j}$ denotes the Kronecker delta.
 By applying Lemma \ref{lem:O_C-basis} to $(x_1,\ldots,x_r)$, after replacing $L$ with its translate by
 $(\diag(h,J_r{}^t\!h^{-1}J_r)\oplus \id_{\Lambda_{\mathrm{an}}})\in \SO(\Lambda,q)$ for a suitable $h\in\GL_r(\mathcal{O}_F)$,
 we may assume that $(x_1,\ldots,x_r)$ satisfies the following: for some integers $s$, $s'$ with $1\le s\le s'\le r$, 
 \begin{itemize}
  \item $x_{s+1},\ldots,x_{s'}$ lie in $\varpi C^+$, and the images of $x_1,\ldots,x_s,\varpi^{-1}x_{s+1},\ldots,\varpi^{-1}x_{s'}$ in
	$C^+/\varpi C^+$ are linearly independent over $k_F$, and
  \item $x_{s'+1},\ldots,x_r$ lie in $\varpi^2 C^+$.
 \end{itemize}
 Note that the new $L$ still belongs to $U_m(C,C^+)\smallsetminus U_{m-1,\langle e_1,\ldots,e_r\rangle_F}(C,C^+)$ by
 Proposition \ref{prop:orthogonal-tubular-nbhd} (iv).

 We put $t=\diag(\underbrace{\varpi^{-1},\ldots,\varpi^{-1}}_{s},\underbrace{1,\ldots,1}_{2r-2s},\underbrace{\varpi,\ldots,\varpi}_s)\oplus \id_{V_{\mathrm{an}}}\in \SO(V,q)$. As 
 \[
  \varpi^{-1}t^{-1}x=(x_1,\ldots,x_s,\varpi^{-1}x_{s+1},\ldots,\varpi^{-1}x_{2r-s},\varpi^{-2}x_{2r+1-s},\ldots,\varpi^{-2}x_{2r},\varpi^{-1}x_{\mathrm{an}}),
 \]
 we have $\Lambda_{t^{-1}L}=\langle \varpi^{-1}t^{-1}x\rangle_{C^+}$ (here we use $\lvert x_i\rvert\le \lvert\varpi^2\rvert$ for $r+1\le i\le 2r$
 and $x_{\mathrm{an}}\in \varpi\Lambda_{\mathrm{an},C^+}$).
 We prove that for every $W\in\OGr(V)$ we have $t^{-1}L\in U_{m-1,W}(C,C^+)$, that is, 
 $\llparenthesis t^{-1}L,W^\perp\rrparenthesis\ge \lvert\varpi^{m-1}\rvert$.
 Let $y_1,\ldots,y_r$ be an $\mathcal{O}_F$-basis of $\Lambda_W$. Further, we take $y_{r+1},\ldots,y_{r'}\in \Lambda_{W^\perp}$ so that
 $y_1,\ldots,y_{r'}$ form a basis of $\Lambda_{W^\perp}$ (in particular $r'=r+\dim V_{\mathrm{an}}$).
 We write $y_j=(y_{1j},\ldots,y_{2r,j},y_{\mathrm{an},j})$ for each $1\le j\le r'$. 

 First suppose that $y_{i_0j_0}\in \mathcal{O}_F^\times$ for some $2r+1-s\le i_0\le 2r$ and $1\le j_0\le r'$.
 Then we have $b_q(\varpi^{-1}t^{-1}x,y_{j_0})\in \sum_{i=1}^s y_{2r+1-i,j_0}x_i+\sum_{i=s+1}^{s'} y_{2r+1-i,j_0}\varpi^{-1}x_i+\varpi C^+$;
 here we use $\lvert\varpi^{-1}x_i\rvert\le \lvert\varpi\rvert$ for $s'+1\le i\le 2r-s$, $\lvert\varpi^{-2}x_i\rvert\le \lvert\varpi\rvert$
 for $2r+1-s\le i\le 2r$ and $\varpi^{-1}x_{\mathrm{an}}\in \varpi\Lambda_{\mathrm{an},C^+}$.
 Since the images of $x_1,\ldots,x_s,\varpi^{-1}x_{s+1},\ldots,\varpi^{-1}x_{s'}$ in $C^+/\varpi C^+$ are linearly independent over $k_F$,
 we have $b_q(\varpi^{-1}t^{-1}x,y_{j_0})\in C^+\smallsetminus \varpi C^+$, in other words,
 $\lvert b_q(\varpi^{-1}t^{-1}x,y_{j_0})\rvert>\lvert\varpi\rvert$.
 Therefore we have $\llparenthesis t^{-1}L,W^\perp\rrparenthesis>\lvert\varpi\rvert>\lvert\varpi^{m-1}\rvert$, as desired.

 Next we consider the case where $y_{ij}\in\mathfrak{p}_F$ for every $2r+1-s\le i\le 2r$ and $1\le j\le r'$.
 By Lemma \ref{lem:matrix-SO}, we can replace $y_1,\ldots,y_r$ so that $(y_{ij})_{1\le i,j\le s}$ is the identity matrix and
 $y_{ij}=0$ for every $1\le i\le s$ and $s+1\le j\le r$. Further, we can also replace $y_{r+1},\ldots,y_{r'}$ so that $y_{ij}=0$ for every
 $1\le i\le s$ and $r+1\le j\le r'$. We prove $\Lambda_{tW^\perp}=\langle \varpi ty_1,\ldots,\varpi ty_s,ty_{s+1},\ldots,ty_{r'}\rangle_{\mathcal{O}_F}$.
 Clearly we have $\Lambda_{tW^\perp}\supset\langle \varpi ty_1,\ldots,\varpi ty_s,ty_{s+1},\ldots,ty_{r'}\rangle_{\mathcal{O}_F}$.
 Let us prove the reverse inclusion.
 For an arbitrary element $z\in \Lambda_{tW^\perp}$, there exist $a_1,\ldots,a_{r'}\in F$ such that $z=a_1(\varpi ty_1)+\cdots+a_s(\varpi ty_s)+a_{s+1}(ty_{s+1})+\cdots+a_{r'}(ty_{r'})$. It suffices to show $a_1,\ldots,a_{r'}\in \mathcal{O}_F$. 
 As $a_1(\varpi ty_1)+\cdots+a_s(\varpi ty_s)+a_{s+1}(ty_{s+1})+\cdots+a_{r'}(ty_{r'})\in \Lambda=\mathcal{O}_F^{2r}\oplus \Lambda_{\mathrm{an}}$,
 by looking at the first $s$ entries, we obtain $a_1,\ldots,a_s\in\mathcal{O}_F$. 
 Therefore we have $a_{s+1}(ty_{s+1})+\cdots+a_{r'}(ty_{r'})\in \Lambda$. Let $d\ge 0$ be the minimal non-negative integer
 such that $\varpi^da_{s+1},\ldots,\varpi^da_{r'}\in \mathcal{O}_F$. If $d\ge 1$, then we have
 $\varpi^da_{s+1}(ty_{s+1})+\cdots+\varpi^da_{r'}(ty_{r'})\in \varpi^d\Lambda\subset \varpi\Lambda$.
 Therefore we obtain the equality $\overline{\varpi^da_{s+1}}\cdot \overline{y_{s+1}}+\cdots+\overline{\varpi^da_{r'}}\cdot \overline{y_{r'}}=0$
 in $\Lambda\otimes_{\mathcal{O}_F}k_F$, where $\overline{a}$ (resp.\ $\overline{y}$) denotes the image of $a\in \mathcal{O}_F$ (resp.\ $y\in \Lambda$)
 in $k_F$ (resp.\ $\Lambda\otimes_{\mathcal{O}_F}k_F$). Note that we have $\overline{ty_j}=\overline{y_j}$ for every $s+1\le j\le r'$,
 as $y_{1j}=\cdots=y_{sj}=0$ and $y_{2r+1-s,j},\ldots,y_{2r,j}\in \mathfrak{p}_F$. Since $\overline{y_{s+1}},\ldots,\overline{y_{r'}}\in\Lambda\otimes_{\mathcal{O}_F}k_F$ are linearly independent over $k_F$, we have $\overline{\varpi^da_{s+1}}=\cdots=\overline{\varpi^da_{r'}}=0$, in other words,
 $\varpi^da_{s+1},\ldots,\varpi^da_{r'}\in \mathfrak{p}_F$. Therefore $\varpi^{d-1}a_{s+1},\ldots,\varpi^{d-1}a_{r'}$ lie in $\mathcal{O}_F$,
 which contradicts the minimality of $d$. Hence we have $d=0$ and $a_{s+1},\ldots,a_{r'}\in\mathcal{O}_F$.
 This completes the proof of $\Lambda_{tW^\perp}=\langle \varpi ty_1,\ldots,\varpi ty_s,ty_{s+1},\ldots,ty_{r'}\rangle_{\mathcal{O}_F}$.
 Therefore we have
 \begin{align*}
  &\llparenthesis t^{-1}L,W^\perp\rrparenthesis=\max\{\lvert b_q(\varpi^{-1}t^{-1}x,y_1)\rvert,\ldots,\lvert b_q(\varpi^{-1}t^{-1}x,y_{r'})\rvert\}\\
  &\qquad =\max\{\lvert\varpi^{-2}\rvert\cdot \lvert b_q(x,\varpi ty_1)\rvert,\ldots,\lvert \varpi^{-2}\rvert\cdot \lvert b_q(x,\varpi ty_s)\rvert,\\
  &\qquad\qquad\qquad\qquad\qquad\qquad \lvert \varpi^{-1}\rvert\cdot \lvert b_q(x,ty_{s+1})\rvert,\ldots,\lvert \varpi^{-1}\rvert\cdot \lvert b_q(x,ty_{r'})\rvert\}\\
  &\qquad\ge \lvert\varpi^{-1}\rvert\cdot\max\{\lvert b_q(x,\varpi ty_1)\rvert,\ldots,\lvert b_q(x,\varpi ty_s)\rvert,\lvert b_q(x,ty_{s+1})\rvert,\ldots,\lvert b_q(x,ty_{r'})\rvert\}\\
  &\qquad=\lvert\varpi^{-1}\rvert\cdot\llparenthesis L,tW^\perp\rrparenthesis=\lvert\varpi^{-1}\rvert\cdot\llparenthesis L,(tW)^\perp\rrparenthesis\ge \lvert\varpi^{m-1}\rvert,
 \end{align*}
 as desired. In the last inequality we use $L\in U_{m,tW}(C,C^+)$.
 Now we conclude that $t^{-1}L\in U_{m-1}(C,C^+)$.
 
 Finally, we consider the case where $g\in \Or(\Lambda,q)\smallsetminus \SO(\Lambda,q)$. By applying the argument above to
 $gL\in U_m(C,C^+)\smallsetminus U_{m-1,\langle e_1,\ldots,e_r\rangle_F}(C,C^+)$, we can find $h\in \SO(V,q)$ such that $hgL\in U_{m-1}(C,C^+)$.
 By Proposition \ref{prop:orthogonal-tubular-nbhd} (iv), we have $g^{-1}hgL\in U_{m-1}(C,C^+)$.
 This completes the proof, since $g^{-1}hg\in \SO(V,q)$.
\end{prf}

\subsection{The unitary case}\label{subsec:unitary}
The unitary case can be treated in a similar way as in the orthogonal case.

\begin{defn}
 Let $E$ be a quadratic extension of $F$ contained in $\overline{F}$ and $c$ the non-trivial element of $\Gal(E/F)$. Let $n\ge 3$ be an integer.
 We put $\Lambda=\mathcal{O}_E^n$ and define a Hermitian form $h_n\colon \Lambda\times \Lambda\to \mathcal{O}_E$ by 
 $h_n((x_i),(y_i))=\sum_{i=1}^nx_ic(y_{n+1-i})$.

 We set $V=E^n$. We also write $h_n$ for the Hermitian form on $V$ induced from $h_n$.
 We put $\U(V)=\{g\in \GL(V)\mid h_n(gx,gy)=h_n(x,y)\ (x,y\in V)\}$ and $\U(\Lambda)=\U(V)\cap \GL(\Lambda)$. 
\end{defn}

\begin{defn}
 An $E$-subspace $W$ of $V$ is said to be totally isotropic if $h_n\vert_{W\times W}=0$.
 It is well-known that the dimension of any maximal totally isotropic subspace
 is equal to $\lfloor\frac{n}{2}\rfloor$. We write $\UGr(V)$ for the set of maximal totally isotropic subspaces of $V$.
\end{defn}

The group $\U(V)$ is quasi-split over $F$ and its type is as follows:
\begin{itemize}
 \item ${}^2\!A'_{n-1}$ if $E/F$ is unramified.
 \item $C\text{-}BC_r$ if $E/F$ is ramified and $n=2r+1$ is odd.
 \item $B\text{-}C_r$ if $E/F$ is ramified and $n=2r$ is even.
\end{itemize}
Therefore, to prove Proposition \ref{prop:case-by-case} (iii) we have only to consider the case $G=\U(V)$ and
$\mu\colon \mathbb{G}_m\to G_{\overline{F}}\cong \GL_{n,\overline{F}}$; $z\mapsto \diag(z,1,\ldots,1)$
(see also Proposition \ref{prop:passing-to-adjoint}). By changing the isomorphism $G_{\overline{F}}\cong \GL_{n,\overline{F}}$,
it suffices to consider $\mathcal{F}(G,1,\{\mu\})^{\mathit{wa}}$
in place of $\mathcal{F}(G,1,\{\mu^{-1}\})^{\mathit{wa}}$. For simplicity we write $\mathcal{F}$ and $\mathcal{F}^{\mathit{wa}}$ for 
$\mathcal{F}(G,\{\mu\})$ and $\mathcal{F}(G,1,\{\mu\})^{\mathit{wa}}$, respectively.
Since $P(\mu)$ is the stabilizer of $\langle e_1\rangle_E\subset V$ (here $e_1,\ldots,e_n$ denote the standard basis of $V=E^n$),
we have $\mathcal{F}=\P(V)_{\C_p}$.
For $W\in \UGr(V)$, we put $\mathcal{F}_W=\P(W)_{\C_p}$.
By using Proposition \ref{prop:ss=wa} (iii) (see also Remark \ref{rem:ss=wa-without-adjointness}), we can observe that
$\mathcal{F}^{\mathit{wa}}=\mathcal{F}^{\ad}\smallsetminus\bigcup_{W\in\UGr(V)}\mathcal{F}_W^{\ad}$.

Let $(C,C^+)$ be a complete affinoid field over $(\C_p,\mathcal{O}_{\C_p})$, and let $\lvert-\rvert$ denote the continuous valuation on $C$
corresponding to the closed point of $\Spa(C,C^+)$. We write $V_C=V\otimes_EC$ and $\Lambda_{C^+}=\Lambda\otimes_{\mathcal{O}_E}C^+$.
For $x\in V_C$ and $a\in V$, we put $\llparenthesis x,a\rrparenthesis=\lvert h_n(x,a)\rvert$.
For $L\in \mathcal{F}(C)$ (resp.\ an $E$-subspace $W$ of $V$), we put $\Lambda_L=\Lambda_{C^+}\cap L$ (resp.\ $\Lambda_W=\Lambda\cap W$),
which is a free $C^+$-module (resp.\ $\mathcal{O}_E$-module) of rank $1$ (resp.\ $\dim_EW$).
For $L\in \mathcal{F}(C)$ and $W\in\UGr(V)$, we define
$\llparenthesis L,W^\perp\rrparenthesis=\max\{\llparenthesis x,a\rrparenthesis \mid x\in \Lambda_L,a\in \Lambda_{W^\perp}\}$,
where $W^\perp$ denotes the orthogonal complement of $W$ with respect to $h_n$.
Since $\Lambda_L$ and $\Lambda_{W^\perp}$ are finitely generated, the maximum always exists.
We fix a uniformizer $\varpi_E$ of $\mathcal{O}_E$.
The following proposition can be proven in the same way as Proposition \ref{prop:Drinfeld-tubular-nbhd}.

\begin{prop}\label{prop:unitary-tubular-nbhd}
 \begin{enumerate}
  \item For $W\in\UGr(V)$ and an integer $m\ge 1$, there exists a unique quasi-compact open subset $U_{m,W}$ of $\mathcal{F}^{\ad}$
	satisfying the following: for every complete affinoid field $(C,C^+)$ over $(\C_p,\mathcal{O}_{\C_p})$, we have 
	\[
	 U_{m,W}(C,C^+)=\{L\in \mathcal{F}^{\ad}(C,C^+)\mid \llparenthesis L,W^\perp\rrparenthesis\ge \lvert\varpi_E^m\rvert\}.
	\]
  \item The subset $U_m=\bigcap_{W\in \UGr(V)}U_{m,W}$ of $\mathcal{F}^{\ad}$ is a quasi-compact open subset of $\mathcal{F}^{\mathit{wa}}$.
  \item We have $\mathcal{F}^{\mathit{wa}}=\bigcup_{m\ge 1}U_m$.
  \item For $g\in\U(\Lambda)$, we have $gU_{m,W}=U_{m,gW}$ and $gU_m=U_m$.
 \end{enumerate}
\end{prop}

To show Proposition \ref{prop:case-by-case} (iii), it suffices to prove the following proposition.

\begin{prop}\label{prop:unitary-case}
 We have $\mathcal{F}^{\mathit{wa}}=\bigcup_{g\in\U(V)}gU_3$.
\end{prop}

We need some preparation for the proof of this proposition.

\begin{lem}\label{lem:max-isotropic-transitive-unitary}
 For every $W_1,W_2\in \UGr(V)$, there exists $g\in\U(\Lambda)$ such that $g\Lambda_{W_1}=\Lambda_{W_2}$.
\end{lem}

\begin{prf}
 We allow the cases $n=0,1,2$. We put $r=\lfloor\frac{n}{2}\rfloor$ and prove the claim by induction on $r$.
 If $r=0$, the claim is clear since $W_1=W_2=0$. Assume $r\ge 1$ and that the claim holds for $r-1$.
 We may assume that $W_2=\langle e_1,\ldots,e_r\rangle_E$. We take an $\mathcal{O}_E$-basis $f_1,\ldots,f_r$ of $\Lambda_{W_1}$.
 We will prove that there exists $g\in \U(\Lambda)$ such that $g(f_1)=e_1$.
 If $p\neq 2$ or $n$ is even or $E/F$ is unramified, it follows immediately from \cite[\S 9.1, Remark 2]{MR0195863}.
 However, if $p=2$, $n$ is odd and $E/F$ is ramified, then $\Lambda$ is not maximal in the sense of \cite[\S 9.1]{MR0195863},
 so we cannot apply \cite[\S 9.1, Remark 2]{MR0195863} directly. Here we give a direct proof. We write $f_1=(x_1,\ldots,x_n)$.
 We have $x_{i_0}\in\mathcal{O}_E^\times$ for some $1\le i_0\le n$. If $n=2r+1$ is odd, we may take $i_0$ so that $i_0\neq r+1$.
 Indeed, if $x_i\in \mathfrak{p}_E$ for every $1\le i\le n$ with $i\neq r+1$, 
 then $x_{r+1}c(x_{r+1})=-\sum_{1\le i\le n,i\neq r+1}x_ic(x_{n+1-i})+h_n(f_1,f_1)=-\sum_{1\le i\le n,i\neq r+1}x_ic(x_{n+1-i})\in\mathfrak{p}_E^2$,
 which is a contradiction.
 Therefore, by replacing $W_1$ with $\sigma W_1$ for some permutation matrix $\sigma$ in $\U(\Lambda)$,
 we may assume that $i_0=1$. Moreover, by replacing $f_1$ with $x_1^{-1}f_1$, we may assume that $x_1=1$.
 Let $g'\colon \Lambda\to \Lambda$ be the $\mathcal{O}_E$-homomorphism defined by $e_1\mapsto f_1$, $e_i\mapsto e_i-h_n(e_i,f_1)e_n$ $(2\le i\le n-1)$
 and $e_n\mapsto e_n$. It is easy to see that $g'\in \U(\Lambda)$. Therefore, $g=g'^{-1}\in \U(\Lambda)$ satisfies $g(f_1)=e_1$.
 The remaining part is exactly the same as in the proof of Lemma \ref{lem:max-isotropic-transitive}.
\end{prf}

\begin{lem}\label{lem:matrix-unitary}
 We put $r=\lfloor\frac{n}{2}\rfloor$. For $W\in \UGr(V)$, we take an $\mathcal{O}_E$-basis $y_1,\ldots,y_r$ of $\Lambda_W$ and write 
 $y_j=(y_{1j},\ldots,y_{nj})$ for each $1\le j\le r$. Let $s$ be an integer satisfying $1\le s\le r$, and assume that 
 $y_{ij}\in\mathfrak{p}_E$ for every $n+1-s\le i\le n$ and $1\le j\le r$.
 Then, the rank of the matrix $(y_{ij}\bmod\mathfrak{p}_E)_{1\le i\le s,1\le j\le r}$ is equal to $s$.
 In particular, we can replace $y_1,\ldots,y_r$ so that $(y_{ij})_{1\le i,j\le s}$ is the identity matrix and
 $y_{ij}=0$ for every $1\le i\le s$ and $s+1\le j\le r$.
\end{lem}

\begin{prf}
 We allow the cases $n=0,1,2$. We prove the claim by induction on $s\ge 1$. Assume that the claim holds for $s-1$.

 First assume that $y_{11},\ldots,y_{1r}\in\mathfrak{p}_E$. We write $\overline{\Lambda}_W$ for the image of $\Lambda_W\subset \Lambda=\mathcal{O}_E^n$
 in $k_E^n$. By the assumption, $\overline{\Lambda}_W$ is contained in $0\oplus k_E^{n-2}\oplus 0$.
 If $E/F$ is unramified, then $h_n$ induces a non-degenerate Hermitian form $\overline{h}_n$ on $k_E^{n-2}$, and $\overline{\Lambda}_W$ is
 totally isotropic with respect to $\overline{h}_n$. Therefore we have $r=\dim \overline{\Lambda}_W\le r-1$, which is a contradiction.
 If $E/F$ is ramified, then $h_n$ induces a non-degenerate symmetric form $\overline{h}_n$ on $k_E^{n-2}$,
 and $\overline{\Lambda}_W$ is contained in $(\overline{\Lambda}_W)^{\perp}$.
 Therefore we have $r=\dim \overline{\Lambda}_W\le \dim (\overline{\Lambda}_W)^{\perp}=n-2-r\le (2r+1)-2-r=r-1$, which is a contradiction again.

 Hence we conclude that one of $y_{11},\ldots,y_{1r}$ belongs to $\mathcal{O}_E^\times$. By permuting $y_1,\ldots,y_r$, we may assume that
 $y_{11}\in\mathcal{O}_E^\times$. By replacing $y_1,\ldots,y_r$ further, we may assume that $y_{11}=1$ and $y_{12}=\cdots=y_{1r}=0$.
 Now we can apply the induction hypothesis to $\langle (y_{2j},\ldots,y_{n-1,j})\mid 2\le j\le r\rangle_E\in\UGr(E^{n-2})$ to conclude the claim.
\end{prf}

\begin{prf}[of Proposition \ref{prop:unitary-case}]
 Let $(C,C^+)$ be a complete affinoid field over $(\C_p,\mathcal{O}_{\C_p})$ and $L\in \mathcal{F}^{\mathit{wa}}(C,C^+)$.
 We show that there exists $g\in \U(V)$ such that $gL\in U_3(C,C^+)$. By Proposition \ref{prop:unitary-tubular-nbhd} (iii),
 there exists an integer $m\ge 1$ such that $L\in U_m(C,C^+)$. 
 We have only to find $g\in \U(V)$ satisfying $gL\in U_{m-1}(C,C^+)$ under the assumption $m\ge 4$.
 We set $I_L=\{W\in \UGr(V)\mid L\notin U_{m-1,W}(C,C^+)\}$. If $I_L=\varnothing$, then $L\in U_{m-1}(C,C^+)$, so we may assume $I_L\neq \varnothing$.
 For $W\in I_L$, we can find $g\in \U(\Lambda)$ such that $gW=\langle e_1,\ldots,e_r\rangle_E$ by Lemma \ref{lem:max-isotropic-transitive-unitary},
 where we put $r=\lfloor\frac{n}{2}\rfloor$ as in Lemma \ref{lem:matrix-unitary}.
 By Proposition \ref{prop:unitary-tubular-nbhd} (iv), we have $gL\notin gU_{m-1,W}(C,C^+)=U_{m-1,gW}(C,C^+)$.
 Replacing $L$ with $gL$, we may assume that $W=\langle e_1,\ldots,e_r\rangle_E$.
 Then we have $\Lambda_{W^\perp}=\langle e_1,\ldots,e_{n-r}\rangle_{\mathcal{O}_E}$.

 Let $x=(x_1,\ldots,x_n)$ be a $C^+$-basis of $\Lambda_L$.
 By $L\notin U_{m-1,\langle e_1,\ldots,e_r\rangle_E}(C,C^+)$ we have $\lvert x_i\rvert<\lvert\varpi_E^{m-1}\rvert\le \lvert\varpi_E^3\rvert$
 for every $r+1\le i\le n$. Therefore, one of $x_1,\ldots,x_r$ belongs to $(C^+)^\times$.
 Let $J_r$ be the $r\times r$ matrix $(\delta_{i,r+1-j})_{i,j}$, where $\delta_{i,r+1-j}$ denotes the Kronecker delta.
 By applying Lemma \ref{lem:O_C-basis} to $(x_1,\ldots,x_r)$, after replacing $L$ with its translate by 
 $\diag(h,\underbrace{1,\ldots,1}_{n-2r},J_rc({}^t\!h^{-1})J_r)\in \U(\Lambda)$ for a suitable $h\in\GL_r(\mathcal{O}_E)$,
 we may assume that $(x_1,\ldots,x_r)$ satisfies the following: for some integers $s$, $s'$ with $1\le s\le s'\le r$, 
 \begin{itemize}
  \item $x_{s+1},\ldots,x_{s'}$ lie in $\varpi_EC^+$, and the images of $x_1,\ldots,x_s,\varpi_E^{-1}x_{s+1},\ldots,\varpi_E^{-1}x_{s'}$ in
	$C^+/\varpi_E C^+$ are linearly independent over $k_E$, and
  \item $x_{s'+1},\ldots,x_r$ lie in $\varpi_E^2 C^+$.
 \end{itemize}
 Note that the new $L$ still belongs to $U_m(C,C^+)\smallsetminus U_{m-1,\langle e_1,\ldots,e_r\rangle_E}(C,C^+)$
 by Proposition \ref{prop:unitary-tubular-nbhd} (iv).

 We put $t=\diag(\underbrace{\varpi_E^{-1},\ldots,\varpi_E^{-1}}_{s},\underbrace{1,\ldots,1}_{n-2s},\underbrace{c(\varpi_E),\ldots,c(\varpi_E)}_s)\in \U(V)$. As $\varpi_E^{-1}t^{-1}x$ equals
 \[
  (x_1,\ldots,x_s,\varpi_E^{-1}x_{s+1},\ldots,\varpi_E^{-1}x_{n-s},N_{E/F}(\varpi_E)^{-1}x_{n+1-s},\ldots,N_{E/F}(\varpi_E)^{-1}x_n),
 \]
 we have $\Lambda_{t^{-1}L}=\langle \varpi_E^{-1}t^{-1}x\rangle_{C^+}$ (here we use $\lvert x_i\rvert\le \lvert\varpi_E^2\rvert$ for $r+1\le i\le n$).
 We prove that for every $W\in\UGr(V)$ we have $t^{-1}L\in U_{m-1,W}(C,C^+)$, that is, 
 $\llparenthesis t^{-1}L,W^\perp\rrparenthesis\ge \lvert\varpi_E^{m-1}\rvert$.
 Let $y_1,\ldots,y_r$ be an $\mathcal{O}_E$-basis of $\Lambda_W$. Further, we take $y_{r+1},\ldots,y_{n-r}\in \Lambda_{W^\perp}$ so that
 $y_1,\ldots,y_{n-r}$ form a basis of $\Lambda_{W^\perp}$. We write $y_j=(y_{1j},\ldots,y_{nj})$ for each $1\le j\le n-r$. 

 First suppose that $y_{i_0j_0}\in \mathcal{O}_E^\times$ for some $n+1-s\le i_0\le n$ and $1\le j_0\le n-r$.
 Then we have $h_n(\varpi_E^{-1}t^{-1}x,y_{j_0})\in \sum_{i=1}^s c(y_{n+1-i,j_0})x_i+\sum_{i=s+1}^{s'}c(y_{n+1-i,j_0})\varpi_E^{-1}x_i+\varpi_EC^+$;
 here we use $\lvert\varpi_E^{-1}x_i\rvert\le \lvert\varpi_E\rvert$ for $s'+1\le i\le n-s$ and $\lvert N_{E/F}(\varpi_E)^{-1}x_i\rvert\le \lvert\varpi_E\rvert$ for $n+1-s\le i\le n$.
 Since the images of $x_1,\ldots,x_s,\varpi_E^{-1}x_{s+1},\ldots,\varpi_E^{-1}x_{s'}$ in $C^+/\varpi_EC^+$ are linearly independent over $k_E$,
 we have $h_n(\varpi_E^{-1}t^{-1}x,y_{j_0})\in C^+\smallsetminus \varpi_EC^+$, in other words,
 $\lvert h_n(\varpi_E^{-1}t^{-1}x,y_{j_0})\rvert>\lvert\varpi_E\rvert$.
 Therefore we have $\llparenthesis t^{-1}L,W^\perp\rrparenthesis>\lvert\varpi_E\rvert>\lvert\varpi_E^{m-1}\rvert$, as desired.

 Next we consider the case where $y_{ij}\in\mathfrak{p}_E$ for every $n+1-s\le i\le n$ and $1\le j\le n-r$.
 By Lemma \ref{lem:matrix-unitary}, we can replace $y_1,\ldots,y_r$ so that $(y_{ij})_{1\le i,j\le s}$ is the identity matrix and
 $y_{ij}=0$ for every $1\le i\le s$ and $s+1\le j\le r$. Further, we can also replace $y_{r+1},\ldots,y_{n-r}$ so that $y_{ij}=0$ for every
 $1\le i\le s$ and $r+1\le j\le n-r$. We prove $\Lambda_{tW^\perp}=\langle \varpi_Ety_1,\ldots,\varpi_Ety_s,ty_{s+1},\ldots,ty_{n-r}\rangle_{\mathcal{O}_E}$.
 Clearly we have $\Lambda_{tW^\perp}\supset\langle \varpi_Ety_1,\ldots,\varpi_Ety_s,ty_{s+1},\ldots,ty_{n-r}\rangle_{\mathcal{O}_E}$.
 Let us prove the reverse inclusion.
 For an arbitrary element $z\in \Lambda_{tW^\perp}$, there exist $a_1,\ldots,a_{n-r}\in E$ such that 
 $z=a_1(\varpi_E ty_1)+\cdots+a_s(\varpi_E ty_s)+a_{s+1}(ty_{s+1})+\cdots+a_{n-r}(ty_{n-r})$. It suffices to show $a_1,\ldots,a_{n-r}\in \mathcal{O}_E$. 
 As $a_1(\varpi_E ty_1)+\cdots+a_s(\varpi_E ty_s)+a_{s+1}(ty_{s+1})+\cdots+a_{n-r}(ty_{n-r})\in \Lambda=\mathcal{O}_E^n$,
 by looking at the first $s$ entries, we obtain $a_1,\ldots,a_s\in\mathcal{O}_E$. 
 Therefore we have $a_{s+1}(ty_{s+1})+\cdots+a_{n-r}(ty_{n-r})\in \mathcal{O}_E^n$. Let $d\ge 0$ be the minimal non-negative integer
 such that $\varpi_E^da_{s+1},\ldots,\varpi_E^da_{n-r}\in \mathcal{O}_E$. If $d\ge 1$, then we have
 $\varpi_E^da_{s+1}(ty_{s+1})+\cdots+\varpi_E^da_{n-r}(ty_{n-r})\in \mathfrak{p}_E^n$.
 Therefore we obtain the equality $\overline{\varpi_E^da_{s+1}}\cdot \overline{y_{s+1}}+\cdots+\overline{\varpi_E^da_{n-r}}\cdot \overline{y_{n-r}}=0$
 in $k_E^n$, where $\overline{a}$ (resp.\ $\overline{y}$) denotes the image of $a\in \mathcal{O}_E$ (resp.\ $y\in \mathcal{O}_E^n$)
 in $k_E$ (resp.\ $k_E^n$). Note that we have $\overline{ty_j}=\overline{y_j}$ for every $s+1\le j\le n-r$,
 as $y_{1j}=\cdots=y_{sj}=0$ and $y_{n+1-s,j},\ldots,y_{nj}\in \mathfrak{p}_E$. Since $\overline{y_{s+1}},\ldots,\overline{y_{n-r}}\in k_E^n$ 
 are linearly independent over $k_E$, we have $\overline{\varpi_E^da_{s+1}}=\cdots=\overline{\varpi_E^da_{n-r}}=0$, in other words,
 $\varpi_E^da_{s+1},\ldots,\varpi_E^da_{n-r}\in \mathfrak{p}_E$. Therefore $\varpi_E^{d-1}a_{s+1},\ldots,\varpi_E^{d-1}a_{n-r}$ lie in $\mathcal{O}_E$,
 which contradicts the minimality of $d$. Hence we have $d=0$ and $a_{s+1},\ldots,a_{n-r}\in\mathcal{O}_E$.
 This completes the proof of $\Lambda_{tW^\perp}=\langle \varpi_E ty_1,\ldots,\varpi_E ty_s,ty_{s+1},\ldots,ty_{n-r}\rangle_{\mathcal{O}_E}$.
 Since $\varpi_E^{-1}c(\varpi_E)\in \mathcal{O}_E^\times$, we also have 
 $\Lambda_{tW^\perp}=\langle c(\varpi_E)ty_1,\ldots,c(\varpi_E)ty_s,ty_{s+1},\ldots,ty_{n-r}\rangle_{\mathcal{O}_E}$.
 Therefore we have
 \begin{align*}
  \llparenthesis t^{-1}L,W^\perp\rrparenthesis&=\max\{\lvert h_n(\varpi_E^{-1}t^{-1}x,y_1)\rvert,\ldots,\lvert h_n(\varpi_E^{-1}t^{-1}x,y_{n-r})\rvert\}\\
  &=\max\{\lvert\varpi_E^{-2}\rvert\cdot \lvert h_n(x,c(\varpi_E)ty_1)\rvert,\ldots,\lvert \varpi_E^{-2}\rvert\cdot \lvert h_n(x,c(\varpi_E)ty_s)\rvert,\\
  &\qquad\qquad\qquad\qquad\qquad \lvert \varpi_E^{-1}\rvert\cdot \lvert h_n(x,ty_{s+1})\rvert,\ldots,\lvert \varpi_E^{-1}\rvert\cdot \lvert h_n(x,ty_{n-r})\rvert\}\\
  &\ge \lvert\varpi_E^{-1}\rvert\cdot\max\{\lvert h_n(x,c(\varpi_E)ty_1)\rvert,\ldots,\lvert h_n(x,c(\varpi_E)ty_s)\rvert,\\
  &\qquad\qquad\qquad\qquad\qquad\qquad\qquad\lvert h_n(x,ty_{s+1})\rvert,\ldots,\lvert h_n(x,ty_{n-r})\rvert\}\\
  &=\lvert\varpi_E^{-1}\rvert\cdot\llparenthesis L,tW^\perp\rrparenthesis=\lvert\varpi_E^{-1}\rvert\cdot\llparenthesis L,(tW)^\perp\rrparenthesis\ge \lvert\varpi_E^{m-1}\rvert,
 \end{align*}
 as desired. In the last inequality we use $L\in U_{m,tW}(C,C^+)$.
 Now we conclude that $t^{-1}L\in U_{m-1}(C,C^+)$.
\end{prf}

\begin{rem}\label{rem:Shen-gap}
 Shen \cite{XuShen-RZ} considered the existence of Siegel domains in the unitary case, under the conditions that $F=\Q_p$, $E/F$ is unramified and $p>2$.
 However, the proof in \cite{XuShen-RZ} contains several gaps. Here we freely use the notation in \cite{XuShen-RZ}, and assume that
 $n$ is odd for simplicity.
 First, in the proof of \cite[Theorem 10.8]{XuShen-RZ}, Shen used the following two inputs
 to conclude that $\{T.g\mathcal{D}'\}_{T,g}$ covers $\mathcal{M}$:
 \begin{itemize}
  \item $\mathcal{M}^{\mathrm{rig}}=\bigcup_{T\in G(\Z_p)\backslash G(\Q_p)/G(\Z_p),g\in J^{\mathrm{der}}_b(\Q_p)/\Stab(\Lambda)}T.g\mathcal{D}'^{\mathrm{rig}}$.
  \item the family $\{T.g\mathcal{D}'\}_{T,g}$ is locally finite.
 \end{itemize}
 However, the notion of local finiteness considered in \cite[p.~873]{XuShen-RZ} is weaker than the usual one, and one cannot conclude the desired equality
 $\mathcal{M}=\bigcup_{T,g}T.g\mathcal{D}'$ from it.
 For example, we take sequences of rational numbers $\{a_n\}_{n\ge 1}$, $\{b_n\}_{n\ge 1}$ satisfying $0=a_1<a_2<\cdots<\sqrt{2}<\cdots<b_2<b_1=2$
 and $\lim_{n\to\infty}a_n=\lim_{n\to\infty}b_n=\sqrt{2}$, and consider quasi-compact open subsets of 
 $\mathbb{D}^1=\Spa(\C_p\langle T\rangle,\mathcal{O}_{\C_p}\langle T\rangle)$ defined as follows:
 $U_n=\{z\in\mathbb{D}^1\mid \lvert p^{a_{n+1}}\rvert\le \lvert z\rvert\le \lvert p^{a_n}\rvert\}$ and 
 $V_n=\{z\in\mathbb{D}^1\mid \lvert p^{b_n}\rvert\le \lvert z\rvert\le \lvert p^{b_{n+1}}\rvert\}$ for $n\ge 1$,
 and $V_0=\{z\in \mathbb{D}^1\mid \lvert z\rvert\le \lvert p^2\rvert\}$.
 Clearly all classical points of $\mathbb{D}^1$ are contained in $\bigcup_{n\ge 1}U_n\cup\bigcup_{n\ge 0}V_n$,
 and the family $\{U_n\}_{n\ge 1}\sqcup \{V_n\}_{n\ge 0}$ is locally finite in the sense of \cite[p.~873]{XuShen-RZ},
 but the point $x$ in $\mathbb{D}^1$ corresponding to ``the disk of radius $\lvert p^{\sqrt{2}}\rvert$'' is not contained in 
 $\bigcup_{n\ge 1}U_n\cup\bigcup_{n\ge 0}V_n$. Note that every open neighborhood of this $x$ intersects infinitely many of
 the open sets $U_n$ and $V_n$, so $\{U_n\}_{n\ge 1}\sqcup \{V_n\}_{n\ge 0}$ is not locally finite in the usual sense.

 Moreover, GPT-6.2 astra pointed out that Proposition 10.1 and Corollary 10.6 in \cite{XuShen-RZ} are false. 
 To explain it, we assume that $n=3$ and take a $\Q_{p^2}$-basis $e_1,e_2,e_3$ of $\mathbf{N}_0$
 such that the Gram matrix $(\{e_i,e_j\})_{i,j}$ is equal to $J_3=\Bigl(\begin{smallmatrix}0&0&1\\0&1&0\\1&0&0\end{smallmatrix}\Bigr)$. 
 Then, the group $J_b(\Q_p)$ is identified with the unitary similitude group with respect to the Hermitian form corresponding to $J_3$.
 Further, the weakly admissible locus $\mathcal{F}(G,b,\{\mu\})^{\mathit{wa}}$ attached to the triple $(G,b,\{\mu\})$ in \cite[p.~835]{XuShen-RZ}
 is identified with $(\P^2)^{\ad}\smallsetminus \bigcup_{W\in \UGr(\Q_{p^2}^3)}\P(W)^{\ad}_{\C_p}$ (see \cite[p.~341]{MR1810589}). 
 Now, consider a point $y=[0:1:0]\in (\P^2)^{\ad}\smallsetminus \bigcup_{W\in \UGr(\Q_{p^2}^3)}\P(W)^{\ad}_{\C_p}=\mathcal{F}(G,b,\{\mu\})^{\mathit{wa}}$.
 Since it is a classical point, it lies in $\mathcal{F}(G,b,\{\mu\})^a$. Therefore, by \cite[p.~843]{XuShen-RZ}
 we can take an element $x\in\mathcal{M}^0$ such that $\pi_{\mathrm{dR}}(x)=y$.
 Since $y$ is fixed by $h=\diag(p,1,p^{-1})\in J_b(\Q_p)$, we have $\pi_{\mathrm{dR}}(hx)=\pi_{\mathrm{dR}}(x)$. 
 By \cite[Proposition 2.3]{XuShen-RZ}, there exists a Hecke operator $T$ such that $hx\in Tx$.
 As $x$ and $hx$ lie in $\mathcal{M}^0$, there exists an integer $m$ such that
 $T=G(\Z_p)\diag(p^m,1,p^{-m})G(\Z_p)$. We take $\Lambda_x\in \mathcal{L}_0$ such that $\mathrm{sp}(x)\in \mathcal{M}_{\Lambda_x}^0$.
 If \cite[Proposition 10.1]{XuShen-RZ} holds, 
 then we have $hx\in Tx\subset T.\mathrm{sp}^{-1}(\mathcal{M}_{\Lambda_x}^0)\subset \mathrm{sp}^{-1}(\mathcal{M}_{\Lambda_x}^0)$.
 On the other hand, we have $hx\in \mathrm{sp}^{-1}(\mathcal{M}_{h\Lambda_x}^0)$. 
 Therefore the intersection $\mathrm{sp}^{-1}(\mathcal{M}_{\Lambda_x}^0)\cap \mathrm{sp}^{-1}(\mathcal{M}_{h\Lambda_x}^0)$ is non-empty, hence
 $h\Lambda_x=\Lambda_x$. In particular, $h$ belongs to a compact-mod-center subgroup $\Stab(\Lambda_x)$ of $J_b(\Q_p)$, which contradicts
 the form of $h=\diag(p,1,p^{-1})$.

 Next we consider \cite[Corollary 10.6]{XuShen-RZ}. 
 By \cite[Proposition 6.7]{XuShen-RZ} and $\pi_{\mathrm{dR}}(\mathcal{C})=\pi_{\mathrm{dR}}(\mathcal{C}^0)$ (see \cite[p.~866]{XuShen-RZ}),
 there exists $x\in \mathcal{C}^0$ such that $\pi_{\mathrm{dR}}(x)=y$.
 Further, by $\mathcal{C}^0=\bigcup_{g\in J_b^{\mathrm{der}}(\Q_p)/\Stab(\Lambda)}g\mathcal{D}$ (see \cite[p.~873]{XuShen-RZ}), 
 we can find $g\in J_b^{\mathrm{der}}(\Q_p)$ such that $x\in g\mathcal{D}$.
 By the same argument as above, for each integer $r\ge 1$, we can find an integer $m_r$ such that the Hecke operator
 $T_r=G(\Z_p)\diag(p^{m_r},1,p^{-m_r})G(\Z_p)$ satisfies $h^rx\in T_rx$. 
 Therefore we have $x\in T_r.h^{-r}g\mathcal{D}$. Since the lattices $h^{-r}g\Lambda$ for $r\ge 1$ are pairwise distinct,
 the point $x$ belongs to infinitely many members of
 $\{T.g\mathcal{D}\}_{T\in G(\Z_p)\backslash G(\Q_p)/G(\Z_p),g\in J_b^{\mathrm{der}}(\Q_p)/\Stab(\Lambda)}$.
 This contradicts the claim of \cite[Corollary 10.6]{XuShen-RZ}. In the same way we can also see that $\{T.g\mathcal{D}'\}_{T,g}$ in the proof of
 \cite[Theorem 10.8]{XuShen-RZ} is not locally finite in the sense of \cite[p.~873]{XuShen-RZ}.
\end{rem}

\subsection{The Hilbert--Blumenthal case}\label{subsec:HB}
Finally we consider the Hilbert--Blumenthal case. Let $n\ge 2$ be an integer, $E$ a finite extension of $F$, and $G=\Res_{E/F}\PGL_n$.
Let $\tau_1$ and $\tau_2$ be two distinct elements of $\Sigma=\Hom_{\text{$F$-alg}}(E,\overline{F})$, and 
$\mu=(\mu_\tau)_{\tau\in\Sigma}$ the cocharacter determined by $\mu_{\tau_1}=\omega_1^\vee$, $\mu_{\tau_2}=\omega_{n-1}^\vee$ and 
$\mu_{\tau'}=0$ for $\tau'\in \Sigma\smallsetminus \{\tau_1,\tau_2\}$. By exchanging $\tau_1$ and $\tau_2$, 
it suffices to consider $\mathcal{F}(G,1,\{\mu\})^{\mathit{wa}}$ in place of $\mathcal{F}(G,1,\{\mu^{-1}\})^{\mathit{wa}}$.
We regard $E$ as a subfield of $\C_p$ via $\tau_1$, and omit $\tau_1$ from the notation. We write $\tau$ for $\tau_2$.
For simplicity we write $\mathcal{F}$ and $\mathcal{F}^{\mathit{wa}}$ for $\mathcal{F}(G,\{\mu\})$ and
$\mathcal{F}(G,1,\{\mu\})^{\mathit{wa}}$, respectively. We have $\mathcal{F}=\P^{n-1}\times (\P^{n-1})^\vee$,
where $(\P^{n-1})^\vee$ denotes the variety classifying $(n-1)$-dimensional subspaces of $\C_p^n$.
The group $G(F)=\PGL_n(E)$ (or rather $\GL_n(E)$) acts on $\mathcal{F}$ by $g(L,H)=(gL,\tau(g)H)$.
For an $E$-subspace $0\neq W\subsetneq E^n$, we put $\mathcal{F}_W=\{(L,H)\in \mathcal{F}\mid L\subset W\otimes_E\C_p, \tau(W)\otimes_{\tau(E)}\C_p\subset H\}$. By using Proposition \ref{prop:ss=wa} (iii), we can observe that
$\mathcal{F}^{\mathit{wa}}=\mathcal{F}^{\ad}\smallsetminus\bigcup_{0\neq W\subsetneq E^n}\mathcal{F}_W^{\ad}$.

Let $(-,-)$ denote the symmetric pairing $E^n\times E^n\to E$; $((a_i),(b_i))\mapsto \sum_{i=1}^n a_ib_i$.
We use the same symbol to denote the base change of the pairing to an extension of $E$.
It gives an identification $(\mathbb{P}^{n-1})^\vee\xrightarrow{\cong}\mathbb{P}^{n-1}$; $H\mapsto H^\perp$, where $H^\perp$ denotes the orthogonal
complement of $H$ with respect to $(-,-)$. Under this identification, we have $\mathcal{F}=\P^{n-1}\times\P^{n-1}$,
$\mathcal{F}_W=\{(L_1,L_2)\in\mathcal{F}\mid L_1\subset W\otimes_E\C_p,L_2\subset \tau(W^\perp)\otimes_{\tau(E)}\C_p\}$ and
$g(L_1,L_2)=(gL_1,\tau({}^t\!g^{-1})L_2)$.
Let $(C,C^+)$ be a complete affinoid field over $(\C_p,\mathcal{O}_{\C_p})$, and let $\lvert-\rvert$ denote the continuous valuation on $C$
corresponding to the closed point of $\Spa(C,C^+)$.
For $x\in (C^+)^n$ and $a\in \mathcal{O}_E^n\cup \tau(\mathcal{O}_E)^n$, we put $\llparenthesis x,a\rrparenthesis=\lvert (x,a)\rvert$.
For $L\in \P^{n-1}(C)$ (resp.\ an $E$-subspace $W$ of $E^n$), we put $\Lambda_L=(C^+)^n\cap L$ (resp.\ $\Lambda_W=\mathcal{O}_E^n\cap W$),
which is a free $C^+$-module (resp.\ $\mathcal{O}_E$-module) of rank $1$ (resp.\ $\dim_EW$).
For $(L_1,L_2)\in \mathcal{F}(C)$ and an $E$-subspace $0\neq W\subsetneq E^n$, we define
$\llparenthesis L_1,W^\perp\rrparenthesis=\max\{\llparenthesis x,a\rrparenthesis \mid x\in \Lambda_{L_1},a\in \Lambda_{W^\perp}\}$,
and $\llparenthesis L_2,\tau(W)\rrparenthesis=\max\{\llparenthesis x,\tau(a)\rrparenthesis \mid x\in \Lambda_{L_2},a\in \Lambda_W\}$.
Since $\Lambda_{L_1}$, $\Lambda_{L_2}$, $\Lambda_{W^\perp}$ and $\Lambda_W$ are finitely generated, the maximum always exists.
We fix a uniformizer $\varpi_E$ of $\mathcal{O}_E$.
The following proposition can be proven in the same way as Proposition \ref{prop:Drinfeld-tubular-nbhd}.

\begin{prop}\label{prop:HB-tubular-nbhd}
 \begin{enumerate}
  \item For an $E$-subspace $0\neq W\subsetneq E^n$ and an integer $m\ge 1$, there exists a unique quasi-compact open subset $U_{m,W}$
	of $\mathcal{F}^{\ad}$ satisfying the following: for every complete affinoid field $(C,C^+)$ over $(\C_p,\mathcal{O}_{\C_p})$, we have 
	\[
	 U_{m,W}(C,C^+)=\{(L_1,L_2)\in \mathcal{F}^{\ad}(C,C^+)\mid \max\{\llparenthesis L_1,W^\perp\rrparenthesis,\llparenthesis L_2,\tau(W)\rrparenthesis\}\ge \lvert\varpi_E^m\rvert\}.
	\]
  \item The subset $U_m=\bigcap_{0\neq W\subsetneq E^n}U_{m,W}$ of $\mathcal{F}^{\ad}$ is a quasi-compact open subset of $\mathcal{F}^{\mathit{wa}}$.
  \item We have $\mathcal{F}^{\mathit{wa}}=\bigcup_{m\ge 1}U_m$.
  \item For $g\in\GL_n(\mathcal{O}_E)$, we have $gU_{m,W}=U_{m,gW}$ and $gU_m=U_m$.
 \end{enumerate}
\end{prop}

\begin{prop}\label{prop:HB-stable}
 We put $\Omega_E^{n-1}=(\P^{n-1})^{\ad}\smallsetminus \bigcup_{\textnormal{$H\subset \P^{n-1}$: $E$-rational hyperplane}}H^{\ad}_{\C_p}$.
 Let $(C,C^+)$ be a complete affinoid field over $(\C_p,\mathcal{O}_{\C_p})$ and $(L_1,L_2)\in \mathcal{F}^{\mathit{wa}}(C,C^+)$. 
 If $L_1\in \Omega_E^{n-1}(C,C^+)$, we have $(L_1,L_2)\in \bigcup_{g\in \GL_n(E)}gU_2(C,C^+)$.
\end{prop}

\begin{prf}
 We write $U_2^{\Omega}$ for the quasi-compact open subset $U_2\subset \Omega_E^{n-1}$ defined in Proposition \ref{prop:Drinfeld-tubular-nbhd} (ii).
 By Proposition \ref{prop:Drinfeld-case}, there exists $g\in \GL_n(E)$ such that $gL_1\in U_2^{\Omega}(C,C^+)$. For this $g$ we prove
 $g(L_1,L_2)\in U_2(C,C^+)$.
 We take an arbitrary $E$-subspace $0\neq W\subsetneq E^n$. Further, we take $w\in \Lambda_W\smallsetminus \{0\}$ such that $\varpi_E^{-1}w\notin\Lambda_W$.
 As $gL_1\in U_2^{\Omega}$, there exists $x\in \Lambda_{gL_1}$ such that $\llparenthesis x,w\rrparenthesis\ge \lvert\varpi_E^2\rvert$.
 In particular we have $\llparenthesis gL_1,W\rrparenthesis\ge \lvert\varpi_E^2\rvert$, hence $g(L_1,L_2)\in U_{2,W^\perp}(C,C^+)$. Varying $W$,
 we conclude that $g(L_1,L_2)\in U_2(C,C^+)$, as desired.
\end{prf}

Next we consider the case which is not covered by Proposition \ref{prop:HB-stable}. To do it, it is convenient not to use the identification
$(\P^{n-1})^\vee\xrightarrow{\cong}\P^{n-1}$.

\begin{defn}
 Let $Y$ be the locally closed subvariety of $\mathcal{F}$ defined by $Y=\{(L,H)\in\mathcal{F}\mid L\subset \langle e_1,\ldots,e_{n-1}\rangle_{\C_p}\nsubset H\}$. We put $Y^{\mathit{wa}}=Y^{\ad}\cap \mathcal{F}^{\mathit{wa}}$. Since $Y^{\mathit{wa}}=\{(L,H)\in\mathcal{F}^{\mathit{wa}}\mid L\subset\langle e_1,\ldots,e_{n-1}\rangle_{\C_p}\}$, it is a closed adic subspace of $\mathcal{F}^{\mathit{wa}}$. We write $\mathcal{F}_{n-1}$
 (resp.\ $\mathcal{F}_{n-1}^{\mathit{wa}}$) for the algebraic variety (resp.\ adic space) obtained by replacing $n$ with $n-1$
 in the definition of $\mathcal{F}$ (resp.\ $\mathcal{F}^{\mathit{wa}}$).

 We define morphisms of algebraic varieties $r\colon Y\to \mathcal{F}_{n-1}$ and $s\colon\mathcal{F}_{n-1}\to Y$ by
 $(L,H)\mapsto (L,H\cap \langle e_1,\ldots,e_{n-1}\rangle_{\C_p})$ and $(L,H')\mapsto (L,H'+\langle e_n\rangle_{\C_p})$, respectively.
 Since $r\circ s=\id_{\mathcal{F}_{n-1}}$, $r$ is surjective.
 Let $P$ be the parabolic subgroup of $\GL_{n,E}$ obtained as the stabilizer of $\langle e_1,\ldots,e_{n-1}\rangle_E$.
 Then, for $x=(L,H)\in Y(\C_p)$ and $g=\begin{pmatrix}h&v\\0&a\end{pmatrix}\in P(E)$, we have $r(gx)=hr(x)$.
 Indeed, we have $\tau(h)(H\cap \langle e_1,\ldots,e_{n-1}\rangle_{\C_p})=\tau(g)(H\cap \langle e_1,\ldots,e_{n-1}\rangle_{\C_p})=\tau(g)H\cap \tau(g)\langle e_1,\ldots,e_{n-1}\rangle_{\C_p}=\tau(g)H\cap \langle e_1,\ldots,e_{n-1}\rangle_{\C_p}$.
\end{defn}

\begin{lem}\label{lem:HB-reduction}
 We have $(r^{\ad})^{-1}(\mathcal{F}^{\mathit{wa}}_{n-1})=Y^{\mathit{wa}}$.
 In particular we have $s^{\ad}(\mathcal{F}^{\mathit{wa}}_{n-1})\subset Y^{\mathit{wa}}$.
\end{lem}

\begin{prf}
 Let $(C,C^+)$ be a complete affinoid field over $(\C_p,\mathcal{O}_{\C_p})$ and $x=(L,H)\in Y^{\ad}(C,C^+)=Y(C)$.
 We put $H'=H\cap \langle e_1,\ldots,e_{n-1}\rangle_C$, and prove that $x=(L,H)$ belongs to $Y^{\mathit{wa}}(C,C^+)$ if and only if
 $r^{\ad}(x)=(L,H')$ belongs to $\mathcal{F}^{\mathit{wa}}_{n-1}(C,C^+)$.
 First we assume $(L,H)\in Y^{\mathit{wa}}(C,C^+)$. If an $E$-subspace $0\neq W\subsetneq E^{n-1}$ satisfies 
 $L\subset W\otimes_EC$ and $\tau(W)\otimes_{\tau(E)}C\subset H'$, then $0\neq W\subsetneq E^n$ satisfies 
 $L\subset W\otimes_EC$ and $\tau(W)\otimes_{\tau(E)}C\subset H$, which contradicts $(L,H)\in Y^{\mathit{wa}}(C,C^+)$.
 Therefore $(L,H')$ belongs to $\mathcal{F}^{\mathit{wa}}_{n-1}(C,C^+)$.
 Next we assume $(L,H')\in \mathcal{F}^{\mathit{wa}}_{n-1}(C,C^+)$.
 If an $E$-subspace $0\neq W\subsetneq E^n$ satisfies 
 $L\subset W\otimes_EC$ and $\tau(W)\otimes_{\tau(E)}C\subset H$, then $W'=W\cap \langle e_1,\ldots,e_{n-1}\rangle_E$ satisfies
 $L\subset W'\otimes_EC$ and $\tau(W')\otimes_{\tau(E)}C\subset H\cap \langle e_1,\ldots,e_{n-1}\rangle_C=H'$. As $L\neq 0$ and $H'\subsetneq C^{n-1}$,
 we have $0\neq W'\subsetneq E^{n-1}$. Therefore this contradicts $(L,H')\in \mathcal{F}^{\mathit{wa}}_{n-1}(C,C^+)$.
 Now we conclude that $(r^{\ad})^{-1}(\mathcal{F}^{\mathit{wa}}_{n-1})=Y^{\mathit{wa}}$. Since $r^{\ad}$ is surjective, we have
 $r^{\ad}(Y^{\mathit{wa}})=\mathcal{F}^{\mathit{wa}}_{n-1}$, hence $s^{\ad}(\mathcal{F}^{\mathit{wa}}_{n-1})\subset Y^{\mathit{wa}}$, as desired.
\end{prf}

To show Proposition \ref{prop:case-by-case} (iv), it suffices to prove the following proposition.

\begin{prop}\label{prop:HB-fundamental-domain}
 There exists a quasi-compact open subset $U$ of $\mathcal{F}^{\textit{wa}}$ such that $\mathcal{F}^{\textit{wa}}=\bigcup_{g\in\GL_n(E)}gU$.
\end{prop}

\begin{prf}
 We allow the case $n=1$ and prove the claim by induction on $n\ge 1$. If $n=1$, we have $\mathcal{F}^{\mathit{wa}}=\mathcal{F}^{\ad}$,
 so we can take $\mathcal{F}^{\mathit{wa}}$ as $U$.

 Assume $n\ge 2$, and that there exists a quasi-compact open subset $U_{n-1}$ of 
 $\mathcal{F}^{\mathit{wa}}_{n-1}$ such that $\mathcal{F}^{\mathit{wa}}_{n-1}=\bigcup_{h\in \GL_{n-1}(E)}hU_{n-1}$.
 By Lemma \ref{lem:HB-reduction}, $s^{\ad}(U_{n-1})$ is a quasi-compact subset of $\mathcal{F}^{\mathit{wa}}$. Therefore, we can take
 a quasi-compact open subset $U'$ of $\mathcal{F}^{\mathit{wa}}$ containing $s^{\ad}(U_{n-1})$. We prove that 
 $Y^{\mathit{wa}}\subset \bigcup_{g\in \GL_n(E)}gU'$. Let $(C,C^+)$ be a complete affinoid field over $(\C_p,\mathcal{O}_{\C_p})$
 and $x=(L,H)\in Y^{\mathit{wa}}(C,C^+)$. By Lemma \ref{lem:HB-reduction}, we have $r^{\ad}(x)\in \mathcal{F}^{\mathit{wa}}_{n-1}(C,C^+)$.
 Since $\mathcal{F}^{\mathit{wa}}_{n-1}=\bigcup_{h\in \GL_{n-1}(E)}hU_{n-1}$,
 we can take $h\in \GL_{n-1}(E)$ such that $hr^{\ad}(x)\in U_{n-1}(C,C^+)$. Replacing $x$ with $\diag(h,1)x$,
 we may assume that $r^{\ad}(x)\in U_{n-1}(C,C^+)$.
 For an integer $m$, we put $t_m=\diag(p^m,\ldots,p^m,1)\in \GL_n(E)$. Then we have $\lim_{m\to \infty}t_mx=s^{\ad}(r^{\ad}(x))$.
 Indeed, if we put $H'=H\cap \langle e_1,\ldots,e_{n-1}\rangle_C$, then there exist $a_1,\ldots,a_{n-1}\in C$ such that
 $H=H'+\langle a_1e_1+\cdots+a_{n-1}e_{n-1}+e_n\rangle_C$, and we have $t_mH=H'+\langle p^ma_1e_1+\cdots+p^ma_{n-1}e_{n-1}+e_n\rangle_C$.
 As $s^{\ad}(r^{\ad}(x))\in s^{\ad}(U_{n-1})\subset U'$, we have $t_mx\in U'(C,C^+)$ for sufficiently large $m$.
 Therefore we obtain $x\in \bigcup_{g\in \GL_n(E)}gU'(C,C^+)$, as desired.

 Let $(C,C^+)$ be as above and $(L,H)\in \mathcal{F}^{\mathit{wa}}(C,C^+)$ such that $L\notin \Omega^{n-1}_E(C,C^+)$. We prove that
 $(L,H)\in \bigcup_{g\in \GL_n(E)}gU'(C,C^+)$. As $L\notin \Omega^{n-1}_E(C,C^+)$, there exists an $(n-1)$-dimensional $E$-subspace $W\subset E^n$
 such that $L\subset W\otimes_EC$. We take $g'\in \GL_n(E)$ such that $g'W=\langle e_1,\ldots,e_{n-1}\rangle_E$.
 Then we have $g'(L,H)\in Y^{\mathit{wa}}(C,C^+)\subset \bigcup_{g\in \GL_n(E)}gU'(C,C^+)$ by the argument above.
 Hence we obtain $(L,H)\in \bigcup_{g\in \GL_n(E)}gU'(C,C^+)$, as desired.
 Combining this with Proposition \ref{prop:HB-stable}, we conclude that $U=U_2\cup U'$ satisfies the condition in the proposition.
\end{prf}

\def\cftil#1{\ifmmode\setbox7\hbox{$\accent"5E#1$}\else
  \setbox7\hbox{\accent"5E#1}\penalty 10000\relax\fi\raise 1\ht7
  \hbox{\lower1.15ex\hbox to 1\wd7{\hss\accent"7E\hss}}\penalty 10000
  \hskip-1\wd7\penalty 10000\box7}
  \def\cftil#1{\ifmmode\setbox7\hbox{$\accent"5E#1$}\else
  \setbox7\hbox{\accent"5E#1}\penalty 10000\relax\fi\raise 1\ht7
  \hbox{\lower1.15ex\hbox to 1\wd7{\hss\accent"7E\hss}}\penalty 10000
  \hskip-1\wd7\penalty 10000\box7}
  \def\cftil#1{\ifmmode\setbox7\hbox{$\accent"5E#1$}\else
  \setbox7\hbox{\accent"5E#1}\penalty 10000\relax\fi\raise 1\ht7
  \hbox{\lower1.15ex\hbox to 1\wd7{\hss\accent"7E\hss}}\penalty 10000
  \hskip-1\wd7\penalty 10000\box7}
  \def\cftil#1{\ifmmode\setbox7\hbox{$\accent"5E#1$}\else
  \setbox7\hbox{\accent"5E#1}\penalty 10000\relax\fi\raise 1\ht7
  \hbox{\lower1.15ex\hbox to 1\wd7{\hss\accent"7E\hss}}\penalty 10000
  \hskip-1\wd7\penalty 10000\box7} \def\cprime{$'$} \def\cprime{$'$}
  \newcommand{\dummy}[1]{}
\providecommand{\bysame}{\leavevmode\hbox to3em{\hrulefill}\thinspace}
\providecommand{\MR}{\relax\ifhmode\unskip\space\fi MR }
\providecommand{\MRhref}[2]{%
  \href{http://www.ams.org/mathscinet-getitem?mr=#1}{#2}
}
\providecommand{\href}[2]{#2}

\end{document}